\documentclass[12pt]{article}
\usepackage{amsfonts}
\usepackage{amssymb}
\usepackage{amscd, color}
\usepackage{hyperref}
\usepackage{amsmath}

\newcommand{\pref}[1]{\ref{#1} page \pageref{#1} }
\newcommand{\fref}[1]{(\ref{#1}) }

\def\boxit#1#2{\setbox1=\hbox{\kern#1{#2}\kern#1}%
\dimen1=\ht1 \advance\dimen1 by #1 \dimen2=\dp1 \advance\dimen2 by #1
\setbox1=\hbox{\vrule height\dimen1 depth\dimen2\box1\vrule}%
\setbox1=\vbox{\hrule\box1\hrule}%
\advance\dimen1 by .4pt \ht1=\dimen1
\advance\dimen2 by .4pt \dp1=\dimen2 \box1\relax}

\def\bo#1{\boxit{1pt}{$#1$}}

\newcommand{\BEQ}{\begin{equation}}
\newcommand{\EEQ}{\end{equation}}
\newtheorem{Th}{Theorem}

\newtheorem{Lem}{Lemma}
\newtheorem{Prop}{Proposition}
\newtheorem{Cor}{Corollary}
\newtheorem{Conj}{Conjecture}
\newtheorem{Ex}{Example}
\newtheorem{CEx}{Counter Example}

\newtheorem{Def}{Definition}
\newtheorem{Rem}{Remark}
\newtheorem{Quest}{Question}
\newtheorem{Observ}{Observation}
\def\nn{\nonumber}
\def\one#1{#1^{\raise5pt\hbox{$\scriptstyle\!\!\!\!1$}}\,{}}
\def\two#1{#1^{\raise5pt\hbox{$\scriptstyle\!\!\!\!2$}}\,{}}
\def\oneM#1{#1^{\raise5pt\hbox{$\scriptstyle\!\!\!\!m$}}\,{}}
\def\oneK#1{#1^{\raise5pt\hbox{$\scriptstyle\!\!\!\!k$}}\,{}}

\def\bea{\begin{eqnarray}}
\def\eea{\end{eqnarray}}

\def\O{{\cal O}}

\newcommand{\CC}{{{\mathbb{C}}}}

\def\L{\Lambda}

\def\esf{\sigma} 
\def\csf{S} 
\def\psf{T} 

\def\esft{\sigma(t)} 
\def\csft{S(t)} 
\def\psft{T(t)} 

\def\l{\lambda}

\def\p{\partial_z}
\def\ttpsi{\widetilde{\tilde \psi}}

\def\s{sgn}

\def\tpsi{\tilde \psi}

\def\GT{Manin}
\def\MMs{Manin\ matrices}
\def\MM{Manin\ matrix}
\def\BX{$\Box$}
\def\PRF{ {\bf Proof.} }
\def\SPRF{ {\bf Sketch of the Proof.} }
\def\GIOPS{ D. Gurevich, A. Isaev, O. Ogievetsky, P. Pyatov, P. Saponov papers \cite{GIOPS95}-\cite{GIOPS05}}

\def\PP{ P } 
\newcommand\kol{{{$\ \mathcal{K}$}}}  
\def\RRR{{{\mathcal{K}}}} 
\def\BB{Y} 
\newtheorem{Not}{Notation}
\def\Lcal{{\cal L}}

\begin{document}


\begin{titlepage}
\hfill ITEP-TH-43/08

\vskip 1.5cm

\centerline{\Large Algebraic properties of \MMs~ 1.}

\vskip 1.5cm
 \centerline{A. Chervov \footnote{E-mail:
chervov@itep.ru}
G. Falqui  \footnote{E-mail:
gregorio.falqui@unimib.it}
V. Rubtsov  \footnote{E-mail: volodya@tonton.univ-angers.fr}
}

\centerline{${}^1$ \sf Institute for Theoretical and Experimental Physics, Moscow - Russia}

\centerline{${}^2$\sf Universit\'a di Milano - Bicocca, Milano - Italy }

\centerline{${}^3$\sf Universit\'e D'Angers, Angers, France }

\vskip 1.0cm
\centerline{\large \bf  Abstract}
\vskip 1cm


We study a class of matrices with noncommutative entries, which were first
considered by Yu. I. Manin in 1988 in relation with quantum group theory.
They are defined as ``noncommutative endomorphisms'' of a
polynomial algebra. More explicitly their defining conditions read:
1) elements in the same column commute; 2)
commutators of the cross terms are equal: $[M_{ij} ,M_{kl}] = [M_{kj} ,M_{il}]$
(e.g. $[M_{11},M_{22}] = [M_{21},M_{12}]$). The basic claim is that
despite noncommutativity many theorems of linear algebra hold true
for \MMs~ in a form identical to that of the commutative case.
Moreover in some examples the converse is also true, that is,
\MMs~ are the most general class of matrices
such that linear algebra holds
true for them.
The present paper gives a complete list and detailed proofs of algebraic properties of \MMs\  known up to the moment; many of them are new.
In particular we
present the formulation of \MMs~ in terms of matrix (Leningrad) notations;
provide  complete proofs that an inverse to a \MM~ is again a \MM~
and for the Schur formula for the determinant of a block matrix;
we generalize  the noncommutative Cauchy-Binet
formulas discovered recently
\href{http://arxiv.org/abs/0809.3516}
{[arXiv:0809.3516]}, which includes the classical Capelli and related identities.
We also discuss many other properties, such as
the Cramer formula for the inverse matrix,
the Cayley-Hamilton theorem, Newton and MacMahon-Wronski
 identities, Pl\"{u}cker relations,
Sylvester's theorem,
the Lagrange-Desnanot-Lewis Caroll formula,
the Weinstein-Aronszajn formula,
some multiplicativity properties for the determinant,
relations with quasideterminants,
calculation of the determinant via Gauss decomposition,
conjugation to the second normal (Frobenius) form,
and so on and so forth. Finally  several examples and open question are discussed.
We refer to \cite{CF07,RST} for some applications in the realm of quantum integrable systems.

\vskip 1.0cm
\end{titlepage}

\tableofcontents


\section{Introduction}
It is well-known that  matrices with generically noncommutative
elements do not admit a natural construction of the determinant with values in a ground ring
and basic theorems of the linear algebra fail to hold true. On the other
hand, matrices with noncommutative entries play a basic role in the
theory of quantum integrability
(see, e.g., \cite{F79}),
in Manin's theory of "noncommutative symmetries" \cite{Manin}, and
so on and so forth.
 Further we prove that
many {\em results of commutative linear algebra} can be applied with
minor modifications in the case of "\MMs". 

We will consider the simplest case of those considered by Manin,
namely - in the present paper - we will restrict ourselves to the
case of commutators, and not of (super)-$q$-commutators, etc.
Let us mention that \MMs~ are defined, roughly speaking, by 
{\bf half} of the relations
of the corresponding quantum group $Fun_q(GL(n))$ and taking $q=1$ (see section \ref{QG} page
\pageref{QG}).
\begin{Def}\label{D1}
Let $M$ be an $n\times n'$ matrix with elements $M_{ij}$ in
(not necessarily commutative) ring $\mathcal{R}$. We will call $M$ a \MM\ if the
following two conditions hold:
\begin{enumerate}
\item Elements in the same column commute between themselves.
\item Commutators of cross terms of $2\times 2$ submatrices of $M$ are equal:
\begin{equation}\label{Cross-term-propert} [M_{ij}, M_{kl}]=[M_{kj}, M_{il}],\>
\forall\, i,j,k,l \text{ e. g. } [M_{11}, M_{22}]=[M_{21}, M_{12}].
\end{equation}
\end{enumerate}
\end{Def}
A more intrinsic definition of \MMs\ via coaction on polynomial and
Grassmann algebras will be recalled in Proposition \pref{Coact-pr}
below. (Roughly speaking variables $\tilde x_i = \sum_j M_{ij}x_j$ commute
among themselves if and only if $M$ is a \MM, where
$x_j$ are commuting variables, also commuting with elements of $M$).

In the previous paper \cite{CF07} we have shown that \MMs~ have various applications
in quantum integrability and outlined some of their basic properties.
This paper is devoted solely to algebraic properties providing a complete
amount of facts known up to the moment. Quite probably the properties
established here can be transferred to some other classes of matrices
with noncommutative entries,
for example, to super-q-\MMs~  (see \cite{CFRS})
and quantum Lax matrices of most of the integrable systems.
Such questions seems to be quite important for quantum integrability,
quantum and Lie groups,  as well as in the
geometric Langlands correspondence theory \cite{CT06-1, CF07}.
But before studying these complicated issues, it seems to be worth to understand
the simplest  case in depth, this is one of the main motivations for us.
The other one is that many statements are so simple and natural extension of the classical results,
that can be of some interest just out of curiosity or pedagogical reasons
for wide range of mathematicians.

\subsection{Results and organization of the paper  }

The main aim of the  paper is to argue the following claim:
{\bf linear algebra statements hold true for \MMs~
in a form { identical} to the commutative case.}
Let us give a list of main properties discussed below.
Some of these results are new, some - can be found
in a previous literature:
Yu. Manin has defined  the determinant and has proven a Cramer's inversion rule, Laplace formulas, as well as Pl\"{u}cker
identities; in \href{http://arxiv.org/abs/math.QA/0303319}{(S. Garoufalidis, T. Le, D. Zeilberger)}
\cite{GLZ}
the MacMahon-Wronski formula
was proved; \href{http://arxiv.org/abs/math.CO/0703203}{(M. Konvalinka)}
\cite{Konvalinka07-1, Konvalinka07-2}
found the Sylvester's identity and the Jacobi ratio's theorem, along
with partial results on an inverse matrix and block
matrices. \footnote{These authors actually considered more general classes of matrices
}
 Some other results were announced in \cite{CF07}, where
applications to integrable systems, quantum and Lie groups can be found.

\begin{itemize}
\item Section \ref{Det-ss}: Determinant can be defined in the standard way
and it satisfies standard properties e.g. it is completely antisymmetric function of columns and rows:
\bea det
M=det^{column} M=\sum_{\sigma\in S_n} (-1)^\sigma
\prod^{\curvearrowright}_{i=1,...,n} M_{\sigma(i),i}.
\eea

\item Proposition \pref{DetMultPr}: let $M$ be a \MM~ and $N$ satisfies: $\forall i,j,k,l:~[M_{ij}, N_{kl}]=0$:
\bea
det(M N)=det(M) det(N).
\eea
Let $N$ be additionally a \MM, then $MN$ and $M\pm N$ are \MMs.

Moreover in case $[M_{ij}, N_{kl}]\ne 0$, but obeys certain conditions we prove (theorem
\pref{ThGrasCSS}):
\bea
det^{col}(M\BB + Q~diag(n-1,n-2,...,1,0) )
= det^{col}(M) det^{col}(\BB).
\eea
This generalizes  \cite{CSS08} and the classical \href{http://gdz.sub.uni-goettingen.de/index.php?id=11&no_cache=1&IDDOC=26783&IDDOC=26783&branch=&L=1}
{Capelli}
 identity \cite{Ca87}.

\item
Section \ref{Cramer-s}: Cramer's rule:
\bea
\mbox{}  M^{-1} \mbox{ is a \MM~ and} &&
M^{-1}_{ij}=(-1)^{i+j} det(M)^{-1}det(\widehat{M_{ij}}).
~~~~~~~~~~~~~~~~~~~~~~~~~~~~~~~~~~~~~
\eea
Here as usually
$\widehat{M}_{l k}$ is the $(n-1)\times (n-1)$ submatrix of $M$ obtained
removing the l-th row and the k-th column.
\item
Section \ref{CH-ss}: the Cayley-Hamilton theorem: $det(t-M)|_{t=M}=0$.
\item
Section \ref{Schur-ss}: the formula for the determinant of block matrices:
\bea
 det\left(\begin{array}{cc}
A & B\\
C & D \\
\end{array}\right)=det(A)det(D-CA^{-1}B)=det(D)det(A-BD^{-1}C).
\eea
Also, we show that  $ D-CA^{-1}B$, $A-BD^{-1}C$ are \MMs.
This is equivalent to the so-called Jacobi ratio theorem, stating that
any minor of $M^{-1}$ equals,
up to a sign, to the product of  $(det^{col}M)^{-1}$
and the corresponding complementary minor of the transpose of $M$.
\item
Section \ref{NMWss}: the Newton and MacMahon-Wronski
 identities between $TrM^k$, coefficients of $det(1-tM)$ and $Tr(S^k M)$.
Denote by $\esft, \csft , \psft$ the following generating functions:
\bea
& \esft =det(1-t M), ~~~  \csft = \sum_{k=0,...,\infty } t^k Tr S^{k}M, ~~~
\psft   =   Tr \frac{M}{1-tM}, &\\
 &\mbox{{\bf Then:}}~~~~ 1= \esft \csft, ~~~ -\partial_t \esft = \esft \psft, ~~~  \partial_t \csft = \psft \csft
.
\eea
\item
Other  facts are also discussed:
 Pl\"{u}cker relations (section \ref{Pluck-ss}),
Sylvester's theorem (section \ref{Sylv-ss}),
Lagrange-Desnanot-Lewis Caroll formula (section \ref{DJ-ss}),
Weinstein-Aronszajn formula (section \ref{WA-ss}),
calculation of the determinant via Gauss decomposition (section \ref{Gaus-ss}),
conjugation to the second normal (Frobenius) form (section \ref{CH-ss}),
some multiplicativity properties for the determinant (proposition \pref{block-mult-prop1}),
etc.
\item
Section \ref{PLss}: matrix (Leningrad) notations form of the definition:
\bea
& &\mbox{$M$ is a \MM} ~~ \Leftrightarrow\\
& &\Leftrightarrow [ M\otimes 1, 1 \otimes M ] = \PP [ M\otimes 1, 1 \otimes M ] \Leftrightarrow \\
& &\Leftrightarrow \frac{(1-\PP)}{2} (M\otimes 1)~ (1 \otimes M) \frac{(1-\PP)}{2}=
\frac{(1-\PP)}{2} (M\otimes 1)~ (1 \otimes M).
\eea
\item
No-go facts: $M^k$ is not a \MM;
elements $Tr(M)$, $det(M)$, etc.  are not central, moreover $[Tr(M),det(M)]\ne 0$ (section
\ref{nogosect});
 $det(e^M)\ne e^{Tr(M)}, det(1+M)\ne e^{Tr(ln(1+M))}$
(section \ref{nogoNewtsect}).
\end{itemize}

We also discuss relations with the quantum groups (section \ref{QG}) and mention
some examples which are related to integrable systems, Lie algebras and quantum groups
(section \ref{CF-ss}).

Organization of the paper: 
section \ref{DefSect} contains main definitions and  properties.
This material is crucial for what follows. The other sections
can be read in an arbitrary order. We tried to make the
exposition in each section independent of the others
at least at formulations of the main theorems and notations.
Though the proofs sometimes use the results from the previous sections.
A short section \ref{WarmUpSect} is
a kind if warm-up, which gives some simple examples
to get the reader interested and to demonstrate some
of  results of the paper on the simplest examples.
The content of each section can be seen from the table
of contents.


\subsection{Context, history and related works  }
\MMs~ appeared first in \cite{Manin},
see also \cite{Manin87,ManinBook91, Manin92, ManinDemid},
where some basic facts like the determinant, Cramer's rule, Pl\"{u}cker identities etc.
were established. The lectures \cite{Manin} is the main  source on the subject.
 Actually Manin's construction  defines "noncommutative endomorphisms"
 of an arbitrary ring (and in principle of any algebraic structure).
Here we restrict ourselves with the simplest case of the commutative polynomial ring $\CC[x_1,...,x_n]$,
its "noncommutative endomorphisms" will be called \MMs.
Some linear algebra facts were also established for a class of "good" rings
(we hope to develop this in future).
The main attention and application  of the original works were on quantum groups,
which are defined by "doubling" the set of the relations.
The literature on "quantum matrices" (N. Reshetikhin, L. Takhtadzhyan, L. Faddeev)
\cite{FRT}  is enormous - let
us only mention \cite{KL94}, \cite{DL03} and especially \GIOPS,
where related linear algebra facts has been established
for various quantum matrices.

Concerning "not-doubled" case let us mention
 papers
\href{http://arxiv.org/abs/math.OA/9807091}{(S. Wang)}
\cite{Wang98},
\href{http://arxiv.org/abs/math/0612724}{(T. Banica, J. Bichon, B. Collins)}
\cite{BBC06}.
 They
investigate Manin's construction applied to finite-dimensional
algebras for example to $\CC^n$ (quantum permutation group).
Such algebras appears to be $C^*$-algebras and are related
to various questions in operator algebras.
Linear algebra of such matrices is not known, (may be it does not exist).

The simplest case which is considered here
was somehow forgotten for many years after Manin's work (see however \cite{RRT02, RRT05}).
The situation changed after \href{http://arxiv.org/abs/math.QA/0303319}{(S. Garoufalidis, T. Le, D. Zeilberger)}
\cite{GLZ}, who discovered MacMahon-Wronski relations
for (q)-\MMs.
This result was followed by a flow of papers
\cite{EtingofPak06,HaiLorenz06,KonvalinkaPak06,FoataHan06} etc.;
let us in particular mention \href{http://arxiv.org/abs/math.CO/0703203}{(M. Konvalinka)}
\cite{Konvalinka07-1, Konvalinka07-2}
which contains Sylvester's identity and  Jacobi's ratio theorem, along
with partial results on an inverse matrix and block
matrices.

We came to this story from the other direction.
In \cite{CF07}  we observed that some examples of quantum Lax
matrices in quantum integrability are exactly examples of \MMs.
Moreover linear algebra of \MMs~ appears to have various applications in quantum
integrability, quantum and Lie group theories.
Numerous other  quantum integrable systems provide various examples of matrices
with noncommutative elements - quantum Lax matrices.
This rises the question how to define the proper determinant and
to develop linear algebra for such matrices. And more generally
proper determinant  exists or not? If yes, how to develop the linear algebra?
Such questions seems to be quite important for quantum integrability
they are related to a "quantum spectral curve" which promises
to be the key concept in the theory \cite{CT06-1}, \cite{CM, CFRy, CFRS,RST}.
Manin's approach applied for more general rings provides classes of matrices
 where such questions can be possibly  resolved, however many examples from integrable
systems do not fit in this approach.

From the more general point of view
we deal with the question of the linear algebra of matrices with noncommutative entries.
It should be remarked that the first appearance of a determinant for matrices with  noncommutative entries
goes back to A. Cayley. He was the first person who had applied the notion of what we call \emph{a column determinant}
in the non-commutative setting. (We are thankful to V. Retakh for this remark).
Let us  mention the initial significative difference between our situation and the work
\href{http://arxiv.org/abs/math.QA/0208146}{I. Gelfand, S. Gelfand, V. Retakh, R. Wilson}
\cite{GGRW02},
where  {\em generic } matrices with noncommutative entries are considered.
There is no natural definition of the determinant ($n^2$ of "quasi-determinants" instead)
in the "general non-commutative case" and  their analogs of the
linear algebra propositions are sometimes quite different from the commutative ones.
Nevertheless results of loc.cit. can be fruitfully applied to some questions here.
Our approach is also different from the classical theory of J. Dieudonn\'e
\cite{Died} (see also \cite{Adj93}),
since in this theory the determinant is an element of the
$\mathcal{K}^{*}/[\mathcal{K}^{*},\mathcal{K}^{*}]$,
where \kol~ is basic ring, while for \MMs~ the determinant is an element of the ring \kol~ itself.

We provide more detailed bibliographic notes in the text but we would like to add,
as a
{\em general disclaimer},
that our bibliographic notes are neither exhaustive nor historically ordered. We simply want to comment those
papers and books that are more strongly related to our work.

We refer to \href{http://books.google.com/books?id=lfnyeEv4O1wC}{
``The Theory of Determinants in the Historical Order of Development''}
\cite{MuirBook} for the early history of many results, which
generalization to the noncommutative case are discussed  below.

\subsection{Remarks}
In  \cite{GLZ, Konvalinka07-1,Konvalinka07-2}
the name "right quantum matrices"
was used,
in \cite{RRT02,RRT05, LT07} the names ``left'' and ``right quantum group'',
 in  \cite{CSS08} the name ``row-pseudo-commutative''.
We prefer to use the name "\MMs".

All the considerations below work for an arbitrary field of characteristic
not equal to 2, but we prefer to restrict ourselves with $\CC$.

In subsequent papers  \cite{CFRS}
we plan to generalize the constructions below
to the case of the Manin matrices related to the more general quadratic algebras
as well as there applications to quantum integrable
systems and some open problems.

{\bf Acknowledgments}
G.F. acknowledges support from the ESF
programme MISGAM, and the Marie Curie RTN ENIGMA.
The work of A.C. has been partially supported by the
Russian President Grant MK-5056.2007.1, grant of Support for the Scientific Schools NSh-3035.2008.2,
RFBR grant 08-02-00287a, the ANR grant GIMP (Geometry and Integrability in
Mathematics and Physics). He acknowledges support and
hospitality of Angers and Poitiers Universities.
The work of V.R. has been partially supported by the
grant of Support for the Scientific Schools NSh-3036.2008.2,
RFBR grant 06-02-17382, the ANR grant GIMP (Geometry and Integrability in
Mathematics and Physics) and INFN-RFBR project "Einstein". He acknowledges support and
hospitality of SISSA (Trieste).
The authors are grateful
to D. Talalaev, A. Molev, A. Smirnov,  A. Silantiev, V. Retakh and D. Gurevich  for multiple stimulating discussions,
to P. Pyatov for sharing with us his unpublished results,
pointing out to the paper \cite{GLZ}
 and for multiple  stimulating discussions.
To Yu. Manin for his interest in this work and stimulating discussions.

\section{
\label{WarmUpSect}
Warm-up $2\times 2$ examples:  \MMs~ everywhere}

Here we present some examples of the appearance of the Manin
property  in various very simple and natural questions concerning
$2\times 2$ matrices. The general idea is the following: we consider
well-known facts of linear algebra and look how to relax the
commutativity assumption for matrix elements such that the results
will be still true. The answer is: {\bf if and only if} $M$ is a
\MM.

This section is a kind of warm-up, we hope to get the reader
interested in the subsequent material and to demonstrate some
results in the simplest examples. The expert reader may wish to skip
this section.

Let us consider a $2\times 2$ matrix $M$: \bea M=\left(
\begin{array}{cc}
a & b \\
c & d
\end{array}\right).
\eea From Definition (\ref{D1}) $M$ is a "\MM" if the following
commutation relations hold true:
\begin{itemize}
\item column-commutativity: $[a,c]=0$,  $[b,d]=0$;
\item equality of commutators of the cross-term: $[a,d]=[c,b]$.
\end{itemize}

\vskip 1cm
\noindent The fact below can be considered as Manin's original idea about the subject.

{\Observ Coaction on a plane.} Consider the polynomial ring
$\CC[x_1,x_2]$, and  assume that the matrix elements $a,b,c,d$
commute with $x_1,x_2$. Define $\tilde x_1, \tilde x_2$ by
$\left(\begin{array}{c}
\tilde x_1 \\
\tilde x_2
\end{array}\right) =
\left(\begin{array}{cc}
a & b \\
c & d
\end{array}\right)
\left(\begin{array}{cc}
x_1 \\
x_2
\end{array}\right)
$.
Then $\tilde x_1, \tilde x_2$ commute among themselves {\bf iff} $M$ is a \MM.

\PRF
\bea
[ \tilde x_1, \tilde x_2]=
[a x_1 + b x_2, c x_1 + d x_2]=
[a, c] x_1^2 + [b,d] x_2^2 +  ([a,c] + [b,d] ) x_1 x_2. \hfill \Box
\eea
Similar fact holds true for Grassman variables (see proposition \pref{Coact-pr} below).

{\Observ \label{CramRulObser} Cramer rule. }
The inverse matrix is given by the standard formula\\
$
M^{-1}=
\frac{1}{ad-cb}  \left(\begin{array}{cc}
d & -b \\
-c & a
\end{array}\right)
$ {\bf iff} $M$ is a \MM.

\PRF
\bea  \left(\begin{array}{cc}
d & -b \\
-c & a
\end{array}\right)
\left(\begin{array}{cc}
a & b \\
c & d
\end{array}\right)
=
\left(\begin{array}{cc}
da-bc & db-bd \\
-ca+ac & -cb+ad
\end{array}\right) = \\
 \mbox { {\bf iff} $M$ is a \MM } =
\left(\begin{array}{cc}
ad-cb & 0 \\
0 & ad-cb
\end{array}\right). \hfill \Box
\eea

\vskip 1cm

{\Observ Cayley-Hamilton.} The equality
 $M^2-(a+d)M+(ad-cb)1_{2\times 2} =0$ holds {\bf iff} $M$ is a \MM.

\bea
M^2-(a+d)M+(ad-cb)1_{2\times 2}=\nn\\=
\left(\begin{array}{cc}
a^2 + bc & ab+bd  \\
ca+dc & cb+d^2
\end{array}\right)
-
\left(\begin{array}{cc}
a^2+da & ab+db\\
ac+dc  & ad +d^2
\end{array}\right)
+
\left(\begin{array}{cc}
ad-cb & 0\\
0  & ad -cb
\end{array}\right)
=\nn\\=
\left(\begin{array}{cc}
 (bc-da) +(ad-cb) & bd -db \\
ca-ac & 0
\end{array}\right)
=
\left(\begin{array}{cc}
 [a,d]-[c,b] & [b,d] \\
 ~ [c,a] & 0
\end{array}\right).
\eea This vanishes {\bf iff} $M$ is a \MM.
\BX ~~~~
Let us mention that similar facts can be seen for the Newton
identities, but not in such a strict form (see example
\pref{Newt22}).
\vskip 1cm


{\Observ \label{DetMulCase22} Multiplicativity of Determinants (Binet Theorem).}\\
$
det^{column}(M N)=det^{column}(M)det
(N) $
 holds true for all  $\CC$-valued matrices $N$
 {\bf iff} $M$ is a \MM.
\bea
det^{col}(M N)-det^{col}(M)det
(N)
=[M_{11},M_{21}]N_{11}N_{12} + [M_{12},M_{22}]N_{21}N_{22}
+([M_{11},M_{22}]-[M_{21},M_{12}])N_{21}N_{12}
\eea

\noindent Here and in the sequel, $det^{column}(X)$ or shortly $det^{col}(X)$ is:
$X_{11}X_{22}-X_{21}X_{12}$, i.e. elements from the
first column stand first in each term.

\subsection{Further properties in $2\times 2$ case }
Let us also present some other properties of $2\times 2$ \MMs.

It is well-known that in commutative case a matrix can be conjugated
to the so-called Frobenius normal form. Let us show that the same
is possible for \MMs, see also page \pageref{FrobCor}.

{\Observ Frobenius
form of a matrix. }
\bea \label{FrobFormulIntr}
\left(
\begin{array}{cccc}
1& 0 \\
a & b
\end{array}\right)
\left(\begin{array}{cc}
a & b \\
c & d
\end{array}\right)
\left(
\begin{array}{cccc}
1& 0 \\
a & b
\end{array}\right)^{-1} =
\left(
\begin{array}{cccc}
0& 1 \\
-(ad-cb) & a+d
\end{array}\right),
\eea

\centerline{ {\bf iff} $[d,a]=[b,c]$ and $db=bd$. }

So we got two of three Manin's relations, to get the last third
relation see example \pref{ExFrobForm2}.

Let us denote the matrix at the right hand side of  \fref{ExFrobForm2}
by $M_{Frob}$ and the first matrix at the left hand side by $D$.
To see that  \fref{ExFrobForm2} is true we just write the following.


\bea D~ M= \left(
\begin{array}{cccc}
1& 0 \\
a & b
\end{array}\right)
~ \left(
\begin{array}{cccc}
a & b \\
c & d
\end{array}\right)
= \left(
\begin{array}{cccc}
a & b \\
a^2+bc & ab+bd
\end{array}\right).
\eea

\bea
 M_{Frob} ~ D=
\left(
\begin{array}{cc}
0& 1 \\
-(ad-cb) & a+d
\end{array}\right)
~ \left(
\begin{array}{cc}
1& 0 \\
a & b
\end{array}\right)
= \left(
\begin{array}{cc}
a & b \\
-(ad-cb)+a^2+da & ab+db
\end{array}\right)
= \\   = \left(
\begin{array}{cc}
a & b \\
a^2+bc+([d,a]-[b,c]) & ab+bd+[d,b]
\end{array}\right).
\eea
{\Observ Inversion \label{ObeservInv} }. The two sided inverse of a Manin matrix $M$ is Manin, and $det(M^{-1})=(det\ M)^{-1}$.\\

See theorem \pref{Inverse-conj} and corollary \pref{DetInvCor}.

Let us briefly prove this fact.
From the Cramer's rule above, one knows the formula
for the left inverse, by assumption it is also right inverse. To
prove the theorem one only  needs  to write explicitly that the
right inverse is given by Cramer rule and the desired commutation
relations appear automatically. Explicitly, from Cramer's formula
(see Observation 2) we see that: \bea
\frac{1}{ad-cb} \left(\begin{array}{cc}
d & -b \\
-c & a
\end{array}\right)
\left(\begin{array}{cc}
a & b \\
c & d
\end{array}\right)=
\left(\begin{array}{cc}
1 & 0 \\
0 & 1
\end{array}\right).
\eea One knows that if both left and right inverses exist then
associativity guarantees that they coincide: $a^{-1}_l =a^{-1}_l  (a
a^{-1}_r)= (a^{-1}_l a) a^{-1}_r=  a^{-1}_r$. So assuming that the
right inverse to $A$ exists, and denoting $ad-cb\equiv \delta$ we
have: \bea \label{lfml11113} \left(\begin{array}{cc}
1 & 0 \\
0 & 1
\end{array}\right)
=
\left(\begin{array}{cc}
a & b \\
c & d
\end{array}\right)
(\delta^{-1}) \left(\begin{array}{cc}
d & -b \\
-c & a
\end{array}\right)
=
\left(\begin{array}{cc}
a (\delta)^{-1}d - b (\delta)^{-1} c & -a (\delta)^{-1}b+  b (\delta)^{-1} a \\
c (\delta)^{-1}d - d (\delta)^{-1} c & - c (\delta)^{-1}b + d
(\delta)^{-1} a
\end{array}\right).
\eea Let us multiply the identity  above by $\delta^{-1}$ on the
left. We have: \bea \label{lfml1111}
&&(\delta)^{-1} = (\delta)^{-1}a (\delta)^{-1}d - (\delta)^{-1}b (\delta)^{-1} c, \mbox{ Element (1,1)},\\
&& (\delta)^{-1}  =- (\delta)^{-1}c (\delta)^{-1}b +(\delta)^{-1} d
(\delta)^{-1} a, \mbox{ Element (2,2)}. \label{lfml1112} \eea So we
see that $(det(M))^{-1}$ equals to  $det(M^{-1})$ and the last does
not depend on the ordering of columns.

Moreover, equating (\ref{lfml1111}) to (\ref{lfml1112}) yields \bea
&&(\delta)^{-1}a (\delta)^{-1}d - (\delta)^{-1}b (\delta)^{-1} c =
- (\delta)^{-1}c (\delta)^{-1}b +(\delta)^{-1} d (\delta)^{-1} a, \\
&& \mbox{ hence:~ ~ } [ (\delta)^{-1}a, (\delta)^{-1} d ] =
[(\delta)^{-1}b, (\delta)^{-1}c]. \eea So the commutators of the
cross-terms of $M^{-1}$ are equal.

From the non-diagonal elements of  equality \ref{lfml11113}
multiplied on the left by $\delta^{-1}$ we have: \bea
 (\delta)^{-1}c (\delta)^{-1}d -  (\delta)^{-1}d (\delta)^{-1} c=0,  ~~~~
- (\delta)^{-1}a (\delta)^{-1}b+  (\delta)^{-1} b (\delta)^{-1} a=0, \mbox{hence}  \nn\\
~~ ~ [(\delta)^{-1}c ,(\delta)^{-1}d ]=0 ~~~~ [(\delta)^{-1} a,
(\delta)^{-1}b] = 0. ~~~~~ \eea So we also have the column
commutativity of the elements of $M^{-1}$. Hence the proposition is
proved in $2\times 2$ case. \BX

{\Observ A puzzle with $det(M)=1$.}
 Let $M$ be a $2\times 2$ \MM,
suppose that $det(M)$ is central element and it is invertible (for example $det(M)=1$).
 Then all elements of $M$ commute
  among themselves.

From the observation \ref{ObeservInv}  above one gets that \bea
M^{-1}=\frac{1}{ det(M)} \left(\begin{array}{cc}
d & -b \\
-c & a
\end{array}\right)
\eea is again a \MM. This gives the commutativity. \BX

This is quite a surprising  fact that: imposing only one
condition we ``kill'' the three
commutators: $[a,b], [a,d]=[c,b], [c,d]$.

In the paper we will consider $n\times n$ \MMs~ and
prove that these and other  properties of linear algebra  work for
them.


\section{
\MMs. \label{DefSect}
Definitions and elementary properties}

In this section we recall the definition of \MMs~ and give some
basic properties. The material is rather simple one, but it is
necessary for the sequel. First we will give an explicit definition
of \MMs (in terms of commutation relations), and  then we will
provide a more conceptual point of view which defines them by the
coaction property on the polynomial and Grassman algebras. (This is
the original point of view of Manin). We also explain the relation
of \MMs~ with quantum groups.
As it was shown by Yu. Manin there
exists natural definition of the determinant which satisfies most of
the properties of commutative determinants; this will be also
recalled below. 
The main reference for this part is Yu. Manin's book
\cite{Manin}, as well as \cite{Manin87,ManinBook91, Manin92,
ManinDemid}.

\subsection{Definition} 
{\Def \label{MM-def} Let us call a matrix $M$  with elements in an associative
  ring 
\kol\
a {\em "\MM"}~
if the properties below are satisfied:}
\begin{itemize}
\item elements which belong to the same column of $M$ commute among themselves
\item commutators of  cross terms are equal:
$\forall p,q,k,l~~[M_{pq}, M_{kl} ]= [M_{kq}, M_{pl} ]$, \\ e.g. $[M_{11}, M_{22} ]= [M_{21}, M_{12} ]$,
$[M_{11}, M_{2k} ]= [M_{21}, M_{1k} ]$.
\end{itemize}
{\Rem ~} The second condition for the case $q=l$  obviously implies
the first one. Nevertheless we deem it more convenient for the
reader to formulate it explicitly.

The conditions can be restated as:
 for each $2\times 2$ submatrix:
\bea
\label{Cross-term-propert222}
\left(\begin{array}{ccccc}
... & ...& ...&...&...\\
... & a &... & b& ... \\
... & ...& ...&...&...\\
... & c &... & d& ... \\
... & ...& ...&...&...
\end{array}\right)
\mbox{ ~~~~~~~~~  of $M$ it holds~~~~~} [a,d]=[c,b],~~~
[a,c]=0=[b,d]. \eea Obviously, by a "submatrix" we mean a matrix
obtained as an intersection of two rows (i.e. straight horizontal
lines, no decline) and two columns (i.e. straight vertical lines, no
decline).

{\Rem ~} These relations were written by Yu. Manin \cite{Manin} (see
chapter 6.1,  Formula 1, page 37). Implicitly they are contained in
(\href{http://www.numdam.org/item?id=AIF_1987__37_4_191_0} {Yu.
Manin}) \cite{Manin87} -- the last sentence on page 198 contains a
definition of the algebra
 $end(A)$ for an arbitrary quadratic Koszul algebra $A$. One can show (see the remarks
on the page 199 top) that $end(\CC[x_1,...,x_n])$ is the algebra generated by $M_{ij}$.

{\Rem ~} Actually, a matrix $M$ such that $M^{t}$ is a \MM~ satisfies the same good properties as \MMs~ do,
we will sometimes mention this case explicitly.

\subsubsection{Poisson version of \MMs \label{Pois-mm}}
{\Def An algebra \kol~ over $\CC$ is called a Poisson algebra, if it
is a commutative algebra, endowed with a bilinear antisymmetric
operation, denoted as  $\{ *, * \}:R\otimes R\to R$ and called a
Poisson bracket, such that  the operation satisfies the Leibniz and
the Jacobi identities (i.e. $\forall f,g,h \in R:$  $\{fg, h\}= f
\{g, h\}+ \{f, h\}g$, $\{f ,\{ g, h\} \}+ \{ g ,\{ h, f \} \} +  \{h
,\{ f, g \} \}=0$). }

{\Def We call "Poisson-Manin" a matrix $M$ with elements in
the Poisson algebra \kol, such that $\{ M_{ij}, M_{kl} \} = \{
M_{kj}, M_{il} \} $}.

We briefly discuss Poisson-Manin matrices in section
\pref{PoissonManinSec}.

\subsection{Characterization via coaction. Manin's construction \label{Coact-ss}} 
Here we recall Manin's original definition. It provides a conceptual
approach to \MMs. Let us mention that the construction below is a
specialization of Manin's general considerations. (See
(\href{http://www.numdam.org/item?id=AIF_1987__37_4_191_0} {Yu.
Manin}) \cite{Manin87},
\cite{Manin,ManinBook91}.)

{\Prop \label{Coact-pr} Coaction. Consider a rectangular $n\times
m$-matrix $M$, the polynomial algebra $\CC[x_1,...,x_m]$ and the
Grassman algebra $\CC[\psi_1,...,\psi_n]$ (i.e. $\psi_i^2=0, \psi_i
\psi_j= -\psi_j \psi_i$); let $x_i$ and $\psi_i$ commute with
$M_{pq}$: $\forall i,p,q:~~[x_i, M_{pq}]=0$, $[\psi_i, M_{pq}]=0$.
Consider new variables $\tilde x_i$, $\tilde \psi_i$: \bea
\left(\begin{array}{c}
\tilde x_1 \\
...  \\
\tilde x_n
\end{array}\right)
=
\left(\begin{array}{ccc}
 M_{11} & ... &  M_{1m} \\
... & ... & ... \\
 M_{n1} & ... &  M_{nm}
\end{array}\right)
\left(\begin{array}{c}
 x_1 \\
... \\
 x_m
\end{array}\right),
~~~~
(\tilde \psi_1, ... ,  \tilde \psi_m) =
(\psi_1, ... ,   \psi_n)
\left(\begin{array}{ccc}
 M_{11} & ... &  M_{1m} \\
...  & ... & ... \\
 M_{n1} & ... &  M_{nm}
\end{array}\right).
\eea
Then the following three conditions are equivalent:
\begin{itemize}
\item $M$ is a \MM 
\item the variables $\tilde x_i$ commute among themselves: $[\tilde x_i,\tilde x_j]=0$
\item the variables $\tilde \psi_i$ anticommute among themselves:
$\tilde \psi_i\tilde \psi_j+\tilde \psi_j \tilde \psi_i =0$.
\end{itemize}
} {\Rem ~} The conditions $\tpsi_i^2=0$ are equivalent to column
commutativity, and  $\tpsi_i \tpsi_j=-\tpsi_j \tpsi_i$, $i<j$, to
the cross term relations.

\subsection{q-analogs and RTT=TTR quantum group matrices  \label{QG} }
One can define q-analogs of \MMs~ and characterize their relation to
quantum group theory. Actually q-\MMs~ are defined by half of the
relations of the corresponding quantum group
$Fun_q(GL_n)$\footnote{More  precisely we should write
$Fun_q(Mat_n)$, since we do not localize the q-determinant.}
(\cite{FRT}). The remaining half consists of relations insuring that
also $M^t$ is a q-\MM. 

{\Def \label{D1a22}
Let us call an $n\times n'$ matrix $M$ by {\em q-\MM},
if the following conditions hold true.
For any $2\times 2$  submatrix
$(M_{ij,kl})$, consisting of rows $i$ and $k$,
and columns $j$ and $l$ (where $1 \leq i < k \leq n$,
and $1 \leq j < l \leq n'$):
\bea
\left(\begin{array}{ccccc}
... & ...& ...&...&...\\
... & M_{ij} &... & M_{il} & ... \\
... & ...& ...&...&...\\
... & M_{kj} &... & M_{kl}& ... \\
... & ...& ...&...&...
\end{array}\right)
\equiv
\left(\begin{array}{ccccc}
... & ...& ...&...&...\\
... & a &... & b& ... \\
... & ...& ...&...&...\\
... & c &... & d& ... \\
... & ...& ...&...&...
\end{array}\right)
\eea
the following commutation relations hold:
\begin{eqnarray}
\label{eq.commutation1}
ca &= & qac, \quad \text{($q$-commutation of the entries in a column)} \\
\label{eq.commutation2}
 db &=& qbd, \quad \text{($q$-commutation of the entries in a column)} \\
\label{eq.commutation3}
ad - da &=& +q^{-1}cb-q bc, \qquad \text{(cross commutation relation)}.
\end{eqnarray}
}
In terms of $M_{ij}$ this reads ($i<k, j<l$) : \bea M_{kj} M_{ij}
= q M_{ij} M_{kj}, ~~~~ M_{ij}M_{kl} - M_{kl}  M_{ij}=q^{-1}
M_{kj}M_{il}-q M_{il}  M_{kj}. \eea

For $q=1$ this definition reduces to the definition \pref{MM-def} of
\MMs.

{\Def \label{QGdef} An $n \times n$ matrix $T$ belongs to the
quantum group $Fun_q(GL_n)$ if the following conditions hold true.
For any $2\times 2$  submatrix $(T_{ij,kl})$, consisting of rows $i$
and $k$, and columns $j$ and $l$ (where $1 \leq i < k \leq n$, and
$1 \leq j < l \leq n$): \bea \left(\begin{array}{ccccc}
... & ...& ...&...&...\\
... & T_{ij} &... & T_{il} & ... \\
... & ...& ...&...&...\\
... & T_{kj} &... & T_{kl}& ... \\
... & ...& ...&...&...
\end{array}\right)
\equiv \left(\begin{array}{ccccc}
... & ...& ...&...&...\\
... & a &... & b& ... \\
... & ...& ...&...&...\\
... & c &... & d& ... \\
... & ...& ...&...&...
\end{array}\right)
\eea
the following commutation relations hold:
\begin{eqnarray}
ca &= & qac, \quad \text{($q$-commutation of the entries in a column)} \\
 db &=& qbd, \quad \text{($q$-commutation of the entries in a column)} \\
ba &= & qab, \quad \text{($q$-commutation of the entries in a row)} \\
 dc &=& qcd, \quad \text{($q$-commutation of the entries in a row)} \\
ad - da &=& +q^{-1}cb-q bc, \qquad \text{(cross commutation relation 1)} \\
bc  &=& cb, \qquad \text{(cross commutation relation 2)}.
\end{eqnarray}
}

As quantum groups are usually defined within the so-called matrix
(Leningrad) formalism, let us briefly recall it. (We will further
discuss this issue in section \pref{MatrSect}). {\Lem The
commutation relations for quantum group matrices can be described in
matrix (Leningrad) notations as follows: \bea R (T\otimes 1)
(1\otimes T) =(1\otimes T) (T\otimes 1) R, \eea where R-matrix can
be given, for example, by the formula: \bea R= q^{-1}
\sum_{i=1,..,n} E_{ii}\otimes E_{ii} +\sum_{i, j=1,..,n; i\ne j}
E_{ii}\otimes E_{jj} +(q^{-1}-q) \sum_{i,j=1,..,n;i>j} E_{ij}\otimes
E_{ji}, \eea where $E_{ij}$ are standard matrix units - zeroes
everywhere except 1 in the intersection of the $i$--th row with the
$j$--th column.} For example in $2 \times 2$ case the R-matrix is:
\bea R= \left(\begin{array}{ccccc}
q^{-1}  & 0        & 0 & 0 \\
0       & 1        & 0 & 0 \\
0       & q^{-1}-q & 1 & 0 \\
0       & 0        & 0 & q^{-1}
\end{array}\right).
\eea

{\Rem ~} This R-matrix differs by the change $q$ to $q^{-1}$ from
the one in \cite{FRT} formula 1.5, page 185.

The relation between q-\MMs~ and quantum groups consists in the
following simple  proposition: {\Prop A matrix $T$ is a matrix in
the quantum group $Fun_q(GL_n)$ if and only if  ~ $T$ and
simultaneously the transpose matrix $T^t$ are q-\MMs.}

So one sees that \MMs\ can be seen as characterized by a "half" of
the conditions that characterize the corresponding  quantum matrix
group.

q-\MMs~ can be characterized by the coaction on a q-polynomial and a
q-Grassman algebra in the same way as in $q=1$ case. Here is the
analogue of proposition \pref{Coact-pr}.
{\Prop\label{Qproposizione}
Consider a rectangular $n\times m$-matrix $M$, the $q$-polynomial
algebra $\CC[x_1,...,x_m]$, where $\forall i<j:~~x_i x_j=q^{-1} x_j
x_i$,
 i.e. $\forall i,j:~~x_i x_j=q^{\s(i-j)} x_j x_i$
 and the $q$-Grassmann
algebra $\CC[\psi_1,...,\psi_n]$ (i.e. $\psi_i^2=0, \psi_i \psi_j=
-q\psi_j \psi_i$, for $i<j$; i.e. $\forall i,j:~~\psi_i \psi_j=-q^{-\s(i-j)} \psi_j \psi_i$);
suppose $x_i$ and $\psi_i$ commute with the matrix
elements $M_{pq}$.
Consider the variables $\tilde x_i$, $\tilde
\psi_i$ defined by: \begin{eqnarray} \left(\begin{array}{c}
\tilde x_1 \\
... \\
\tilde x_n
\end{array}\right)
=
\left(\begin{array}{ccc}
 M_{11} & ... &  M_{1m} \\
... & ... & ... \\
 M_{n1} & ... &  M_{nm}
\end{array}\right)
\left(\begin{array}{c}
 x_1 \\
... \\
 x_m
\end{array}\right),
~~~~
(\tilde \psi_1, ... ,  \tilde \psi_m) =
(\psi_1, ... ,   \psi_n)
\left(\begin{array}{ccc}
 M_{11} & ... &  M_{1m} \\
... & ... & ... \\
 M_{n1} & ... &  M_{nm}
\end{array}\right),
\end{eqnarray} that is the new variables are obtained via left action (in the
polynomial case) and right action (in the Grassmann case) of $M$ on
the old ones. Then the following three conditions are equivalent:
\begin{itemize}
\item The matrix $M$ is a q-\MM.
\item The variables
$\tilde x_i$ $q$-commute among themselves: $\forall i<j:~~\tilde x_i \tilde x_j=q^{-1} \tilde x_j \tilde x_i$,
 i.e. $\forall i,j:~~\tilde x_i \tilde x_j=q^{\s(i-j)} \tilde x_j \tilde x_i$.
\item The variables $\tilde \psi_i$ $q$-anticommute among themselves:
$\tpsi_i^2=0, \tpsi_i \tpsi_j=
-q\tpsi_j \tpsi_i$, for $i<j$; i.e. $\forall i,j:~~\tpsi_i \tpsi_j=-q^{-\s(i-j)} \tpsi_j \tpsi_i$
.
\end{itemize}}

\label{rem3} The conditions $\tpsi_i^2=0$ are equivalent
to the relations~\eqref{eq.commutation1}, \eqref{eq.commutation2},
and the conditions $\tpsi_i \tpsi_j=-q\tpsi_j \tpsi_i$, $i<j$,
are equivalent to the relations~\eqref{eq.commutation3}.

We plan to discuss q-\MMs~ and
some of their applications in the theory of integrability in a
subsequent publications (see \cite{CFRS}).

\subsection{The determinant \label{Det-ss}} 
Here we recall a definition of the determinant following  Manin's
ideas. It is well-known that, for generic matrices over a
noncommutative ring (possibly, algebra), one cannot develop a full theory of
determinants with values in the same ring. However, for some
specific matrices there may exist a "good" notion of determinant.
In particular 
for \MMs~ one can define the determinant just taking
the column expansion as a definition. Despite its simplicity
such a definition is actually a good one. It satisfies almost all the properties of the
determinants in the commutative case and is consistent
with the other concepts of noncommutative determinants
(quasideterminants of
\href{http://arxiv.org/abs/q-alg/9705026}{I. Gelfand, V. Retakh} \cite{GR97}
(see section \pref{SectQuasiDet})
 and 
Dieudonn\'e \cite{Died} determinant).
Lemma \ref{DetTopForm} provides a more conceptual approach to the
notion of determinant. It states that the determinant equals to a
coefficient of the action on the top form, exactly in the same as in
the commutative case.

{\Def
Let $M$ be a \MM.
Define the determinant of $M$ by  column expansion: \bea det
M=det^{column} M=\sum_{\sigma\in S_n} (-1)^\sigma
\prod^{\curvearrowright}_{i=1,...,n} M_{\sigma(i)i}, \eea where
$S_n$ is the group of
permutations of $n$ letters, and the
symbol $\curvearrowright$ means that in the product
$\prod_{i=1,...,n} M_{\sigma(i),i}$ one writes at first the elements
from the first column, then from the second column and so on and so
forth.}

{\Ex ~} For the case $n$=2, we have
\bea
det^{col}\left(\begin{array}{cc}
a & b\\
c & d \\
\end{array}\right)
\stackrel{def}{=} ad-cb\stackrel{lemma~ \ref{colExch} ~below
}{=}da-bc. \eea The second equality is a restatement of the second
condition of Definition \ref{MM-def}.

Let us now recall the setting of section \ref{Coact-ss} page
\pageref{Coact-ss}. Consider a Grassman algebra
$\CC[\psi_1,...,\psi_n]$ (i.e. $\psi_i^2=0, \psi_i \psi_j= -\psi_j
\psi_i$); let $\psi_i$ commute with $M_{pq}$: $\forall
i,p,q:~~[\psi_i, M_{pq}]=0$. Consider  the new variables $\tilde
\psi_i$: \bea (\tilde \psi_1, ... ,  \tilde \psi_n) = (\psi_1, ... ,
\psi_n) \left(\begin{array}{ccc}
 M_{11} & ... &  M_{1n} \\
... \\
 M_{n1} & ... &  M_{nn}
\end{array}\right).
\eea \bea \mbox{By proposition \ref{Coact-pr}, page
\pageref{Coact-pr} we have that ~~} \tilde \psi_i \tilde \psi_j =-
\tilde \psi_j \tilde \psi_i. \eea
%
%

{\Lem \label{detTopFormDef} For an arbitrary matrix $M$ (not
necessarily \MM) it holds: \label{DetTopForm}
\bea
 det^{column} (M) ~ \psi_1 \wedge ... \wedge \psi_n
= \tpsi_1 \wedge ... \wedge \tpsi_n, \eea if $M$ is a \MM, it is
true that: \bea \forall p\in S_n ~~~  det^{column} (M) ~ \psi_1
\wedge ... \wedge \psi_n = (-1)^{sgn(p)} \tpsi_{p(1)} \wedge ...
\wedge \tpsi_{p(n)}. \eea } The proof is straightforward.

{\Lem \label{colExch} The exchange of any two columns in a \MM~
changes only the sign of the determinant. More generally: an
arbitrary permutation $p$ of columns changes the determinant only by
multiplication on $(-1)^{sgn(p)}$, that is, the determinant is a
fully antisymmetric function of columns of \MMs.}

This property is specific for \MMs. For generic matrices the column
determinant is defined, but it does not satisfy this basic property.

\PRF
Let us denote a \MM~ by $M$ and by $M^p$ the matrix obtained by the permutation
of columns with respect to the permutation $p\in S_n$.
It is quite easy to see that for an arbitrary matrix $M$, (not necessarily \MM)
it is true that:
\bea
\tpsi_{p(1)} \wedge ... \wedge \tpsi_{p(n)} =  det^{column }(M^{p})
\psi_1 \wedge ... \wedge \psi_n.
\eea
For \MMs~ due to anticommutativity of $\tpsi_i$ (which is guaranteed by
proposition \ref{Coact-pr}, page \pageref{Coact-pr})
$\tpsi_{p(1)} \wedge ... \wedge \tpsi_{p(n)}
= (-1)^{sgn(p)} \tpsi_{1} \wedge ... \wedge \tpsi_{n}$
which equals to $(-1)^{p} det(M) \tpsi_{1} \wedge ... \wedge \tpsi_{n}$ by lemma \ref{detTopFormDef}.
So we conclude that $ det (M^p)= (-1)^{sgn(p)} det(M)$, for an arbitrary permutation $p\in S_n$.

\hfill \BX Another proof of this fact goes as follows. Since any
permutation can be presented as a product of transpositions of
neighbors $(i,i+1)$ it is enough to prove the proposition for such
transpositions. But for them it follows from the equality of the
commutators of cross elements (formula \ref{Cross-term-propert}).

{\Lem \label{Ex-rows-lem}  The exchange of two {\bf rows} in an {\bf
arbitrary } matrix (not necessarily \MM)
 changes only  the  sign of the $det^{{\bf column} }$.
Also if two rows in a matrix $M$ coincide, then $det^{{\bf column}
}(M)=0$. Also, one can add any row  to another row and this does not
change the determinant. More generally one can add a row multiplied
by a scalar (or more generally multiplied by an element which
commutes with $M_{ij}$) and the same is true. } The proof of these
statements is immediate. Let us stress that no conditions of
commutativity is necessary.

The lemmas above imply the following:
{\Cor ~ \label{DetPopCor} }
\begin{enumerate}
\item { 
\label{col-zero} Assume that two columns or two rows in a 
\MM~ $M$ coincide,  then $det(M)=0$.}
\item
{ 
\label{el-trans-col}
One can add any column of a \MM~ to another column and this does change the determinant.
More generally one can add a column multiplied by a constant (or more generally multiplied
by an element which commutes with $M_{ij}$) and the same is true.
}
\item { 
One can easily see that any submatrix of a \MM~ is a \MM~. So one has   natural definition
of minors and again one can choose an arbitrary order of columns (rows) to define
minors.
}
\item
{ The determinant of a \MM~ does not depend on the order of columns
in the column expansion, i.e., \bea \forall p\in S_n ~~~
det^{column} M=\sum_{\sigma\in S_n} (-1)^\sigma
\prod^{\curvearrowright}_{i=1,...,n} M_{\sigma(p(i))\ p(i)}. \eea }
\end{enumerate}



{\Prop \label{DetMultPr} Multiplicativity. Let $A$ be a \MM\ and $A'$
a generic matrix with elements in the same ring \kol. Suppose that
all elements of $A'$ commute with all elements of $A$, i.e. $\forall
i,j,k,l:~ [A_{ij},A_{kl}']=0$.
Then $det^{column}(A A')=det^{column}(A)det^{column}(A')$. If $A'$
is also a \MM, then $AA'$ is a \MM. }

{\Rem ~} If $A,B$ are the matrices such that  $\forall i,j,k,l ~[A_{ij},B_{kl}]=0$,
(for example  $B$ - $\CC$-valued matrix),
it nevertheless does not follow that $det^{column} (AB)=det^{column}(A) det^{column}( B)$,
even in $2\times 2$ case.

\PRF  The proposition
is a direct consequence of the coaction characterization of \MMs (proposition
\ref{Coact-pr}, page \pageref{Coact-pr}). The details are as follows.
\\
Consider the Grassman algebra $\Lambda[\psi_1,...,\psi_n]$ and
introduce $\tpsi_1,...,\tpsi_n$, $\ttpsi_1,...,\ttpsi_n$ as \bea
(\tilde \psi_1, ... ,  \tilde \psi_n) = (\psi_1, ... , \psi_n)
\left(\begin{array}{ccc}
 A_{11} & ... &  A_{1n} \\
... \\
 A_{n1} & ... &  A_{nn}
\end{array}\right),
~~~~~
(\ttpsi_1, ... ,  \ttpsi_n) =
(\tpsi_1, ... ,   \tpsi_n)
\left(\begin{array}{ccc}
 A_{11}' & ... &  A_{1m}' \\
... \\
 A_{n1}' & ... &  A_{nm}'
\end{array}\right).
\eea

It is easy to see that (see lemma \pref{detTopFormDef}): \bea
\label{ZZZ} \tpsi_1 \wedge ... \wedge \tpsi_n = det^{col} (A) \psi_1
\wedge ... \wedge \psi_n. \eea This equality does not require
anything except anticommutativity of $\psi_i$ and $[A_{ij},
\psi_l]=0$. (In particular we do not need $A$ to be a \MM). However,
when $A$ is a \MM, by proposition \ref{Coact-pr}, page
\pageref{Coact-pr}
 $\tpsi_i$ also anticommute; since we require that
 $\forall i,j,k,l:~ [A_{ij},A_{kl}']=0$ and $[\psi_k,A_{kl}']=0$,
we can use the same lemma again:
\bea 
\ttpsi_1 \wedge ... \wedge \ttpsi_n = det^{col} (A') \tpsi_1 \wedge ... \wedge \tpsi_n, \\
\mbox{ Using \ref{ZZZ} : ~~~~}
\ttpsi_1 \wedge ... \wedge \ttpsi_n = det^{col} (A') det^{col} (A) \psi_1 \wedge ... \wedge \psi_n.
\eea

On the other hand one can apply lemma \pref{detTopFormDef} directly to the product $(A'A)$:
\bea \ttpsi_1 \wedge ... \wedge \ttpsi_n = det^{col} (A'A) \psi_1
\wedge ... \wedge \psi_n. \eea Combining the  equalities one gets $
det^{col} (A'A)=det^{col} (A') det^{col} (A) $, which is equal to $
det^{col} (A) det^{col} (A')$ since $\forall i,j,k,l:~
[A_{ij},A_{kl}']=0$. So the first part of the proposition is proved.

Whenever  $A'$ is a \MM~ as well, one sees that $\ttpsi_i$ are again
Grassman variables (proposition \ref{Coact-pr}, page
\pageref{Coact-pr}). On the other hand $ \ttpsi_i = \sum_l \psi_l
(AA')_{li}$. So by the same proposition $AA'$ is a \MM. \BX

We can also argue $det(A A')=det( A) det( A')$ in more direct way.
One should observe that {\em all}
 elements of  $det( A) det( A')$ are contained in $det( A A')$,
but generally written in different order - the property
that $det(A)$ does not depend
on the order of  column expansion provides that
one can reorder in an appropriate way.
The property of the  column commutativity of elements of $A$
provides that unwanted terms in $det( A A')$ cancel each other.

We have already seen (observation \pref{DetMulCase22}) that
in $2\times 2$ case $det(A A')=det( A) det( A')$ for any $A'$  implies that $A$
is a \MM.
The straightforward generalization is not true in $n\times n$ case.
Indeed, consider a matrix $A$ such that all elements in some row are equal to zero.
Then, clearly $det^{col} (A)=0=det^{col}(AA')$ for any matrix $A'$.
In $2\times 2$ case any matrix with row of zeroes is a \MM, however
this is clearly not true for $3\times 3$ matrices, etc.
However for generic enough matrices $A$, such that $det^{col} (A)=det^{col}(AA')$
for any $\CC$-valued $A'$, it is true that $A$ is a \MM.
This can be seen considering $A'$ to be transposition matrices and matrices
$1+E_{ii+1}$, where as usually $E_{ij}$ matrix unit with zeroes everywhere except
$1$ at position $(ij)$.
However we do not see how this "generality condition" can be formulated in
a compact form.

{\Rem ~} Since $det^{column} M^t=det^{row} M$, where
$det^{row} M=\sum_{\sigma\in S_n} (-1)^\sigma \prod_{i=1,...,n}
M_{i\,\sigma(i)}$, the statements above can be easily  reformulated
for the case $M^t$ is a \MM.

\subsection{The permanent}
The permanent of a matrix is a
polylinear function of its columns and rows,
similar to the determinant, without sign factors $(-1)^{sign
(\sigma)}$ in its definition. We will make use of it in section
\ref{MacMahSect} below. {\Def Let $M$ be a \MM. We define its
permanent by {\bf row} expansion \footnote{Remark the difference
with the definition of the determinant,
where one uses column expansion} 
as \bea \label{M-perm-f} perm M=perm^{row} M=\sum_{\sigma\in S_n}
\prod^\curvearrowright_{i=1,...,n} M_{i\,\sigma(i)}, \eea } {\Lem
\label{perm-Col-indep} The permanent of a \MM~ does not depend on
the order of rows in the row expansion: \vskip -10mm \bea \forall
p\in S_n ~~~ perm^{row} M = \sum_{\sigma\in S_n}
\prod^\curvearrowright_{i=1,...,n} M_{p(i)\,\sigma(p(i))}, \eea or,
in the other words, the permanent of a \MM~ does not change after
arbitrary permutation of rows in a matrix $M$: \bea
 \forall p\in S_n ~~~ perm^{row} M= perm^{row} M^p,
\eea
where $M^p$ is a matrix obtained by the $p$-permutation of rows
in matrix $M$.
}

\PRF Since any permutation can be presented as a product of
transpositions of neighbors $(i,i+1)$ it is enough to prove the
proposition for such transpositions. But for them it follows from
the equality of the commutators of cross elements (formula
\ref{Cross-term-propert}). \hfill\BX
 {\Ex ~} \vskip -10mm \bea
perm^{row}\left(\begin{array}{cc}
a & b\\
c & d \\
\end{array}\right)
\stackrel{def}{=} ad+bc {=}da+cb, \eea
where the last equation follows both form the Lemma above and from the very
definition of \MM.
 {\Rem
\label{perm-Row-indep} It is easy to see that the row-permanent of
an arbitrary matrix $M$ (even without conditions of commutativity)
does not change  under any permutation of  columns. }

\subsection{Elementary properties}
Let us herewith collect some properties described above,
that are simple consequences of the
definition of \MM.

\begin{enumerate}
{\item ~ Any matrix with commuting elements is a \MM. } \vskip 0mm
\item ~ Any submatrix of a \GT\  matrix is again a \GT\ matrix.
\item ~ If $A$, $B$ are  \GT\  matrices and $\forall i,j,k,l:
[A_{ij},B_{kl}]=0$, then $A+B$  is again a \GT\ matrix.
{\item ~ If $A$ is a \GT\  matrix, $c$ is constant,  then $cA$ is a \MM. } \vskip 0mm
{\item ~ If $A$ is a \GT\  matrix, $C$ is constant matrix,  then $CA$
and $AC$ are \GT\ matrices and $det(CA)=det(AC)=det(C)det(A)$. }
\vskip 0mm
{\item ~ If $A,B$ are  \GT\  matrices and $\forall i,j,k,l:
[A_{ij},B_{kl}]=0$, then
 $AB$ is a \GT\ matrix and $det(AB)=det(A)det(B)$.
(Proposition \pref{DetMultPr}).
 }
\vskip 0mm
{\item ~ If $A$ is a
\GT\  matrix, then one can exchange the $i$-th and the $j$-th
columns(rows);
 one can put
  $i$-th column(row) on $j$-th place (erasing $j$-th column(row));
one can add new column(row) to matrix $A$ which is equal to one of
the columns(rows) of the matrix $A$; one can add the $i$-th
column(row) multiplied by any constant to the $j$-th  column(row);
in all cases the resulting matrix will be again a \GT\   matrix. }
\\
 {\item ~ If $A$ and simultaneously $A^t$  are \GT\
matrices, then all elements $A_{ij}$ commute with each other. (A
q-analog of this lemma says that if  $A$ and simultaneously $A^t$
are q-\GT\ , then $A$ is quantum matrix: $"R\one A \two A=\two A
\one A R"$ (See
\href{http://www.numdam.org/item?id=AIF_1987__37_4_191_0}
{Yu. Manin} \cite{Manin87}, \cite{Manin} 
)).}
\item
The exchange of two columns in a \MM~ changes the  sign of the
determinant. If
 two columns or two rows in a 
\MM~ $M$ coincide,  then $det(M)=0$.

This has been already discussed in corollary \ref{DetPopCor} page \pageref{DetPopCor}.
\end{enumerate}

\subsubsection{Some No-Go facts for \MMs.\label{nogosect}}
Let $M$ be a \MM\  with elements in the associative ring \kol.
\\{\bf Fact}  In general $det(M)$ is not a central element of \kol.
{\Rem ~} This should be compared with the quantum matrix group
$Fun_q(GL_n)$, where $det_q$ is central. The reason why this
property does not hold for \MMs\ is that their defining relations
are half of those of quantum matrix groups.

\vskip 0.5cm

{\bf Fact} In general $[Tr M^k, Tr M^m]\neq 0$, ~ $[Tr M, det( M)]\neq 0$.

{\Rem ~} Taking traces of powers of Lax matrices is the standard
procedure for obtaining commuting integrals of motion for integrable
systems. Indeed, Manin's conditions (or their q-analogs) do not
imply commutativity of traces. Although the concept of \MM\ is
related with (quantum) integrable systems, for this commutativity
property one needs stronger conditions like the Yang-Baxter relation
$R\one T \two T= \two T\one TR$.

{\bf Fact} In general $M^k, k=2,\ldots,$ is not a \MM.
We will prove however that $M^{-1}$ is again a \MM.
{\Rem ~} For \cite{FRT} ``quantum matrices'' such that
$T$: $R_q\one T \two T=\two T \one T R_q$, with $R_q$ a solution of the
Yang-Baxter equation it is true that $T^k$ satisfies $R_{\tilde q} \one T \two
T= \two T\one T R_{\tilde q}$, for $\tilde q=q^k$ (see e.g. \cite{ManinBook91}
section 4.2.9 pages 132-133).

{\bf Fact} Let $M$ be a \MM; then in general $det(e^M)\neq
e^{Tr(M)}$, $log(det(M))\neq {Tr(log(M))}$ (see page  \pageref{CEx1}).

\subsection{Examples. 
\label{CF-ss}} {\Def \label{CF-def} A matrix $A$ with the elements
in \kol\ is called a Cartier-Foata (see \cite{CF69,Fo79}) matrix if
elements from different {\em rows}
commute with each other.}\\
\begin{Prop}\label{CaFo} Any Cartier-Foata matrix  is a 
\MM.
\end{Prop}
{\bf Proof}. Clear form the definitions.\BX

Consider arbitrary elements $r_1,...,r_n$ in a ring \kol, then the matrix below
is obviously a $\MM$:
\bea \label{ExManMat11}
\left(\begin{array}{cccc}
r_1 & r_2 & ...& r_n \\
r_1 & r_2 & ...& r_n \\
... &...& ... & ...  \\
r_1 & r_2 & ...& r_n \\
\end{array}\right)
=
\left(\begin{array}{cccc}
1 \\
1 \\
...\\
1 \\
\end{array}\right)
\left(\begin{array}{cccc}
r_1 & r_2 & ...& r_n \\
\end{array}\right).
\eea

Consider arbitrary \MM~$M$,
consider elements $c^{i}_{jk}$ such that they commute with each other and with all elements $M$
(for example $c^{i}_{jk}$ are scalars),
then:
\bea \label{ExManMat22}
C^1 M C^2+C^3, \mbox{~~~ - is a \MM},
\eea
where $(C^i)_{jk}=c^{i}_{jk}$.

Inverse (if it exists) of a \MM~ is again a \MM~
(under natural conditions - see theorem \pref{Inverse-conj}).
Also recall that any submatrix of a \MM~ is again \MM, so
\bea \label{ExManMat33}
submatrix~of~(C^1 M C^2+C^3)^{-1} \mbox{~~~ - is a \MM}.
\eea

However if $M$ is of the form (\ref{ExManMat11}), consideration of inverses and
submatrices does not give something really new - such examples will be close
to ones defined by  (\ref{ExManMat11}).


{\Observ ~} Let $M$ be a \MM~over some ring \kol, consider arbitrary $x\in R$,
then the matrix $(\tilde M)_{ij}=x M_{ij} x^{-1}$ is also obviously a \MM.

\paragraph{Examples related to Lie algebras and integrable systems}
Let us give some remarkable examples of \MMs. They are related to
integrable systems and Lie algebras (see  \cite{CF07} section 3 for
further information).

Let $x_{ij}, y_{ij}$ be commutative variables. Let $X,Y$ be $n\times
k$ matrices with matrix elements $(X)_{ij}=x_{ij}, (Y)_{ij}=
y_{ij}$. Let us denote by $\partial_X, \partial_Y$ the $n\times k$
matrices with matrix elements $\frac{\partial}{ \partial x_{ij}}$
and $\frac{\partial}{ \partial y_{ij}}$. Let $z$ be a variable
commuting with $y_{ij}$.
\begin{description}
\item[i)]
The following $2n\times 2k$, $(n+k)\times (n+k)$ matrices are
\MMs~(the second one is related to the Capelli identities (see
\cite{CF07}, section 4.2.2)):
\bea\label{capmat}
\left(\begin{array}{cc}
X  & \partial_Y \\
Y & \partial_X \\
\end{array}\right), ~~~~~~~
\left(\begin{array}{cc}
z ~ 1_{k\times k}  & (\partial_Y)^t \\
Y & \partial_z ~ 1_{n\times n}\\
\end{array}\right).
\eea
\item[ii)]
Let  $K^1,K^2$ be $n\times n$, respectively $k\times k$ matrices
with elements in $\CC$, one can see that the following matrix is
actually a \MM: \bea \label{GaudSimple}
 \partial_z 1_{n\times n} +K^1 - Y ( z  1_{k\times k} +K^2)^{-1}(\partial_Y)^t.
\eea For the sake of concreteness, we notice that, in the case
$n=2$, $k=1$ the formula above yields the matrix:
\bea\label{Gaudsimple} \left(\begin{array}{cc}
 \p &  0  \\
 0 & \p
\end{array}\right)
+
\left(\begin{array}{cc}
 K^1_{11}  &  K^1_{12}  \\
 K^1_{21} & K^1_{22}
\end{array}\right)
  -
\left(\begin{array}{cc}
 \frac{y_{1} \partial_{y_{1}} } {z-k}  & \frac{y_{1} \partial_{y_{2}} } {z-k}  \\
 \frac{y_{2} \partial_{y_{1}} } {z-k}  & \frac{y_{2} \partial_{y_{2}} } {z-k}
\end{array}\right).
\eea
\item[iii)]
Consider the standard matrix units $e_{ij}$ (i.e. $n\times n$
matrices defined as $(e_{ij})_{kl}=
\delta_{ik}\delta_{jl}$\footnote{In other words, $e_{ij}$ has $1$ in
the position $i,j$ and $0$ everywhere else.}) Consider a variable
$z$ and the operator $\p$.  These elements commute with $e_{ij}$ and
satisfy $[\p,z]=1, e^{-\p} f(z)=f(z-1) e^{-\p}$. \bea
\label{Lax-gl} \p 1_{n\times n} -
 \frac{1}{z} \left(\begin{array}{ccc}
 e_{11} &  ... & e_{1n}  \\
 ... & ... &...  \\
 e_{n1}  & ... & e_{nn}
\end{array}\right),
~~~~~
f(z)~e^{-\p}\Bigl(  1_{n\times n} +
 \frac{1}{z} \left(\begin{array}{ccc}
 e_{11} &  ... & e_{1n}  \\
 ... & ... &...  \\
 e_{n1}  & ... & e_{nn}
\end{array}\right) \Bigr).
\eea are \MMs, where $f(z)$ is an arbitrary function.
\end{description}
 To check that the matrices in (\ref{capmat}) are indeed \MMs, one
only needs to use the standard commutation relations
$[\frac{\partial}{
\partial x_{ij}}, x_{kl}]=\delta_{ik} \delta_{jl}$, $[x_{ij},
x_{kl}]=[\frac{\partial}{ \partial x_{ij}}, \frac{\partial}{
\partial x_{kl}}]=0 $, (the same for $y_{ij}$ and $z$, while it is
understood that  operators referring to different "letters" commute,
e.g. $[\frac{\partial}{ \partial x_{ij}}, \frac{\partial}{ \partial
y_{kl}}]=0$. Similarly, for the matrix (\ref{Gaudsimple}).

The matrices in (\ref{Lax-gl}) read, in the $2\times 2$ and $f(z)=1$ case:
\bea
\left(\begin{array}{ccc}
 \p - \frac{1}{z} e_{11} &  - \frac{1}{z} e_{12} \\
 - \frac{1}{z}e_{21}  & \p- \frac{1}{z} e_{22}
\end{array}\right),
~~~~~~ \label{Lax-gl22}
\left(\begin{array}{ccc}
 e^{-\p} (1+\frac{1}{z} e_{11})  & e^{-\p}  \frac{1}{z} e_{12} \\
 e^{-\p} \frac{1}{z}e_{21}  & e^{-\p} (1+ \frac{1}{z} e_{22})
\end{array}\right).
\eea Let us check the Manin relations for the leftmost of these
matrices: Column 1 commutativity reads\footnote{Column 2
commutativity is completely analogous.}:\bea [\p - \frac{1}{z}
e_{11}, - \frac{1}{z}e_{21} ] = - [\p , \frac{1}{z} ] e_{21} +
\frac{1}{z^2} [e_{11}, e_{21}] = \frac{1}{z^2} e_{21} +
\frac{1}{z^2} (- e_{21}) =0. \eea  The cross-term relation $[M_{11},
M_{22}] = [M_{21}, M_{12}]$ is: \bea
 && [M_{11}, M_{22}] = [\p - \frac{1}{z} e_{11}, \p- \frac{1}{z} e_{22}]=
-[\p, \frac{1}{z} ] e_{22} - e_{11} [ \frac{1}{z} , \p]=
 \frac{1}{z^2} e_{22} -  \frac{1}{z^2} e_{11},\\
 && ~[M_{21}, M_{12}]= [- \frac{1}{z}e_{21}, - \frac{1}{z} e_{12}] =
\frac{1}{z^2} ( e_{22} - e_{12}). \eea Notice that the matrices
above are matrices with elements in $Diff(z)\otimes Mat_n$, i.e.
they belong to $Mat_n[ Diff(z)\otimes Mat_n]$, where $Diff(z)$ is an
algebra of differential operators in $z$. Actually one only needs
that elements $e_{ij}$ satisfy the commutation relations $ [e_{ij},
e_{kl}]= e_{il} \delta_{kj}-e_{kj}\delta_{il}$. For example, if one
regards the $e_{ij}$ as elements of the universal enveloping algebra
of $gl_n$ (or their images in an arbitrary representation), then the
matrices (\ref{Lax-gl}) will still be \MMs. We should stress   that
in the examples above we consider $n\times n$ matrices with elements
in $Mat_n \otimes Diff(z)$, and {\em not}  $n^2\times n^2$ matrices
with elements in $Diff(z)$. Indeed, in the second case, they do not
satisfy the Manin properties.

{\Rem ~} The examples above are intimately related to the Gaudin and
XXX-Heisenberg integrable systems (see \cite{CF07} and also \cite{CFRy}). According to a
specialization of the remarkable result of
\href{http://arxiv.org/abs/hep-th/0404153}{D. Talalaev} \cite{T04}
the differential operators $H_i(z)$ defined as $det(M)=
\sum_{k=0,...,n} H_k(z) \p^k$ (where $M$ is defined by
\ref{GaudSimple}) commute among themselves (i.e. $\forall
i,j:[H_i(z), H_j(w)]=0$). The operators $H_i(z)$ provide a full set
of quantum commuting integrals of motion for the Gaudin integrable
system. The $gl_n$ Gaudin system was introduced in \cite{Gaudin76},
but the full set of quantum integrals of motion (the so--called
"higher Gaudin Hamiltonians"), whose existence was proved in
\cite{FFR94}) and whose commutativity proved in \cite{ER95} was not
explicitly constructed before \cite{T04}. This construction has
far-reaching applications to Bethe ansatz and separation of
variables (see \cite{CT06-1,CF07}). Similar constructions also play
an important role in the Langlands correspondence and Kac-Moody
algebra theory (explicit description of the center of $U(\widehat{
gl_n})$ (\cite{CT06-1, CF07,CM})).
Further considerations concerning q- and elliptic analogs will
appear in \cite{CFRS}, \cite{RST}.

\subsection{Hopf structure 
\label{Hopf-s}}

Let us consider the algebra over $\CC$ generated by $M_{ij}$ $1\le
i,j \le n$ with relations: $[M_{ij}, M_{kl}]= [M_{kj}, M_{il}]$. One
can see that  it is a bialgebra with the coproduct
$\Delta(M_{ij})=\sum_k M_{ik}\otimes M_{kj}$. This is usually
denoted as follows: \bea \Delta (M) = M\stackrel{.}{\otimes}  M.
\eea {\Rem ~} It is easy to see that this coproduct is coassociative
(i.e. $ (\Delta \otimes 1)\otimes \Delta =  (1 \otimes
\Delta)\otimes \Delta$) and also (from proposition
\pref{DetMultPr}), that it holds $\Delta(det(M))=det(M)\otimes
det(M)$.

The natural antipode for this bialgebra should be $S(M)=M^{-1}$. So
it exists only in some "field of fractions" for the algebra
generated by $M_{ij}$.  {\Rem ~} Let us frame the above property
within the notions of noncommutative geometry (see e.g. \cite{Manin}
for an introduction).\\ Consider a group $G$, and denote by $Fun(G)$
the algebras of functions on $G$. The multiplication map: \bea m:
G\times G \to Fun(G), \eea clearly induces the dual map: \bea
\Delta: Fun(G)\to Fun(G)\times Fun(G) , \eea by the rule $\Delta(f)
[g_1,g_2]= f(g_1g_2)$. Associativity of the multiplication clearly
induces coassociativity of the comultiplication: $ (\Delta \otimes
1)\otimes \Delta =  (1 \otimes \Delta)\otimes \Delta$.

So, according to the noncommutative geometry point of view,
 one should think of bialgebras as some kind of "functions" on "noncommutative
spaces" which are not just spaces, but groups (more precisely
semi-groups, since we have not discussed inversion operation, which
corresponds to the antipode). So the algebra generated by matrix
elements $M_{ij}$ of a \MM~ can be thought as the algebra of
functions on some noncommutative (semi)group-space.  Moreover this
analogy can be continued. Consider a group $G$ acting on some set
$V$, and denote by $Fun(G), Fun(V)$ the algebras of functions on $F$
and, respectively, $V$. The action of the group $G$ on $V$ defines a
dual morphism of commutative algebras: \bea \phi: Fun(V) \to Fun(G)
\otimes Fun(V), ~~~  \phi(f)[g,v]=f(gv). \eea The condition of
action:  $g_1 (g_2 v)= (g_1 g_2) v$, implies: \bea (\Delta \otimes
1)(\phi) =  (1 \otimes \phi )(\phi). \eea

Which is now reformulated only in terms of $\phi, \Delta$, so makes sense
for an arbitrary bialgebra.

The "coaction"-proposition \pref{Coact-pr} implies that there exists
morphisms of algebras: \bea \phi_1: \CC[x_1,...,x_m]\to \CC<M_{ij}>
\otimes \CC[x_1,...,x_m], ~~~
\phi_1(x_i) = \sum_k M_{ik} x_k,\\
\phi_2: \CC[\psi_1,...,\psi_n] \to \CC<M_{ij}> \otimes \CC[\psi_1,...,\psi_n], ~~~
\phi_1(\psi_i) = \sum_k M_{ki} \psi_k
.
\eea
One can check that both of the maps satisfy
the condition: $(\Delta \otimes 1)(\phi_i) =  (1 \otimes \phi_i )(\phi_i)$, $i=1,2$.

So one can consider the maps $\phi_i$ as "coactions" of \MMs~on a
the space $\CC^n$ and its super version.

{\Rem ~} Let us also mention that there exists another coproduct. To
motivate it let us give another look on the algebra generated by
$M_{ij}$. The defining relations for $M_{ij}$ are written entirely
in terms of commutators, so the associative algebra generated by
$M_{ij}$ is the universal enveloping algebra to the Lie algebra
defined by the relations $[M_{ij}, M_{kl}]=[M_{kj}, M_{il}]$. For an
arbitrary universal enveloping algebra the coproduct can be defined
by the formula $\Delta(M_{ij})= M_{ij}\otimes 1+ 1 \otimes M_{ij}$.
It is easy to see directly (or conclude from the general properties
of Lie algebras) that this coproduct is compatible with the defining
relations, coassociative and there exist an antipode
$S(M_{ij})=-M_{ij}$ and counit $\epsilon(M_{ij})=0, \epsilon (1)=1$.
However this coproduct is not natural in Manin's framework (it is
not compatible with the coaction on $\CC[x_1,...,x_n]$ described
above).


\section{Inverse of a \MM \label{InvSect} }
In this section we discuss several facts about the inverse of a \MM.
The first one is Cramer's rule which states that an inverse matrix
can be calculated by the same formula with minors as in the
commutative case; the second one is that inverse matrix is actually
again a \MM~ under natural conditions - this fact is related
to a formula which goes back to
Lagrange, Desnanot, Jacobi and Lewis Caroll in the commutative case.
These should be considered
among the main results of the paper\footnote{These results, together with sketchy proofs,
were announced in \cite{CF07}}. 

\subsection{Cramer's formula \label{Cramer-s} and quasideterminants}
{\Prop\label{ad-inv} \cite{Manin} Let $M$ be a
\MM~ and denote by $M^{adj}$ the adjoint matrix defined in
the standard way, (i.e.
 $M^{adj}_{kl}=(-1)^{k+l}det^{column}(\widehat{M}_{l k})$) where
$\widehat{M}_{l k}$ is the $(n-1)\times (n-1)$ submatrix of $M$ obtained
removing the l-th row and the k-th column.
Then the same formula as in the commutative case holds true, that is,
\begin{equation}
\label{Ead}
 M^{adj} M= det^{column}(M) ~1_{n\times n}.
\end{equation}
Here $1_{n\times n}$ is the identity matrix of size $n$.  If $M^t$
is a \MM, then  $M^{adj}$ is defined by row-determinants and $M
M^{adj} = det^{column}(M^t) ~1_{n\times n}=det^{row}(M) ~1_{n\times
n}$.} \\
\PRF One can see that the equality $\forall~i:(M^{adj}
M)_{i\,i}=det^{col}(M)$, follows from the fact that $det^{col}(M)$
does not depend on the order of the column expansion in the
determinant. This independence was proved above (corollary
\ref{DetPopCor} page \pageref{DetPopCor}). Let us introduce a matrix
$\tilde M$ as follows. Take the matrix $M$ and set the $i$-th column
equal to the $j$-th column; denote the resulting matrix by $\tilde
M$.
Note that $det^{col}(\tilde M)=0$ precisely gives $(M^{adj}
M)_{i\,j}=0$ for $i\ne j$. To prove that $det^{col}(\tilde M)=0$ we
argue as follows. Clearly $\tilde M$ is a \MM. Lemma \ref{colExch}
page \pageref{colExch} allows to calculate the determinant taking
the elements first from $i$-th column, then $j$-th, then other
elements from the other columns. This yields that $det^{col}(\tilde
M)=0$, since it is the sum of the elements of the form
$(xy-yx)(z)=0$, where $x,y$ are the elements from the $i$-th and
$j$-th of $\tilde M$, so from $j$-th column of $M$. By column
commutativity of a \MM,  $xy-yx=0$, so $det^{col}(\tilde M)=0$.
\hfill\BX {\Rem ~} The only difference with the commutative case is
that, in the equality (\ref{Ead}) the order of the products of
$M^{adj}$ and $M$ has to be kept in mind.

{\Rem ~} In the works by Manin (see, e.g., \cite{Manin}) one can find
wider classes of matrices with noncommutative entries with
properly defined determinants and versions of the Cramer rule.\\

The question how far and whether the property
 $M^{adj} M= det^{column}(M) ~1_{n\times n}$ characterizes is open.
Observation  \pref{CramRulObser} shows that it is indeed the case for $2\times
2$ matrices. Since in higher rank case the relations coming from the
Cramer rule are of order $n$, while Manin's relations are always
quadratic, it is not obvious at all how to settle the matter.
\subsubsection{ Relation with quasideterminants 
 \label{SectQuasiDet} }
We will herewith recall a few constructions from the theory of
quasideterminants  of I. Gelfand and V. Retakh, R. Wilson and
discuss their counterparts in the case of \MMs. It is fair to say
that the general theoretical set-up of quasideterminants
\href{http://arxiv.org/abs/math.QA/0208146}{(I. Gelfand, S. Gelfand,
V. Retakh, R. Wilson)} \cite{GR91}, \cite{GR97}, \cite{GGRW02} can be
briefly presented as follows: {\em many facts of linear algebra can
be reformulated with the only use of an inverse matrix}. Thus it can
be extended to the noncommutative setup and can be  applied, for
example, to some questions considered here. We must stress the
difference between our set-up and that of \cite{GGRW02}: we consider
{\em a special class} of  matrices with noncommutative entries (the
\MMs), and for this class we can extend many facts of linear algebra
basically in the same form as in the commutative case, (in
particular, as we have seen,  there exists a well-defined notion of
determinant). On the other hand, in \cite{GGRW02} {\em generic
matrices} are considered; thus there is no natural notion of the
determinant, and facts of linear algebra are not exactly given in
the same form as in the commutative case.

Let us recall (\cite{GGRW02} definition 1.2.2 page 9) that the
$(p,q)$-th quasideterminant $|A|_{pq}$ of an invertible matrix $A$
is defined as $|A|_{pq}= (A^{-1}_{qp})^{-1}$, i.e. the inverse to
the $(q,p)$-element of the matrix inverse to A. It is also denoted
by: \bea |A|_{pq}= \left| \begin{array}{cccccc} A_{11} & A_{12} &
\ldots &  \ldots    & \ldots & A_{1n} \cr
 ...  & ...  & ... & ... & ...  & ... \cr
...        & ... & ... & \bo {A_{pq} } &... &... \cr
 ...  & ...  & ... & ... & ...  & ... \cr
\end{array} \right| , 
\eea From the Cramer rule we have {\Lem   \bea |M|_{pq}= (-1)^{p+q}
det(\widehat{M}_{pq})^{-1}det(M), \eea $\widehat{M}_{pq}$ is the
$(n-1)\times (n-1)$ submatrix of $M$ obtained by removing the $p$-th
row and the $q$-th column. } \\
Also, from Lemma \pref{blockInvLem} below, one can deduce the {\Lem
\bea |A|_{pq}= A_{pq} - A_{p*} (\widehat{A}_{pq})^{-1} A_{*q}, \eea
where $\widehat{A}_{pq}$ is the $(n-1)\times (n-1)$ submatrix of $A$
obtained removing the $p$-th row and the $q$-th column, $A_{p*}$ is
$p$-th row of $A$ without the element $A_{pq}$ and $A_{*q}$ is
$q$-th column of $A$ without the element $A_{pq}$. This is contained
\cite{GGRW02} proposition 1.2.6 page 10.}
 {\Ex ~} For $n = 2$ there
are four quasi-determinants: \bea \left|
\begin{array}{cc}
\bo{a_{11}}& a_{12} \cr
a_{21}& a_{22} \cr
\end{array}
\right|
=a_{11} - a_{12}\, a_{22}^{-1}\, a_{21}\ , \ \
\left|
\begin{array}{cc}
a_{11}& \bo{a_{12}} \cr
a_{21}& a_{22} \cr
\end{array}
\right|
=a_{12} - a_{11}\, a_{21}^{-1}\, a_{22}\ ,
\eea
\bea
\left|
\begin{array}{cc}
a_{11}& a_{12} \cr
\bo{a_{21}}& a_{22} \cr
\end{array}
\right|
=a_{21} - a_{22}\, a_{12}^{-1}\, a_{11}\ , \ \
\left|
\begin{array}{cc}
a_{11}& a_{12} \cr
a_{21}& \bo{a_{22}} \cr
\end{array}
\right| =a_{22} - a_{21}\, a_{11}^{-1}\, a_{12}\ . \eea

The following lemma is often useful in applications of
quasideterminants to determinants \cite{GGRW02}. It holds thanks to
the Cramer  rule for \MMs.
{\Lem \label{DetQuasDetLem}  \bea && det \left|
\begin{array}{cccc}
M_{11}& M_{12} &\ldots & M_{1n}\cr
M_{21}& M_{22} &\ldots & M_{2n}\cr
\vdots &\vdots &\ddots & \vdots\cr
M_{n1}& M_{n2} &\ldots & M_{nn}\cr
\end{array}
\right|
= \\ = &&
M_{nn}
\left|\begin{array}{cccc}
\bo{M_{n-1\,n-1}} & M_{n-1\,n}\cr
M_{n\,n-1} & M_{nn}\cr
\end{array}
\right|
\ \ldots \
\left|\begin{array}{cccc}
\bo{M_{22}}& \ldots & M_{2n}\cr
\vdots &\ddots & \vdots\cr
M_{n2}& \ldots & M_{nn}\cr
\end{array}
\right|
\left|\begin{array}{cccc}
\bo{M_{11}}& M_{12} &\ldots & M_{1n}\cr
M_{21}& M_{22} &\ldots & M_{2n}\cr
\vdots &\vdots &\ddots & \vdots\cr
M_{n1}& M_{n2} &\ldots & M_{nn}\cr
\end{array}
\right|
= \\ = &&
M_{11}
\left|\begin{array}{cccc}
{M_{1\,1}} & M_{1\,2}\cr
M_{2\,1} & \bo{M_{22}}\cr
\end{array}
\right|
\ \ldots \
\left|\begin{array}{cccc}
M_{11}& \ldots & M_{1 n-1}\cr
\vdots &\ddots & \vdots\cr
M_{n-11}& \ldots &\bo{ M_{n-1n-1} }\cr
\end{array}
\right|
\left|\begin{array}{cccc}
{M_{11}}& M_{12} &\ldots & M_{1n}\cr
M_{21}& M_{22} &\ldots & M_{2n}\cr
\vdots &\vdots &\ddots & \vdots\cr
M_{n1}& M_{n2} &\ldots & \bo{ M_{nn}}\cr
\end{array}
\right| \ . \eea}

{\Ex ~} For $2\times 2$ \MMs: \bea ad-cb= det^{col} \left|
\begin{array}{cc}
a & b \cr
c & d \cr
\end{array}
\right|
=
a \left| \begin{array}{cc}
a & b \cr
c & \bo{ d}  \cr
\end{array}
\right|
=
d \left| \begin{array}{cc}
\bo{ a}  & b \cr
c &  d  \cr
\end{array}
\right|
=da-bc.
\eea

\subsection{Lagrange-Desnanot-Jacobi-Lewis Caroll formula \label{DJ-ss}}
Below we present two identities for  \MMs.
The first of the identities is trivial in the
commutative case, while the second  has a long story: according to
({D. Bressoud})
\cite{Bres99}
 page 111 (remarks after theorem 3.12)
Lagrange found this identity for $n=3$, Desnanot proved it for $n\le
6$, Jacobi proved the general theorem (see theorem \ref{JR-th} page
\pageref{JR-th} here), C. L. Dodgson -- better known as Lewis Caroll
-- used it to derive an algorithm for calculating determinants that
required only $2\times 2$ determinants
(\href{http://www.jstor.org/pss/112607}{"Dodgson's condensation"}
method \cite{D1866}). It is quite surprising  how widely  such a
simple identity appears in various fields of mathematics \cite{DI}.

The proof of the both identities consists in a small extension and
rephrasing of the arguments in the proof of Lemma 1 page 5
{\href{http://arxiv.org/abs/hep-th/0209071} {(O. Babelon, M.
Talon)}} \cite{BT02}.

{\Prop ~\label{MainLemInv}} Let $M$ be a 
\MM\  and assume that a two sided inverse matrix $M^{-1}$ exists
(i.e. $M^{-1} M=MM^{-1}=1$). Then:
\begin{enumerate}
\item
column commutativity for $M^{-1}$ holds:
\begin{equation}
 \label{MainLemInv-fml1}
(M^{-1})_{ij} (M^{-1})_{kj}- (M^{-1})_{kj} (M^{-1})_{ij}=0.
\end{equation}
\item The
Lagrange-Desnanot-Jacobi-Lewis Caroll formula holds for \MMs, that
is: \begin{equation} \label{MainLemInv-fml2}
\begin{split}&(M^{-1})_{ij} (M^{-1})_{kl}- (M^{-1})_{kj}
(M^{-1})_{il}= \\=& (-1)^{i+j+k+l} ~ (det(M))^{-1}
det^{col}(\hat{M}_{jl,ik}),\quad j\ne l, i \ne k\end{split}
\end{equation}
where we use the notation
\[
\hat{M}_{jl,ik}:= M_{without~j-th~and~l-th~rows;
i-th~and~k-th~columns}, \] and,  in the case $n=2$ we set by
definition $det(\hat{M}_{12,12})=1$.
\end{enumerate}
{\Rem ~} The standard formulation of the
Lagrange-Desnanot-Jacobi-Lewis Caroll formula: \begin{equation}
 det(M^{ji}) det (M^{lk})- det(M^{jk}) det(M^{li}) =
(det(M))   det(\hat{M}^{jl; ik}), \end{equation} can be retrieved
(assuming that $det(M)$ commutes with principal minors)
multiplying the identity \ref{MainLemInv-fml2} by $det(M)^2$ and
with the help of the Cramer rule. As usual, we denoted by $M^{ij}$
the submatrix without $i$-row and $j$-th column, and with
$\hat{M}^{jl; ik}$ the $(n-2)\times (n-2)$ submatrix obtained
removing the $j$-th and $l$-th rows and the $i$-th and $k$-th
columns of $M$.\\
\PRF
Let us denote by $\Delta = det^{col} (M),$ and $ \Delta_{ij} =
(-1)^{i+j} det^{col} (M_{without~i-th~row,~j-th~column})$ the
cofactor of $M_{ij}$ in $det(M)$. Consider the Grassman algebra
$\CC[\psi_1,...,\psi_n]$ (i.e. $\psi_i^2=0, \psi_i \psi_j= -\psi_j
\psi_i$); let $\psi_i$ commute with $M_{p\,q}$, i.e., $\forall
i,p,q:~~[\psi_i, M_{pq}]=0$. Consider the new variables $\tilde
\psi_i$: \bea (\tilde \psi_1, ... , \tilde \psi_n) = (\psi_1, ... ,
\psi_n) \left(\begin{array}{ccc}
 M_{11} & ... &  M_{1n} \\
... \\
 M_{n1} & ... &  M_{nn}
\end{array}\right).
\eea By proposition \ref{Coact-pr}, page \pageref{Coact-pr} it holds
$\tilde \psi_i \tilde \psi_j =- \tilde \psi_j \tilde \psi_i$
%
%

We have the equalities
\bea
&
\Delta \cdot \psi_1 \wedge ... \wedge \psi_n \stackrel{def}{=}
 det(M) ~ \psi_1 \wedge ... \wedge \psi_n
\stackrel{lem~\ref{detTopFormDef}}{=} \tpsi_1 \wedge ... \wedge \tpsi_n,& \\
\label{Deltaij-fml}
&\Delta_{ij} ~ \psi_1 \wedge ... \wedge \psi_n = (-1)^{j+1} \psi_i \wedge
\tpsi_1 \wedge ...  \wedge \widehat{\tpsi_j}
\wedge ...   \wedge \tpsi_n\> \mbox{(here $\tpsi_j$ is omitted)},
&  \\
&=\tpsi_1 \wedge ...  \wedge \stackrel{j-th~place}{\psi_i} \wedge ...   \wedge \tpsi_n.
&\eea
It is easy to see that the equalities (\ref{MainLemInv-fml1},
\ref{MainLemInv-fml2}) can be reformulated as follows:
\begin{equation} \label{LemInvReform-fml} \begin{split}&\Delta_{ji} M^{-1}_{kl}
\psi_1 \wedge ... \wedge \psi_n=\\&\\ \Delta_{jk} M^{-1}_{il} \psi_1
\wedge ...& \wedge \psi_n+ \tpsi_1 \wedge \tpsi_2 \wedge ...\wedge
\stackrel{i-th~place}{\psi_j} \wedge ...\wedge
\stackrel{k-th~place}{\psi_l} \wedge ...\wedge \tpsi_n
.
\end{split}
\end{equation} Here we obviously assume $i\neq k$, since for $i=k$ (\ref{MainLemInv}) is
tautological: $\Delta_{ji} M^{-1}_{il} = \Delta_{ji} M^{-1}_{il}$.
Let us prove \ref{LemInvReform-fml}. The definition of $\tpsi_i$ is
$ (\tilde \psi_1, ... ,  \tilde \psi_m) = (\psi_1, ... ,   \psi_n)
M.$ Multiplying this relation by $M^{-1}$ on the right we get
\[
(\tilde \psi_1, ... ,  \tilde \psi_n)M^{-1} = (\psi_1, ... , \psi_n)
,\mbox{ that is, } \sum_v \tpsi_v M^{-1}_{vl} = \psi_l.
\]
Multiplying   $\tpsi_1 \wedge \tpsi_2 \wedge ...\wedge
\stackrel{i-th~place}{\psi_j} \wedge ...\wedge
\stackrel{k-th~place}{ empty } \wedge ...\wedge \tpsi_n$ and using
the Grassmann relations $\tpsi_m^2=0, i=1,\ldots, n$ we get \[ (
\tpsi_i M^{-1}_{il} +  \tpsi_k M^{-1}_{kl}  - \psi_l) \tpsi_1 \wedge
\tpsi_2 \wedge ...\wedge \stackrel{i-th~place}{\psi_j} \wedge
...\wedge \stackrel{k-th~place}{ empty } \wedge ...\wedge \tpsi_n
=0.\] By \ref{Deltaij-fml} it gives \ref{LemInvReform-fml} and hence
proposition
\ref{MainLemInv} is proved. \BX \\


\subsection{The inverse of a \MM~ is again  a \MM \label{Inv-c-ss}} 

{\Th \label{Inverse-conj}
Let $M$ be a 
\MM, assume that two sided inverse matrix $M^{-1}$ exists
(i.e. $\exists M^{-1}:$  $M^{-1} M=MM^{-1}=1$); then $M^{-1}$ is again 
a \MM.
}

{\Rem ~} We will discuss in section \pref{LefRigInvSec}, that
for a reasonable class of rings  (which includes
main examples) left invertibility of a matrix (not necessarily Manin) implies
right invertibility, and hence for \MMs~ invertibility is implied by the invertibility
of the determinant. Let us also remark that
an analogue of the theorem above holds true for Poisson-\MMs~(see
section \pref{PoissonManinSec}).

\PRF This statement follows from
proposition \ref{MainLemInv} page \pageref{MainLemInv}\footnote{So we
provide a slightly different proof with respect to that of
\cite{CF07}.}. Column commutativity for $ M^{-1}$  is just
formula (\ref{MainLemInv-fml1}) - so it is already established.

For the cross-term relation we notice that \bea \label{MInvCrTerm}
[(M^{-1})_{ij} , (M^{-1})_{kl} ] = [(M^{-1})_{kj} ,  (M^{-1})_{il}
],
\mbox{ ~~~ can be rewritten as:} \\
\label{lll1}(M^{-1})_{ij}   (M^{-1})_{kl} - (M^{-1})_{kj}   (M^{-1})_{il} =
(M^{-1})_{kl}  (M^{-1})_{ij}   - (M^{-1})_{il} (M^{-1})_{kj},
\eea
 according to (\ref{MainLemInv-fml2})
both sides of (\ref{lll1})  are equal to
\bea (det(M))^{-1}
(-1)^{i+j+k+l} det^{col} (M_{without~j-th~and~l-th~rows;
i-th~and~k-th~columns}).
\eea
The Theorem
\ref{Inverse-conj} is proved. \BX
{\Rem ~} We will derive the formula $det(M)^{-1}=(det\ M)^{-1}$
for the $n\times n$ in the next Section.

{\Rem ~} As it is was remarked in \cite{CF07} (section 4.2.1 page
17) this theorem implies a result by
{\href{http://arxiv.org/abs/math.AG/0112276} {B. Enriquez, V.
Rubtsov}} \cite{EnriquezRubtsov01} (theorem 1.1 page 2) and
{\href{http://arxiv.org/abs/hep-th/0209071} {O. Babelon, M. Talon}}
\cite{BT02} (theorem 2 page 4) about "commutativity" of quantum
Hamiltonians satisfying separation relations, which has important
applications in the theory of quantum integrable systems.

\subsection{On left and right inverses of a matrix \label{LefRigInvSec} }

The main theorems on the inverse of a
\MM~(theorem \pref{Inverse-conj}) and  on the Schur complement (theorem \pref{det-block-prop1})
are formulated under the condition that left
and right inverse matrices exist.
The lemma below shows that for a reasonable class of rings  (which includes
main examples) left invertibility of a matrix (not necessarily Manin) implies
right invertibility, and hence for \MMs~ invertibility is implied by the invertibility
of the determinant.

{\Lem Assume that ring \kol~ is a noncommutative field (i.e.
$r \ne 0 \in \mathcal{K}: \exists~ r^{-1}: r^{-1}r=  r r^{-1} =1$),
then  for any matrix $X\in Mat_n(\mathcal{K})$  for any $n$,
if the left (right) inverse exists, then the right (respectively left) inverse exists also
and they coincide.  }

If both left and right inverses exist then for any ring
associativity guarantees that they coincide: $a_l^{-1}=a_l^{-1}(aa_r^{-1})=(a_l^{-1}a)a_r^{-1}=a_r^{-1}$.

\PRF Let us prove by induction by the size $n$ of matrix $X$. For $n=1$ the
lemma is obviously true. Consider general $n$.
At least one element in first column of $X$ is nonzero, otherwise $X$ is not left invertible.
Assume it is the element $X_{11}$, otherwise multiplying by the permutation matrix we put
non-zero element to the position $(11)$.

Denote the corresponding  blocks of the matrix $X$ as $B,C,D$ as follows:
\bea
X=
\left(\begin{array}{cccc}
X_{11} & B_{1\times n-1} \\
C_{n-1\times 1} & D_{n-1\times n-1}
\end{array}\right).
\eea
Clearly using the only condition of invertibility of $X_{11}$, one can write:
\bea
X\left(\begin{array}{cccc}
1 & 0 \\
-C X_{11}^{-1} & 1
\end{array}\right)
=
\left(\begin{array}{cccc}
X_{11} & B \\
0  & D-CX_{11}^{-1} B
\end{array}\right)
.
\eea
Matrix $X$ is invertible from the left by the assumption of the lemma,
the other matrix at the left hand side of the formula above is obviously
two-sided invertible, so  left hand side in the formula above is left invertible.
So the right hand side is left invertible and it clearly implies
that $n-1\times n-1$ matrix $D-CX_{11}^{-1} B$ is left invertible.
So by the induction it is right invertible also.
Clearly:
\bea
\left(\begin{array}{cccc}
X_{11} & B \\
0  & D-CX_{11}^{-1} B
\end{array}\right)^{-1} =
\left(\begin{array}{cccc}
X_{11}^{-1} & - X_{11}^{-1} B (D-CX_{11}^{-1}B)^{-1}  \\
0  & (D-CX_{11}^{-1}B)^{-1}
\end{array}\right),
\eea
moreover it is two-sided inverse, since element $X_{11}$ is two-sided invertible as any element in
\kol~ and $(D-CX_{11}^{-1}B)$ is  two-sided invertible  by induction.
Hence we can present matrix $X$ itself as a product of two-sided invertible matrices:
\bea
X= \left(\begin{array}{cccc}
1 & 0 \\
C X_{11}^{-1} & 1
\end{array}\right)
\left(\begin{array}{cccc}
X_{11} & B \\
0  & D-CX_{11}^{-1} B
\end{array}\right)
,
\eea
hence $X$ is two-sided invertible. \BX


\section{Schur's complement and Jacobi's ratio theorem \label{SchurSect}  } 
The main result of this section is a formula for the determinant of
a \MM~ in terms of the determinant of a submatrix and the
determinant of the so-called "Schur complement". This theorem is
equivalent to the Jacobi's ratio theorem, which expresses a minor of
an inverse matrix in term of a complementary minor of matrix itself. The formulations of the
both results are exactly the same as in the commutative case. We
start with some result on multiplicativity property of the
determinant for \MMs, which is actually a key point in the proofs of
the main theorems. We also show that the so-called
Weinstein-Aronszajn and Sylvester formulas hold true for \MMs~ and actually
follow from the main theorems.

\subsection{Multiplicativity of the Determinant for special matrices of block form}
The proposition(s) below is instrumental in the proof of the Schur
formula for the determinant of a block matrix to be discussed in the
next subsection. It can be proven in a more general form than that
strictly needed for the Schur formula, and, in our opinion, is of
some interest on its own.

{\Prop \label{block-mult-prop1} Let $M$ be an  $n\times n$ \MM, with
elements in an associative ring \kol.  Let $X$ be a $k\times (n-k)$
matrix, ($k<n$), with arbitrary matrix elements in \kol.  Then
\bea
det^{column} ( M \left(\begin{array}{cc}
1_{k\times k} &  X_{k\times n-k}\\
0_{n-k\times k}  & 1_{n-k \times n-k}
\end{array}\right) ) =det^{column} M.
\eea}

Pay attention that elements of $X$ do not need to commute with elements of $M$
- they are absolutely arbitrary and so the matrix at the left hand side is not
a \MM~ in general.

\PRF Let us first state the following simple lemmas.

{\Lem Consider elements $a_i$ such that $[a_i, a_j]=0$.
Consider a matrix with the only condition that elements in some columns $i$ and $i+1$ have
a form below, then
\bea det^{column}
\left(\begin{array}{cccc}
... &  a_1 & b_1 + a_1 x & ... \\
... &  a_2 & b_2 + a_2 x & ... \\
... &  ... & ...  & ... \\
... &  a_n & b_n + a_n x & ... \\
\end{array}\right)  =
det^{column}
\left(\begin{array}{cccc}
... &  a_1 & b_1  & ... \\
... &  a_2 & b_2 & ... \\
... &  ... & ...  & ... \\
... &  a_n & b_n   & ... \\
\end{array}\right)
.
\eea
}
Using the property that one can exchange columns of any \MM,
changing only the sign of the column-determinant, column commutativity of the elements of
\MMs~ and applying the lemma above, one arrives to the following lemma:
{\Lem \label{add-col-lem2} Assume $M$ is $n\times k$, $k<n$ \MM, $n\times (n-k-1)$ matrix $C$ is absolutely
arbitrary, as well as $n\times 1$ column $b$ and $k\times 1$ column $x$, then
\bea
det^{column}
\left(\begin{array}{cccc}
M & b+Mx & C
\end{array}\right)  =
det^{column}
\left(\begin{array}{cccc}
M & b & C
\end{array}\right)
. \eea } Where we have used the notation: \bea
\left(\begin{array}{cccc} M & b+Mx & C
\end{array}\right)  =
\left(\begin{array}{ccccccc}
M_{11} & ... &  M_{1k} & b_1 + \sum_j M_{1j} x_j & C_{11} & ...& C_{1(n-k-1)} \\
M_{21} & ... &  M_{2k} & b_2 + \sum_j M_{2j} x_j & C_{21} & ...& C_{2(n-k-1)} \\
... &  ... & ...  & ... & ... &  ... & ...  \\
M_{n1} & ... &  M_{nk} & b_n + \sum_j M_{nj} x_j & C_{n1} & ...& C_{n(n-k-1)} \\
\end{array}\right)
.
\eea

The proofs of the lemmas are trivial.

Let us also remind that without any conditions on the blocks $X,Y$
it is true that: \bea \label{block-mult} \left(\begin{array}{cc}
1 & X \\
0 & 1
\end{array}\right)
\left(\begin{array}{cc}
1 & Y \\
0 & 1
\end{array}\right)
=
\left(\begin{array}{cc}
1 & Y \\
0 & 1
\end{array}\right)
\left(\begin{array}{cc}
1 & X \\
0 & 1
\end{array}\right)
=
\left(\begin{array}{cc}
1 & X+Y \\
0 & 1
\end{array}\right)
.
\eea

Recall that we need to prove the following, for a Manin matrix $M$:
\bea det^{column} ( M \left(\begin{array}{cc}
1_{k\times k} &  X_{k\times n-k}\\
0_{n-k\times k}  & 1_{n-k \times n-k} \\
\end{array}\right) ) =det^{column} M.
\eea
Let us decompose the matrix $X$ as a sum of its columns:
\bea
X= \left(\begin{array}{cccc}
0 & ... & 0 & X_{1~~ n-k}  \\
0 & ... & 0 & X_{2~~ n-k}  \\
... & ... & ... & ...  \\
0 & ... & 0 & X_{k~~ n-k}  \\
\end{array}\right)
+ ... +
\left(\begin{array}{cccc}
X_{1 1} & 0& ... & 0   \\
X_{2 1} & 0&... & 0   \\
... & ... & ... & ...  \\
X_{k 1} & 0& ... & 0   \\
\end{array}\right)
.
\eea
Let us denote this decomposition
as:
\bea
X = X_{(n-k)}+ ... +   X_{(1)}.
\eea

According to formula \ref{block-mult} let us write corresponding multiplicative
decomposition\footnote{In our inductive argument, the chosen order is
  crucial.} :
\bea
\left(\begin{array}{cc}
1 &  X\\
0  & 1  \\
\end{array}\right)
=
\left(\begin{array}{cc}
1 &  X_{(n-k)} \\
0  & 1  \\
\end{array}\right)
\times ... \times
\left(\begin{array}{cc}
1 &  X_{(1)} \\
0  & 1  \\
\end{array}\right).
\eea

After these preliminaries the proof of  the proposition follows immediately.
Observe the following:
\bea
M \left(\begin{array}{cc}
1 &  X_{(n-k)} \\
0  & 1  \\
\end{array}\right)
\times ... \times
\left(\begin{array}{cc}
1 &  X_{(n-k-l+1)} \\
0  & 1  \\
\end{array}\right)
=
\left(\begin{array}{cccccc}
M_{11} & ... & M_{1(n-l)} &  * & ... & * \\
... & ... & ... &  ... & ... & ... \\
M_{n1} & ... & M_{n(n-l)} &  * & ... & * \\
\end{array}\right).
\eea
This means that first $n-l$ columns have not been changed,
so they satisfy Manin's properties.

Now we can apply lemma \ref{add-col-lem2}:
\bea
det^{column}( \left(\begin{array}{cccccc}
M_{11} & ... & M_{1n-l} &  * & ... & * \\
... & ... & ... &  ... & ... & ... \\
M_{n1} & ... & M_{nn-l} &  * & ... & * \\
\end{array}\right)
\left(\begin{array}{cc}
1 &  X_{n-k-l} \\
0  & 1  \\
\end{array}\right) )
=\nn \\
=det^{column}( \left(\begin{array}{cccccc}
M_{11} & ... & M_{1n-l} &  * & ... & * \\
... & ... & ... &  ... & ... & ... \\
M_{n1} & ... & M_{nn-l} &  * & ... & * \\
\end{array}\right)
).
\eea

Applying this equality for $l=1,...,k$, one finishes the proof:
\bea
det^{column} (M \left(\begin{array}{cc}
1 &  X_{(n-k)} \\
0  & 1  \\
\end{array}\right)
\times ... \times
\left(\begin{array}{cc}
1 &  X_{(1)} \\
0  & 1  \\
\end{array}\right) )
= ... = \nn\\ =
det^{column} ( M \left(\begin{array}{cc}
1 &  X_{(n-k)} \\
0  & 1  \\
\end{array}\right)
\times ... \times
\left(\begin{array}{cc}
1 &  X_{(n-k-l+1)} \\
0  & 1  \\
\end{array}\right) )
=...= det^{column} M.
\eea The proposition is proven. \BX\\
In the same manner on can prove the
following {\Prop
\begin{itemize}\item For $M^t$ a Manin matrix, and an arbitrary block $X$:
\begin{equation}
 det^{row} ( \left(\begin{array}{cc}
1_{k\times k} &  X_{k\times n-k}\\
0_{n-k\times k}  & 1_{n-k \times n-k} \\
\end{array}\right) M ) =det^{row} M.
\end{equation}
\item Defining $det^{column~reverse~order} C=\sum_{\sigma\in S_n}
(-1)^\sigma \prod_{i=n,n-1,n-2,...,1} C_{\sigma(i),i}$, e.g.
$det^{col~reverse~order}\left(\begin{array}{cc}
a & b\\
c & d \\
\end{array}\right)
= da-bc$, and analogously for $det^{row~reverse~order}$, it holds:
\begin{equation}
det^{column~reverse~order} (
\left(\begin{array}{cc} 1_{k\times k} &  0_{k\times n-k}\\
X^t_{n-k\times k}  & 1_{n-k \times n-k} \\
\end{array}\right) M ) =det^{column} M,
\end{equation}
for $M$ Manin, and
\bea det^{row~reverse~order} ( M
\left(\begin{array}{cc}
1_{k\times k} &  0_{k\times n-k}\\
X^t_{n-k\times k}  & 1_{n-k \times n-k} \\
\end{array}\right) M ) =det^{row} M, \\
\eea for $M^t$ Manin.
\end{itemize}
}

\subsection{Block matrices, Schur's formula and Jacobi's ratio theorem \label{Schur-ss} }
Here we prove a formula for the determinant of a  block \MM~ in
terms of the determinant of a submatrix and the determinant of the
Schur complement. It is one the principal results of the paper. This
fact can be equivalently reformulated as Jacobi's ratio theorem
which expresses  minors of an inverse matrix in terms of minors of
the original matrix. The formulations of the theorems are exactly
the same as in the commutative case.

In our argument, we prefer to formulate the theorems, their
corollaries, and discuss their equivalence first. However, we
will provide proofs of each of the two theorems, (see pages
\pageref{SchurPRF1}, \pageref{SchurPRF2} respectively) since these
are quite different in flavor.


{\Th \label{det-block-prop1} Consider an $n\times n$ \MM ~ $M$. Let
us denote by $A,B,C,D$ its submatrices defined by \bea M=
\left(\begin{array}{cc}
A_{k\times k}  & B_{k\times n-k} \\
C_{n-k \times k} & D_{n-k\times n-k} \\
\end{array}\right),
\eea where $k<n$. Assume that $M, A,D$ are invertible on both sides,
i.e. $\exists M^{-1}, A^{-1}, \exists D^{-1}$:
\[
M^{-1}M=MM^{-1}=1,\> A^{-1}A=AA^{-1}=1,\> D^{-1}D=D D^{-1}=1.\]
Then:
\begin{itemize}
\item
It holds: \bea
 && det^{column}(M)=det^{column}\left(\begin{array}{cc}
A & B\\
C & D \\
\end{array}\right) = \nn \\
\label{Schur-fml1}
&&= det^{column}(A) det^{column}(D- C A^{-1} B)
=det^{column}(D) det^{column}(A- B D^{-1} C).
\eea
\item
The matrices $(A- B D^{-1} C)$ and $(D- C A^{-1} B)$ are \MMs.
\end{itemize}
} {\Rem ~} Matrices $(D- C A^{-1} B),(A- B D^{-1} C)$ are called
"Schur's complements". Also, notice that Schur's formula is exactly
that of the usual commutative case.


{\Prop The following more detailed statements hold, as it can be
deduced from the proof of Theorem \ref{det-block-prop1}. \bea
\label{det-bl-prec1}
 && det^{column}(M)= \nn\\
 && =det^{column}(A) det^{column}(D- C A^{-1} B), \mbox{  ~ if ~} \exists A^{-1}: AA^{-1}=1,
\label{DetSchur1}\\
 && =det^{column}(D) det^{column}(A- B D^{-1} C), \mbox{  ~ if ~} \exists D^{-1}: DD^{-1}=1, \label{DetSchur2}\\
 && =det^{column~reverse~order}(D- C A^{-1} B) det^{column}(A), \mbox{  ~ if ~} \exists A^{-1}: A^{-1}A=1,\label{DetSchur3}\\
\label{det-bl-prec4} && =det^{column~reverse~order}(A- B D^{-1} C)
det^{column}(D), \mbox{  ~ if ~} \exists D^{-1}: D^{-1}D=1.
\label{DetSchur4} \eea We recall that $det^{column~reverse~order}
C=\sum_{\sigma\in S_n} (-1)^\sigma \prod_{i=n,n-1,n-2,...,1}
C_{\sigma(i)\,i}$,  while  $det^{col}$ is defined  via the natural
order: $\prod_{i=1,2,...,n}$.}

{\Cor \label{DetInvCor} Assume $M$ is a \MM, then:
\bea
det^{column}(M) det^{column}(M^{-1})= 1,  \mbox{  ~ if ~} \exists M^{-1}: M M^{-1}=1,\\
det^{column~reverse~order}(M^{-1})  det^{column}(M)=1 ,  \mbox{  ~
if ~} \exists M^{-1}: M^{-1} M =1.
\eea
And so if $M$ is two-sided invertible, then $det(M)$ is two-sided invertible and
\bea
det(M)^{-1}=det(M^{-1}).
\eea
}
{\bf Proof of the Corollary.}
Assume that $\exists\, M^{-1}: M  M^{-1}=1$,  consider
the $2n\times 2n$ matrix below and apply formula \ref{DetSchur1}:
\bea (-1)^{n^2}= det^{column} \left(\begin{array}{cc}
M              & 1_{n\times n}\\
1_{n\times n}  & 0_{n\times n} \\
\end{array}\right) \stackrel{by~(\ref{Schur-fml1})}{=} det^{column}(M) det^{column}(-M^{-1}).
\eea From this one concludes the first statement. Similar arguments
prove the second. \BX

Theorem \ref{det-block-prop1} can be reformulated in the form called
"Jacobi's ratio theorem"\footnote{
main formula (\ref{JacRatFml}) has been also proved for \MMs~
of the form $1- tM$, $t$ is a formal parameter, in the remarkable
paper by \href{http://arxiv.org/abs/math.CO/0703203}{M. Konvalinka}
\cite{Konvalinka07-1} (see theorem 5.2 page 13). His proof is based on
combinatorics.
}
:
{\Th \label{JR-th} Consider a \MM~ $M$
that admits a left and right inverse $M^{-1}$. Let
$det^{col}(M^{-1}_{I,J})$ be the  minor of $M^{-1}$ indexed by
$I=(i_1, ...,i_k), J=(j_1, ...,j_k)$; then
\bea \label{JacRatFml}
det^{col}(M^{-1}_{I,J})=(-1)^{\sum_l i_l+\sum_l j_l}
(det^{col}(M))^{-1} det^{col}(M_{(1,...,n)\setminus J,
(1,...,n)\setminus I}),
\eea
where $det^{col}(M_{(1,...,n)\setminus
J, (1,...,n)\setminus I})$ is the minor of the matrix $M$ indexed by
the complementary set of indices.

In other words: any minor of $M^{-1}$ equals, up to a sign, to the
product of  $(det^{col}M)^{-1}$ and the corresponding complementary
minor of the {\em transpose} of $M$\footnote{One should pay
attention to the "transposition" of indexes: the minor of $M^{-1}$
is indexed by $I,J$, but the minor of $M$ is indexed by
$(1,...,n)\setminus J$ then $(1,...,n)\setminus I$.}. }

{\bf Equivalence of Schur's and Jacobi ratio theorems.} The two
theorems are actually equivalent. To see this we need to recall the
following standard lemma, that holds without any assumptions on
commutativity of the matrix elements and matrix blocks involved.
{\Lem \label{blockInvLem} Assume that the matrix
$\left(\begin{array}{cc} A&B\\C&D\end{array}\right)$ is invertible
from both sides, as well as its submatrices $A$ and $D$. Then the
matrices $(A-BD^{-1}C),(D-CA^{-1} B)$ are also invertible from both
sides, and \bea \label{block-inv-fml} \left(\begin{array}{cc}
A & B\\
C & D \\
\end{array}\right)^{-1}
=
\left(\begin{array}{cc}
(A-BD^{-1}C)^{-1} & -A^{-1}B(D-CA^{-1} B)^{-1}\\
-D^{-1}C (A-BD^{-1}C)^{-1} & (D-CA^{-1} B)^{-1} \\
\end{array}\right).
\eea . }

\SPRF The Lemma follows from the factorization formulas below: {
\bea \left(\begin{array}{cc}
A & B\\
C & D \\
\end{array}\right)
=
\left(\begin{array}{cc}
1 & 0 \\
CA^{-1}  &  1 \\
\end{array}\right)
\left(\begin{array}{cc}
A &  B\\
0 & D-CA^{-1} B \\
\end{array}\right)
=\\=
\left(\begin{array}{cc}
A-BD^{-1}C & B \\
0  &  D \\
\end{array}\right)
\left(\begin{array}{cc}
1 &  0\\
D^{-1}C & 1 \\
\end{array}\right),
\\
\left(\begin{array}{cc}
X & Y\\
0 & Z \\
\end{array}\right)^{-1}
=
\left(\begin{array}{cc}
X^{-1} & -X^{-1}YZ^{-1}\\
0 & Z^{-1} \\
\end{array}\right).
\eea }
The lemma is proved. \BX

The equivalence of the two theorems can be given as follows.
Consider the set of indexes $I,J$ in Jacobi's ratio theorem to be
$I=(1,2,..,k), J=(1,2,..,k)$. Denote the corresponding submatrices
$M_{I,J}=A$, $M_{(1,2,..,n)\setminus I,(1,2,..,n)\setminus J}=D$,
and so on and so forth, i.e., write \bea M= \left(\begin{array}{cc}
A_{k\times k}  & B_{k\times n-k}\\
C_{n-k\times k} & D_{n-k\times n-k} \\
\end{array}\right)\label{JR-SCx}.
\eea Assume Jacobi's ratio theorem holds true. This in particular
means that $M^{-1}$ is a \MM~ if $M$ is one. According to formula
(\ref{block-inv-fml}), $M^{-1}_{I,J}= (A-BD^{-1}C)^{-1}$, so we
conclude that $(A-BD^{-1}C)^{-1}$ is a \MM, since it is a submatrix
of the \MM~ $M^{-1}$. Hence $(A-BD^{-1}C)$ is a \MM~ as well by the
first claim of Jacobi ratio theorem. Similarly, $ (D-CA^{-1}B)$ is a
\MM. So the second conclusion of theorem \ref{det-block-prop1} is
derived from Jacobi's ratio theorem.

In order to derive formulas (\ref{Schur-fml1}) in theorem \ref{det-block-prop1} from
Jacobi's ratio theorem we only observe the following.
For the case $I=(1,2,..,k), J=(1,2,..,k)$ it is true that
$M^{-1}_{I,J}= (A-BD^{-1}C)^{-1}$ by (\ref{block-inv-fml}).
So Jacobi's ratio formula (\ref{JacRatFml}) reads, in this case
\[
det(
(A-BD^{-1}C)^{-1}) = (det(M))^{-1} det(D).
\]
Since $(A-BD^{-1}C)$ is
a \MM, $det\left( (
A-BD^{-1}C)^{-1}\right)=\left(det(A-BD^{-1}C)\right)^{-1}$ we arrive
at the first claim in Theorem \ref{det-block-prop1}.

Thus we have derived Schur's complement theorem
\ref{det-block-prop1} from Jacobi's ratio theorem.

\vskip 0.5cm

Let us do the converse. Assume that the Schur's complement theorem
\ref{det-block-prop1} is true for a \MM\ $M$, with a (right and
left) inverse $M^{-1}$. Construct the $2n\times 2n$ block \MM:
\bea
 M^{ext}=\left(\begin{array}{cc}
M & 1_{n\times n}\\
1_{n\times n} & 0 \\
\end{array}\right),
\eea
its Schur's complement $D-CA^{-1}B$ is precisely $-M^{-1}$. So
from theorem \ref{det-block-prop1} we conclude that $M^{-1}$ is a
\MM. Also, we can see that $det(M^{-1})=(det(M))^{-1}$, indeed
$ 
(-1)^{n^2}=det(M^{ext})\stackrel{by~(\ref{JacRatFml})}{=}det(M) det(-M^{-1}),
$ 
and $(-1)^{n^2}=(-1)^{n}$. So the first part of Jacobi's ratio theorem
holds. To derive the second claim for the case $I=(1,2,..,k),
J=(1,2..,k)$ one uses the same arguments as in the discussion after
formula (\ref{JR-SCx}). The statement for arbitrary sets of indexes
follows from this special case by changing the order of rows and
columns, and taking into account that changing the order of rows in
$M$ implies change of order of columns in $M^{-1}$, and vice versa.

Equivalence of Schur's and Jacobi ratio theorems is established.
\BX
\subsubsection{Proof of the Schur's complement theorem. \label{SchurPRF1}}
{\bf Proof 1} Let us prove the theorem \ref{det-block-prop1} page
\pageref{det-block-prop1}. Namely, we first prove the formula
(\ref{DetSchur1}) and (\ref{DetSchur2})  page \pageref{DetSchur1}
which are the more refined statements of the first claim (formula
\ref{Schur-fml1}) in theorem \ref{det-block-prop1}. The proof uses
the same idea as in the commutative case and proposition
\ref{block-mult-prop1} page \pageref{block-mult-prop1}.

Let us consider the standard decomposition\footnote{Observe that
there is no need in any commutativity constrains, but only existence
of the (right) inverse of the upper left block $A$.}
\bea\label{standard-facto} \left(\begin{array}{cc}
A & B\\
C & D \\
\end{array}\right)
\left(\begin{array}{cc}
1 & -A^{-1} B\\
0 &  1 \\
\end{array}\right)
=
\left(\begin{array}{cc}
A & 0\\
C & D-CA^{-1} B \\
\end{array}\right).
\eea
By proposition \ref{block-mult-prop1} page \pageref{block-mult-prop1}
one gets the first equality: $det^{col} M = det^{col}( A) det^{col} (D-CA^{-1} B)$.
This is desired formula \ref{DetSchur1}  page \pageref{DetSchur1}.

To prove the second equality (
$det(M)=det^{column}(D)det^{column}(A-BD^{-1}C)$ formula
\ref{DetSchur2}  page \pageref{DetSchur2}) one observes that, since
$M$ is Manin, \bea det\left(\begin{array}{cc}
A & B\\
C & D \\
\end{array}\right)
=
(-1)^{nk} det \left(\begin{array}{cc}
B & A\\
D & C \\
\end{array}\right).
\eea by column transposition.

Now we can change the order of rows, (which is possible for the
column-determinant of any matrix) to get \bea det
\left(\begin{array}{cc}
B & A\\
D & C \\
\end{array}\right)
= (-1)^{nk}
det
\left(\begin{array}{cc}
D & C\\
B & A \\
\end{array}\right).
\eea
So we one gets:
\bea
det\left(\begin{array}{cc}
A & B\\
C & D \\
\end{array}\right) = det
\left(\begin{array}{cc}
D & C\\
B & A \\
\end{array}\right) = det^{column}(D)det^{column}(A-BD^{-1}C).
\eea The last equality holds to a factorization analogous to the one
of eq. (\ref{standard-facto}).

To prove the remaining formulas \ref{DetSchur3}, \ref{DetSchur4} page \pageref{DetSchur3} 
 one uses the same arguments as above for the decomposition:
\bea
\left(\begin{array}{cc}
1 & 0\\
-CA^{-1} & 1 \\
\end{array}\right)
\left(\begin{array}{cc}
A & B\\
C & D \\
\end{array}\right)
=
\left(\begin{array}{cc}
A & B\\
0 & D-CA^{-1}B \\
\end{array}\right),
\eea where now $A^{-1}$ is a left inverse of $A$.


We are left with proving that the  Schur's complements $(A- B D^{-1}
C)$, $(D- C A^{-1} B)$ are \MMs. To do this we recall  the Lemma
\ref{blockInvLem} page \pageref{blockInvLem}: \bea
\label{block-inv-fmlCpy} \left(\begin{array}{cc}
A & B\\
C & D \\
\end{array}\right)^{-1}
=
\left(\begin{array}{cc}
(A-BD^{-1}C)^{-1} & -A^{-1}B(D-CA^{-1} B)^{-1}\\
-D^{-1}C (A-BD^{-1}C)^{-1} & (D-CA^{-1} B)^{-1} \\
\end{array}\right).
\eea From theorem \ref{Inverse-conj} page \pageref{Inverse-conj} we
know that $M^{-1}$ is a \MM, if $M$ is a \MM ~ and $M$ is two sided
invertible. Trivially any submatrix of a \MM~ is again a \MM. So we
conclude that $(A-BD^{-1}C)^{-1}$ and $(D-CA^{-1} B)^{-1}$ are \MMs.
Applying theorem \ref{Inverse-conj} again we conclude that
$(A-BD^{-1}C)$ and $(D-CA^{-1} B)$ are \MMs. \BX 


\subsubsection{Proof 2. (Proof of the Jacobi Ratio Theorem via quasideterminants)
\label{SchurPRF2} } Let us now  give a proof of Jacobi's ratio
theorem \pref{JR-th}. Our proof is  quite simple and
non-computational, the arguments being borrowed from:
\href{http://arxiv.org/abs/q-alg/9411194}{(D. Krob, B. Leclerc)} \cite{KL94} theorem 3.2
page 17, and the theory of quasideterminants \cite{GGRW02}. Roughly
speaking it goes as follows: one represents the determinants
appearing in Jacobi's formula (\ref{JR-th}) as a product of
quasiminors by lemma \pref{DetQuasDetLem} (essentially by Cramer's
rule). Using lemma \ref{NonComJac} (to be proven below) below, the
result then follows from Cramer's rule. Of course, we treat (and
use) the fact that $M^{-1}$ is also a \MM~ as already established
(theorem \pref{Inverse-conj}).

\PRF First, let us recall the
following lemma, which is called noncommutative Jacobi's ratio
theorem (\cite{GR91}, \cite{GR97} theorem 1.3.3 page 8, \cite{KL94}
theorem 2.4 page 8) or inversion law for quasiminors (\cite{GGRW02}
theorem 1.5.4 page 19):

{\Lem \label{NonComJac} For an arbitrary invertible matrix $A$ (not
necessarily a \MM), it is true that the  $(ji)$-th quasiminor of
$A^{-1}$ is the inverse of the  "almost" complementary $(ij)$-th
quasiminor of $M$. Here by almost complementary we mean the
complementary united with the $i$-th row and $j$-th column. \bea |
A^{-1}_{P,Q}|_{ij}=   | A_{ \{1...n\} - P \bigcup i,\{1...n\}- Q
\bigcup j } |_{ji}^{-1}, \eea where $A_{I,J}$ is a submatrix indexed
by index sets $I,J$.

In particular for $P=Q={1,...,k}$ and $i=j=k$:
\bea
\left| \begin{array}{cccccc}
A_{11}^{-1} &  \ldots &  A_{1k}^{-1} \cr
 ... & ...  & ... \cr
 A_{k1}^{-1} & ... & \bo {A_{kk}^{-1} } \cr
\end{array} \right|_{kk}  
= \left| \begin{array}{cccccc}
\bo {A_{kk}} &  \ldots &  A_{kn} \cr
 ... & ...  & ... \cr
 A_{nk} & ... &  {A_{nn} } \cr
\end{array} \right|_{kk}^{-1}  . 
\eea } The proof of this lemma quite readily follows from lemma
\pref{blockInvLem}\BX.

Coming back to the Jacobi ratio theorem, we need to prove: \bea
\label{JacRatFml2} det^{col}(M^{-1}_{I,J})=(-1)^{\sum_l i_l+\sum_l
j_l} (det^{col}(M))^{-1} det^{col}(M_{(1,...,n)-J, (1,...,n)-I}),
\eea Possibly changing the order of rows and columns (which is
possible for $M$ and $M^{-1}$ since both are \MMs~(theorem
\pref{Inverse-conj})
 we reduce this
identity to the case $I=J=\{1,...,k\}$: \bea \label{JacRatFml3}
det^{col}(M^{-1}_{k})= (det^{col}(M))^{-1} det^{col}(M_{\backslash k
}), \eea where $M^{-1}_{k}$ is the submatrix of $M^{-1}$ made of the
first $k$ rows and columns, and $M_{\backslash k}$ is the submatrix
of $M$ made of the last  $(n-k)$ rows and columns. One should pay
attention to the fact that changing the order of rows implies
changing the order of  columns of the inverse matrix; this implies a
possible sign factor and explains the signs and transposition of
index sets in formula (\ref{JacRatFml2}).

By Theorem \pref{Inverse-conj} $M^{-1}$ is also a \MM, so we can use
Cramer's rule \pref{ad-inv} for it, as well as for $M$ itself.
\bea
&& det^{col}(M^{-1}_{k})=\\
&& \mbox{ Let us multiply this by  } ~ 1= \prod_{i=1,...,k-1}det( M^{-1}_{i}) / det( M^{-1}_{i})  \\
&&=M^{-1}_{11} ~ \left( M^{-1}_{11} det( M^{-1}_{2}) \right) ~ \left( det(M^{-1}_{2})^{-1} det( M^{-1}_{3}) \right)
~...~ \left( det(M^{-1}_{k-1})^{-1} det( M^{-1}_{k}) \right) =\\
&& \mbox{ By Cramer's rule: } ~ |M^{-1}_{k}|_{kk}= det(M^{-1}_{k-1})^{-1} det( M^{-1}_{k})  \\
&&=
M^{-1}_{11} ~ |M^{-1}_{2}|_{22}  ~ |M^{-1}_{3}|_{33}~ ... ~ |M^{-1}_{k}|_{kk} = \\
&& \mbox{ by Lemma \ref{NonComJac}}\nn \\
&& =
|M|_{11}^{-1} ~  |M_{\backslash 1}|_{22}^{-1} ~ |M_{\backslash 2}|_{33}^{-1} ~ ... ~  |M_{\backslash k-1}|_{kk}^{-1}
=\\
&& \mbox{ by Cramer's rule for $M$}\nn \\
&& = \left( det^{col}(M)^{-1} det^{col}(M_{\backslash 1})\right)
\left( det^{col}(M_{\backslash 1} )^{-1} det^{col}(M_{\backslash 2}) \right)
...
\left( det^{col}(M_{\backslash k-1} )^{-1} det^{col}(M_{\backslash k}) \right)
=\\
&& \mbox{ by chain cancelation:}\nn \\
&& = det^{col}(M)^{-1} det^{col}(M_{\backslash k}) . \eea Formula
(\ref{JacRatFml3}) is proved, hence theorem is proved. \BX

\subsection{The Weinstein-Aronszajn formula \label{WA-ss} }
{\Prop\label{W-A-form} Let $A,B$ be $n\times k$ and $k\times n$
\MMs~ with pairwise commuting elements: $\forall i,j,k,l: [A_{ij},
B_{kl}]=0$, then \bea det^{col}(1_{n\times n} - AB)=
det^{col}(1_{k\times k} - B A). \eea } \PRF Consider the following
matrix: \bea \label{uswz}\left(\begin{array}{cc}
1_{k\times k} & B\\
A & 1_{n \times n} \\
\end{array}\right),
\eea
it is clearly a \MM.
Applying formula \ref{Schur-fml1} page \pageref{Schur-fml1} one obtains the result.
\BX

{\Rem.}
The name Weinstein-Aronszajn formula comes from \cite{Kato} Chapter
4, section 6, in the analysis of finite rank perturbation of
operators in (infinite dimensional) Hilbert spaces. The formula is
used in the theory of integrable systems (
\href{http://arxiv.org/abs/math/0108191}{H. Flaschka, J. Millson}
\cite{FlaschkaMillson01} page 23 section 6.1,
\href{http://arxiv.org/abs/nlin/0211021}{K. Takasaki}
\cite{Takasaki02} page 10), as follows. One considers the matrix
\[
M_n:=\mathbf{1}_{n\times n}-\sum_{\alpha=1}^k x_\alpha\otimes
y_\alpha,\quad x_\alpha,\, y_\alpha\in \CC^n.\] It can be considered
as a perturbation of the  identity operator by means of the $k$ rank
$1$ operators $x_\alpha\otimes y_\alpha, \alpha=1,\ldots,k$. The
Weinstein- Aronszajn formula reads
\[
det(M_n)=\det(\mathbf{1}_{k\times k}-S_{k})
\]
where $S_k$ is the $k\times k$ matrix whose element
$S_{\alpha,\beta}$ is the "scalar" product $(x_\alpha,y_\beta)$.

To get this form from our result, one simply sets $B$ to be the
matrix whose $\alpha$-th row is $x_\alpha$, and $A$ is the matrix
whose $\beta$-th column if $y_\beta$ in the expression (\ref{uswz}).

Thus Proposition (\ref{W-A-form}) holds for this case of \MMs\ as
well.

\subsection{Sylvester's determinantal identity  \label{Sylv-ss} }
Sylvester's identity is a classical determinantal identity (e.g.
\href{http://arxiv.org/abs/math.QA/0208146}{I. Gelfand, S. Gelfand,
V. Retakh, R. Wilson} \cite{GGRW02} theorem 1.5.3 page 18). Using
combinatorial methods it has been generalized for \MMs~  of the form
$1+M$ and their's q-analogs by
\href{http://arxiv.org/abs/math.CO/0703213}{(M. Konvalinka)}
\cite{Konvalinka07-2}. Here (following the classical paper by
\href{http://www.jstor.org/pss/2004533}{E. H. Bareiss} \cite{Ba68})
we show that the identity easily follows from the Schur formula
above (theorem \pref{det-block-prop1}).

Let us first recall the commutative case: {\Th (Commutative
Sylvester's identity.)
 Let $A$ be a matrix $(a_{ij})_{m \times m}$; take $n < i,j \leq m$; denote:
\bea
A_0 = \begin{pmatrix} a_{11} & a_{12} & \cdots & a_{1n} \\
               a_{21} & a_{22} & \cdots & a_{2n} \\
               \vdots & \vdots & \ddots & \vdots \\
               a_{n1} & a_{n2} & \cdots & a_{nn} \end{pmatrix}, \quad
 a_{i*} = \begin{pmatrix} a_{i1} & a_{i2} & \cdots & a_{in} \end{pmatrix}, \quad a_{*j} = \begin{pmatrix} a_{1j} \\ a_{2j} \\ \vdots \\ a_{nj}
\end{pmatrix}.
\eea Define the $(m-n)\times (m-n)$ matrix $B$ as follows:
 \bea
B_{ij} = det \begin{pmatrix} A_0 & a_{*j}
\\ a_{i*} & a_{ij} \end{pmatrix}, \quad B = (B_{ij})_{n+1 \leq i,j \leq m}.
\eea
 Then
 \bea\label{csyl}
det B = det A \cdot (det A_0)^{m-n-1} . \eea }
{\Th (Sylvester's
identity for \MMs.)
 Let $M$ be $m\times m$ a \MM~ with right and left inverse; take $n < i,j \leq m$ and denote:
\bea
M_0 = \begin{pmatrix} M_{11} & M_{12} & \cdots & M_{1n} \\
               M_{21} & M_{22} & \cdots & M_{2n} \\
               \vdots & \vdots & \ddots & \vdots \\
               M_{n1} & M_{n2} & \cdots & M_{nn} \end{pmatrix}, \quad
 M_{i*} = \begin{pmatrix} M_{i1} & M_{i2} & \cdots & M_{in} \end{pmatrix},
\quad M_{*j} = \begin{pmatrix} M_{1j} \\ M_{2j} \\ \vdots \\ M_{nj}
\end{pmatrix}.
\eea Define the $(m-n)\times (m-n)$ matrix $B$ as follows:
 \bea\label{bsylv}
B_{ij} = (\det (M_0))^{-1} \cdot
\det \begin{pmatrix} M_0 & M_{*j}
\\ M_{i*} & M_{ij} \end{pmatrix}, \quad B = (B_{ij})_{n+1 \leq i,j \leq m}.
\eea
 Then the matrix $B$ is a \MM~ and
 \bea\label{ncsyl}
det B = (\det M_0)^{-1} \cdot det M.
\eea
}

{\Rem ~} Formula (\ref{ncsyl}) reduces to (\ref{csyl}) in the
commutative case. In the noncommutative case, the Sylvester identity
holds  in the form (\ref{ncsyl}).

\PRF Once chosen $M_0$, we consider the resulting block
decomposition of  $M$, \bea M= \left(
\begin{array}{cc}
M_0 & M_1 \\
M_2 & M_3
\end{array}\right).
\eea

The key observation is the following:

{\Lem The matrix $B$ defined by (\ref{bsylv}) equals to the Schur
complement matrix: $M_0-M_2 (M_3)^{-1} M_1$. }

To see this, we need to use Schur complement theorem  \pref{det-block-prop1} again:
 \bea
B_{ij} = (\det (M_0))^{-1} \cdot
\det \begin{pmatrix} M_0 & M_{*j}
\\ M_{i*} & M_{ij} \end{pmatrix}
= (\det (M_0))^{-1} \cdot ( (\det (M_0)) ( M_{ij} - M_{i*} M_{0}^{-1} M_{*j}) )= \\
( M_{ij} - M_{i*} M_{0}^{-1} M_{*j}) =(M_0-M_2 (M_3)^{-1} M_1)_{ij}.
\eea In particular, we used the Schur formula
$det(M)=det(A)det(D-CA^{-1}B)$ for blocks $M_0=A$, $ M_{*j} =B$,
$M_{i*}=C$, $M_{ij}=D$, the last being a $1\times 1$ matrix. So the
lemma is proved.

The theorem follows from the lemma immediately; indeed, $B=(M_0-M_2
(M_3)^{-1} M_1)$ is a \MM, since a Schur complement is a \MM~ by
theorem  \pref{det-block-prop1}. $det B = (\det M_0)^{-1} \cdot det
M$ follows from the formula for determinant of
 the Schur complements.
\BX

{\Rem Bibliographical notes. }
Sylvester's identity for quasi-determinants has been  found in \cite{GR91}
(see also \href{http://arxiv.org/abs/math.QA/0208146}{(I. Gelfand,
S. Gelfand, V. Retakh, R. Wilson)} \cite{GGRW02} theorem 1.5.2 page
18). The generalization to quantum matrices in
\href{http://arxiv.org/abs/q-alg/9411194}{(D. Krob, B. Leclerc)}
\cite{KL94} theorem 3.5 page 18. Using combinatorial methods it has
been generalized for \MMs~  of the form $1+M$ and their's q-analogs
in \href{http://arxiv.org/abs/math.CO/0703213}{(M. Konvalinka)}
\cite{Konvalinka07-2}. The identity for Yangians and q-affine
algebras can be found in section 2.12 page 18
\href{http://arxiv.org/abs/math.QA/0211288}{(A. I. Molev)}
\cite{Molev02}, section 1.12 \cite{M07} and section 3 page 8
\href{http://arxiv.org/abs/math.QA/0606121}{(M.J. Hopkins, A.I.
Molev)} \cite{MH06} respectively; for twisted Yangians in section 3
page 15 \href{http://arxiv.org/abs/math.QA/0408303}{( A.I. Molev)}
\cite{Mol04}, section 2.14 \cite{M07}. These facts are used in the so-called centralizer
construction for the corresponding algebras and have some other
applications. Commutative version of the identity is discussed in
various texts (e.g. \cite{MG85}), we followed
\href{http://www.jstor.org/pss/2004533}{(E. H. Bareiss)}
\cite{Ba68}.


\subsection{Application to numeric matrices}
Let us discuss a corollary on a calculation of the
usual determinants of specific  numeric matrices, which in principle
 might provide faster algorithm
for calculating the determinant of such matrices.

{\Prop~ Consider $nm\times nm$ matrix $\tilde M$ with elements in a commutative ring  $K$.
Divide it into $m\times m$ square blocks. Denote by
$M$ an $n \times n$ matrix over $Mat_m(K)$, which matrix elements are corresponding blocks of $\tilde M$.
Assume $M$ is an $n\times n$ \MM~ over
$Mat_m(K)$. Then the  determinant of $\tilde M$ can be calculated in two steps: first one calculates
$n\times n$-column-determinant of a corresponding
$n\times n$-matrix $M$ over $Mat_m(K)$, this determinant is itself
an $m\times m$ matrix $B$ over $K$, second
one calculates  the determinant of $B$ in the usual sense:
\bea
det_{nm\times nm}(\tilde M)= det_{m\times m}( det^{col}_{n\times n} (M)),
\eea
we denoted by $det_{r\times r}$ determinants of $r\times r$-matrices.
}

Clearly such a formula is not true in general  without assumption that $M$ is a \MM.

{\Ex Let $n=2$, so we consider $2m\times 2m$ matrix over $K$ which is divided into
$4$ blocks of size $m\times m$, and it is a $2\times 2$ \MM~over $Mat_{m}(K)$
(i.e. $[a,c]=[b,d]=0$, $[a,d]=[c,b]$), then:
\bea
det_{2m\times 2m} \left(\begin{array}{cc}
a_{m\times m}  & b_{m\times m}\\
c_{m\times m} & d_{m\times m} \\
\end{array}\right)
=det_{m \times m} ( a_{m\times m} d_{m\times m} -   c_{m\times m} b_{m\times m}).
\eea
}

\PRF
Let us fix $m$ and prove by induction in $n$. For $n=1$ the statement
is tautology.

Consider general $n$. Let us assume that $det_{nm\times nm}(\tilde M)\ne 0$ and so
it is two-sided invertible, otherwise it is quite easy to see that the  proposition is true.

 Let us denote by $B,C,D$ blocks of the matrix $M$:
\bea
M=
\left(\begin{array}{ccc}
M_{11} &...& M_{1n} \\
... &...& ... \\
M_{n1} &...& M_{nn} \\
\end{array}\right)
=\left(\begin{array}{ccc}
M_{11} & B  \\
C  & D \\
\end{array}\right)
,
\eea
here $M_{ij}$ are themselves $m\times m$ matrices over $K$.
We may assume that $M_{11}$ is invertible, otherwise one should make a permutation of
rows or columns.
The Schur complement formula (\ref{Schur-fml1}) page \pageref{Schur-fml1} can be applied for
matrices over commutative rings:
\bea \label{fmlLLLL1}
det_{nm\times nm} (\tilde M)= det_{m\times m} (M_{11}) det_{(n-1)m\times (n-1)m } (D-CM_{11}^{-1} B),
\eea
on the other hand $M$ is a \MM~over $Mat_m(K)$, so we can use the Schur complement formula
in the following way:
\bea
det_{n\times n}^{col} (M)= M_{11} det_{(n-1)\times (n-1) }^{col} (D-CM_{11}^{-1} B),
\eea
so
\bea
det_{m\times m} ( det_{n\times n}^{col} (M))=
det_{m\times m} ( M_{11} det_{(n-1)\times (n-1) }^{col} (D-CM_{11}^{-1} B))=\\ =
det_{m\times m} ( M_{11})  det_{m\times m} ( det_{(n-1)\times (n-1) }^{col} (D-CM_{11}^{-1} B)) =
\eea
by theorem \pref{det-block-prop1} on Schur complements for \MMs~ one knows that $(D-CM_{11}^{-1} B)$
is a \MM, so by the induction:
\bea
= det_{m\times m} ( M_{11})   det_{(n-1)m\times (n-1)m } (D-CM_{11}^{-1} B)).
\eea
This  coincides with (\ref{fmlLLLL1}), so the proposition is proved.
\BX

{\Rem~} Taking $R=Mat_m(K)$ and considering examples \ref{ExManMat22}
page \pageref{ExManMat22}
one obtains  examples of matrices of the form required in the proposition.
It is however still unclear to us
whether such block-\MMs~may appear in practical numerical applications.

\section{Cauchy-Binet formulae and Capelli-type identities  }
We have already discussed (proposition \pref{DetMultPr}) that $det^{col}(M Y)= det^{col}(M)det^{col}(Y)$
if $[M_{ij}, Y_{kl}]=0$ and $M$ is a \MM, and actually one can prove in the same way, the
Cauchy-Binet formulae:
$det^{col}((M Y)_{IJ})= \sum_{L} det^{col}(M_{IL})det^{col}(Y_{LJ})$.
In a recent remarkable paper
\href{http://arxiv.org/abs/0809.3516}{S. Caracciolo, A. Sportiello, A. Sokal}
 \cite{CSS08} found an unexpected  noncommutative analogue of the
Cauchy-Binet formulae for \MMs,  for the case  $[M_{ij}, Y_{kl}]\ne 0$  -- but
subjected to obey
certain conditions.
Remark that the left hand side of their formul\ae\  contains a correction: $det^{col}((M Y)_{IJ} +H)$.
As a particular case of their identity one obtains
the classical Capelli and related identities.
This section is based on our interpretation of
\cite{CSS08};
we shall give a generalization of their results and provide
a different and, according to us,
more transparent proofs. Also, we obtain similar formulas for permanents.
Our main tool is use of the Grassman algebra for calculations with the determinants
and  respectively polynomial  algebra for the case of the permanents.
The condition found in \cite{CSS08} has a natural reformulation in terms
of these algebras and implies that certain expressions (anti)commute,
as if it would be a commutative case (see lemmas \ref{reformulCSSlem1},\ref{reformulCSSlem2}
page \pageref{reformulCSSlem1} and lemma \pref{XreformulCSSlem2}).

Here and below we use the following quite standard notations.
{\Not \label{NotAIJ} ~}
Let $A$ be an $n\times m$ matrix. Consider multi-indexes $I=(i_1,...,i_{r_1})$,
$J=(j_1,...,j_{r_2})$.
We denote by $A_{IJ}$ the following $r_1\times r_2$ matrix:
\bea
(A_{IJ})_{ab}=
\left\{ \begin{array}{cc}
A_{i_aj_b}, ~~i_a\le n, j_b\le m  \\
0 , ~~ i_a>n~~or~j_b>m \\
\end{array} \right.
\eea
Note that  $I,J$ are not assumed to be ordered and any number $\alpha$ may occur
several times in sequences $I,J$.

{\Ex ~} Even if $A$ is $1\times 1$ matrix, we can construct $2\times 2$, ... matrices from it:
\bea
A_{(11)(11)}=
\left(\begin{array}{cc}
A_{11} & A_{11}\\
A_{11} & A_{11} \\
\end{array}\right),
A_{(12)(11)}=
\left(\begin{array}{cc}
A_{11} & A_{11}\\
0 & 0 \\
\end{array}\right).
\eea

{\Not \label{NotPsiA} ~}
Consider some elements $\psi_1,...,\psi_n$ and an $n \times o$ matrix $A$,
for brevity we  denote by $\psi^A_i$ the element $\sum_{k=1,...,n} \psi_k A_{ki}$
i.e. just the application of the matrix $A$ to row-vector $(\psi_1,...,\psi_n)$:
\bea
(\psi_1^A,...,\psi_m^A)=  (\psi_1,...,\psi_n) A.
\eea

Vice versa: if we have some elements $\psi_1,...,\psi_n$ and $\psi_1^A,...,\psi_m^A$ we
denote by $A$, such a matrix that $(\psi_1^A,...,\psi_m^A)=  (\psi_1,...,\psi_n) A$.

\subsection{Grassman algebra condition for  Cauchy-Binet formulae }
Here we prove Cauchy-Binet formulae
under certain conditions on matrices; furthermore we will derive
results of  \cite{CSS08} as  particular cases.

Consider an   $n\times m$ matrix $M$ and an  $m\times s$ matrix $\BB$ with elements
in some ring \kol.  Consider the Grassman algebra $\Lambda[\psi_1,..., \psi_n]$
(i.e. $\psi_i^2=0, \psi_i\psi_j=-\psi_j \psi_i$), recall that we denote by
$\psi^M_i =\sum_k \psi_k M_{ki} \in  \RRR \otimes \Lambda[\psi_1,...,\psi_n]$.
We will be interested in the situation where the matrices $M, \BB$ satisfy the following two conditions.

{\bf Condition 1.}
\bea \label{grCSS1}
\forall p,j~\exists \psi^Q_j \in \RRR \otimes \Lambda[\psi_1,...,\psi_n]:~~~~\sum_{l=1,...,m} \psi^M_l [ \BB_{lj}, \psi^M_p] = \psi^M_p \psi^Q_j.
\eea
Notice that we require that $\psi^Q_j$ does not depend on $p$.

As usually  we denote by $Q$ the $n \times s$ matrix  corresponding to the elements $\psi_i^Q$:
\bea \label{QfmlLOC}
(\psi_1^Q,...,\psi_s^Q)=(\psi_1,...,\psi_n) Q
\eea

{\bf Condition 2.}
$\psi^M_i$ and  $\psi^Q_j$ anticommute:
\bea \label{grCSS2}
\forall i,j~:~~~~\psi^M_i \psi^Q_j = -  \psi^Q_j \psi^M_i .
\eea

{\Th \label{ThGrasCSS} Assume that $M$ is an $n\times m$ \MM, $\BB$ is an $m\times s$ matrix  (not necessarily Manin),
and matrices $M$ and $\BB$ satisfy conditions (\ref{grCSS1}),(\ref{grCSS2}) above.
If $n=m=s$, then:
\bea \label{CBfmlTh1}
det^{col}(M\BB + Q~diag(n-1,n-2,...,1,0) )
= det^{col}(M) det^{col}(\BB),
\eea
where $Q$ is a matrix corresponding to the elements $\psi_i^Q$ according to formula (\ref{QfmlLOC})
and the elements $\psi_i^Q$ arise in the condition 1 above;
by $diag(a_1,a_2,...)$ we denote the matrix with $a_i$ on the diagonal and $0$'s elsewhere.

More generally the following Cauchy-Binet formulae holds.
Let $I=(i_{1}<i_2<...<i_r)$,  $J=(j_{1},...,j_r)$;
\footnote{conditions  $j_{a}<j_b$, $j_a\ne j_b$ are {\em not}  required}
be two multi-indexes $i_a\le n, j_a\le s$, $r\le n,s$
then:
\bea \label{CBfmlTh2}
det^{col}((M\BB)_{IJ} + Q_{IJ}~diag(r-1,r-2,...,1,0) )
= \sum_{L=(l_1<l_2<l_3<...<l_r\le m)}  det^{col}(M_{IL}) det^{col}(\BB_{LJ}),
\eea
here $A_{IJ}$ is a matrix such that $(A_{IJ})_{ab}=A_{i_aj_b}$,
(see the notation \pref{NotAIJ}).
}

Before giving the proof of this formula, let us present some of its corollaries.

{\Cor ~} Consider the case $m<n,s$,
\bea
M= \left( \begin{array}{ccc}
M_{11} & ... & M_{1m} \\
... & ... & ... \\
... & ... & ... \\
... & ... & ... \\
M_{n1} & ... & M_{nm} \\
\end{array} \right),
\BB=
\left( \begin{array}{ccccc}
\BB_{11} & ... & ... & ...& \BB_{1s} \\
... &  ...& ...& ... & ... \\
\BB_{m1} & ...& ...& ... & \BB_{ms} \\
\end{array} \right),
\eea
then for  any $r>m$:
\bea
det^{col}((M\BB)_{IJ} + Q_{IJ}~diag(r-1,r-2,...,1,0))=0.
\eea
Indeed, for $r>m$  there is no such $L$ that $(l_1<l_2<l_3<...<l_r\le m)$, and so there is no terms
in the sum at the right hand side.
It is the same as in the commutative case, where rank of $M\BB$ is not more than $m$ and so
all minors of the size $r>m$ are zeros.

{\Cor ~ }
Assume that the matrix $\BB$ is also a \MM. Consider the matrix $(M \BB)^{\sigma}$ obtained as an arbitrary
permutation $\sigma$ of columns of $M\BB$, then
\bea
\label{Cor2222222222}
det^{col}\Bigl((M\BB)_{IJ}^{\sigma} + Q_{IJ}^{\sigma}~diag(r-1,r-2,...,1,0) \Bigr) = \\ = (-1)^{sgn(\sigma)}
det^{col}\Bigl( (M\BB)_{IJ} + Q_{IJ}~diag(r-1,r-2,...,1,0) \Bigr).
\eea

Indeed,  it is easy to see that matrices  $M$ and $\BB^\sigma$ satisfy
conditions 1 and 2 with the matrix $Q^\sigma$ and
$(M\BB)^{\sigma}=M(\BB)^{\sigma}$, so
using formulae  (\ref{CBfmlTh2}) in the main theorem
we obtain from the left hand side of (\ref{Cor2222222222})
sum of terms $det(M_{IL}) det (\BB^\sigma_{LJ})$,
and from the right hand side we obtain   $det(M_{IL}) det (\BB_{LJ})$,
since $\BB$ is a \MM, we see that they are equal up to  $(-1)^{sgn(\sigma)}$.

Now, let us prove the theorem. 

\PRF\footnote{ Proposition \pref{ToyId} is a toy model of the theorem above its extremely simple proof
is actually a good illustration for the proof below}
The equality (\ref{CBfmlTh2}) in the theorem
 can be reformulated (with the help of (\ref{GrVsMinFml1}) below)
 in terms of the Grassman algebra as follows:
\bea \label{FmlGras}
(\psi^{M\BB}_{j_1}+(r-1)\psi^Q_{j_1}) (\psi^{M\BB}_{j_2}+(r-2)\psi^Q_{j_2})...(\psi^{M\BB}_{j_r})
= \\ =
\sum_{L=(l_1<l_2<l_3<...<l_r)} \sum_{I=(i_1<i_2<i_3<...<i_r)}
\psi_{i_1}...\psi_{i_r} det^{col}(M_{IL}) det^{col}(\BB_{LJ}).
\eea

The following formulas are well-known and obvious in the commutative case
and extend to the noncommutative case without any problem.
{\Lem \label{DetGrasmLem} Consider Grassman variables $\psi_i$ and  an $n\times m$ matrix $A$;
multi-index $J=(j_1,j_2,...,j_r)$ is arbitrary (i.e. it is {\em not} assumed that
$j_a \ne j_b$, nor $j_a <j_b$).

Assume that $\psi_i$ commute with $A_{kl}$, then:
\bea \label{GrVsMinFml1}
\psi^A_{j_1} \psi^A_{j_2} ... \psi^A_{j_r}=\sum_{L=(l_1<l_2<l_3<...<l_r\le n)} \psi_{l_1} \psi_{l_2}...
\psi_{l_r} det^{col}(A_{LJ}).
\eea

Without assumption $[\psi_i, A_{kl}]=0$ we can write the same in the following way:
\bea \label{GrVsMinFml2}
\sum_{I=(i_1,..,i_r:~ 1\le i_a \le m)} \psi_{i_1} \psi_{i_2} ... \psi_{i_r} A_{i_1j_1} A_{i_2j_2}... A_{i_rj_r}
=\sum_{L=(l_1<l_2<l_3<...<l_r\le n)} \psi_{l_1} \psi_{l_2}...
\psi_{l_r} det^{col}(A_{LJ}).
\eea
The  expansion of the column determinant with respect to the first column implies the following:
\bea \label{GrVsMinFml3}
\sum_{l_1=1,...,n} \psi_{l_1} \sum_{L^-=(l_2<l_3<...<l_r\le n)}  \psi_{l_2}...
\psi_{l_r} A_{l_1j_1} det^{col}(A_{L^-J^-})
=
\sum_{L=(l_1<l_2<l_3<...<l_r\le n)} \psi_{l_1} \psi_{l_2}...
\psi_{l_r} det^{col}(A_{LJ}),
\eea
where $J^{-}=(j_2,j_3,...,j_r)$.
}

By (\ref{GrVsMinFml1}) the right hand side of equality (\ref{FmlGras})   can be also rewritten as:
\bea
\label{RHSrewrFML}
\sum_{L=(l_1<l_2<l_3<...<l_r)}
\psi^{M}_{l_1}...\psi^{M}_{l_r} det^{col}(\BB_{LJ}).
\eea

Transform the left hand side of (\ref{FmlGras}):
\bea
(\psi^{M\BB}_{j_1}+(r-1)\psi^Q_{j_1}) (\psi^{M\BB}_{j_2}+\psi^Q_{j_2})...(\psi^{M\BB}_{j_r}+\psi^Q_{j_r})
\stackrel{induction~and~(\ref{RHSrewrFML})}{=} \\ \label{IndFml000}=
(\psi^{M\BB}_{j_1}+(r-1)\psi^Q_{j_1}) (\sum_{L^-=(l_2<l_3<...<l_r)}
\psi^{M}_{l_2}...\psi^{M}_{l_r} det^{col}(\BB_{L^-J^-}) ) =,
\eea
$\psi^{M\BB}_{j_1}= \sum_{l_1} \psi^{M}_{l_1} \BB_{l_1j_1}$,
commuting  $\BB_{l_1j_1}$ and $\psi^{M}_{l_2}...\psi^{M}_{l_r}$, we get:
\bea \label{calcGrCss}
=\sum_{l_1} \sum_{L^-=(l_2<l_3<...<l_r)} \psi^{M}_{l_1}
\psi^{M}_{l_2}...\psi^{M}_{l_r} \BB_{l_1j_1}  det^{col}(\BB_{L^-J^-})
+ \\ +
 \sum_{L^-=(l_2<l_3<...<l_r)} \Bigl(
\sum_{l_1} \psi^M_{l_1} [\BB_{l_1j_1},  \psi^{M}_{l_2}...\psi^{M}_{l_r} ]
+
(r-1)\psi^Q_{j_1} \psi^{M}_{l_2}...\psi^{M}_{l_r} \Bigr)
 det^{col}(\BB_{L^-J^-}),
\eea
{\Lem \label{AntiCommLem} Conditions (\ref{grCSS1}), (\ref{grCSS2}) in the theorem above guarantee that
\bea
\sum_{l_1} \psi^M_{l_1} [\BB_{l_1j_1},  \psi^{M}_{l_2}...\psi^{M}_{l_r} ]
+(r-1)\psi^Q_{j_1} \psi^{M}_{l_2}...\psi^{M}_{l_r}  =0.
\eea
}

{\bf Proof of the lemma.} Let us transform the first term:
\bea
\sum_{l_1}  \psi^M_{l_1} [\BB_{l_1j_1},  \psi^{M}_{l_2}...\psi^{M}_{l_r} ]
=
\eea
by the Leibniz rule:
\bea
=\sum_{l_1} \psi^M_{l_1} \sum_{p=2,...,r} \psi^{M}_{l_2} ... [\BB_{l_1j_1},\psi^{M}_{l_p}]  ...\psi^{M}_{l_r}
=
\eea
using Manin's property (proposition \pref{Coact-pr}) we know that $\psi^M_i$ anticommute among themselves
so we can move
$ \psi^M_{l_1}$ in front of $[\BB_{l_1j_1},\psi^{M}_{l_p}]$ (gaining $(-1)^{p-2}$).
Using  condition (\ref{grCSS1}) (i.e. $\sum_{l} \psi^M_l [ Y_{lj}, \psi^M_p] = \psi^M_p \psi^Q_j$)
we get:
\bea
=(-1)^{p-2}   \sum_{p=2,...,r} \psi^{M}_{l_2} ... \psi^M_{l_p} \psi^Q_{j_1}  ...\psi^{M}_{l_r}
=
\eea
using anticommutativity of $\psi^M$ and $\psi^Q$ (condition (\ref{grCSS2})),
we put $\psi^Q_{l_p}$ in front of the expression and gain $(-1)=(-1)^{p-2}(-1)^{p-1}$:
\bea
=(-1) \psi^Q_{j_1}  \sum_{p=2,...,r} \psi^{M}_{l_2} ... \psi^M_{l_p}   ...\psi^{M}_{l_r}
=
\eea
so in the sum $\sum_{p=2,...,r}$ we see all the terms are identically the same,
so we have:
\bea
=(-1)(r-1) \psi^Q_{j_1}   \psi^{M}_{l_2} ...\psi^{M}_{l_r}.
\eea
Which is exactly the right hand side in the lemma. Lemma is proved.\BX

So  we continue (\ref{calcGrCss}):
\bea
=\sum_{l_1} \sum_{L^-=(l_2<l_3<...<l_r)} \psi^{M}_{l_1}
\psi^{M}_{l_2}...\psi^{M}_{l_r} Y_{l_1j_1}  det^{col}(\BB_{L^-J^-})
+ \\ +
0 =
\eea
by formula (\ref{GrVsMinFml3})
(i.e. column expansion of the determinant) we have:
\bea
=\sum_{L=(l_1<l_2<l_3<...<l_r)} \psi^{M}_{l_1}
\psi^{M}_{l_2}...\psi^{M}_{l_r}  det^{col}(\BB_{LJ})
=
\eea
by (\ref{GrVsMinFml1})
\bea
=\sum_{L=(l_1<l_2<l_3<...<l_r)} \sum_{I=(i_1<i_2<i_3<...<i_r)}
\psi_{i_1}...\psi_{i_r} det^{col}(M_{IL}) det^{col}(\BB_{LJ}).
\eea
So we transformed the left hand side (\ref{FmlGras}) to the right hand side of (\ref{FmlGras}).
Equality (\ref{FmlGras}) is equivalent  to desired formula (\ref{CBfmlTh2}) in the theorem.
Hence the theorem is proved. \BX

{\Quest ~} Consider two \MMs~$M$, $\BB$, such that they satisfy conditions 1,2 above.
Can one develop some linear algebra (Cramer rule, Cayley-Hamilton
theorem, etc.) for  $M\BB$ (or $M\BB+Q$) ? We will see below that if $Q$ is zero
then $M\BB$ is a \MM, so the answer is affirmative. Also note that
for the Capelli case (i.e.   $M_{ij}=x_{ij}$, $\BB_{ij}=\partial_{ji}$) it is also true.

\subsection{No correction case and new \MMs}
The conditions given above are easy to check for concrete pairs $M,\BB$;
however it is not so clear how to parameterize all the solutions in a simple way.
Let briefly discuss the simplest case.

The Theorem above has the following corollary:
{\Cor  Assume that $M$ is a \MM, and $\BB$ is such that the following is true:
\bea \label{grCSS1zer}
\forall p,j:~~~~\sum_{l=1,...,m} \psi^M_l [ \BB_{lj}, \psi^M_p] = 0,
\eea
then:
\bea
 \label{CBfmlTh2zer00}
det^{col}(M\BB)
= det^{col}(M) det^{col}(\BB), \\ ~
 \label{CBfmlTh2zer}
det^{col}((M\BB)_{IJ} )
= \sum_{L=(l_1<l_2<l_3<...<l_r), l_r\le m)}  det^{col}(M_{IL}) det^{col}(\BB_{LJ}).
\eea

Moreover if additionally $\BB$ is also a \MM, then
$M\BB$ is a \MM.
}

\PRF
The requirement (\ref{grCSS1zer}) implies that condition  1 page \pageref{grCSS1} is true for $\psi^Q_i=0$.
Zero anticommutes with everything, so   condition  2 page \pageref{grCSS2} holds true as well. So we apply the theorem with the matrix $Q$ equal to zero and obtain (\ref{CBfmlTh2zer00}) and (\ref{CBfmlTh2zer}).

To prove that $M\BB$ is a \MM, it is enough to prove that $\psi_j^{M\BB}$  anticommute.
Indeed this is guaranteed by  Manin's property (proposition \pref{Coact-pr}).
\bea
\psi_{j_1}^{M\BB} \psi_{j_2}^{M\BB} \stackrel{by~(\ref{GrVsMinFml1})}{=}
\sum_{i_1,i_2} \psi_{i_1}\psi_{i_2} det^{col}((M\BB)_{(i_1i_2)(j_1j_2)} )
= \\ \stackrel{by~(\ref{CBfmlTh2zer})}{=}
\sum_{i_1,i_2} \sum_{l_1,l_2} \psi_{i_1}\psi_{i_2} det^{col}((M)_{(i_1i_2)(l_1l_2)})det^{col}((\BB)_{(l_1l_2)(j_1j_2)})
=
\eea
$\BB$ is a \MM, so the determinant changes the sign after interchange of columns:
\bea
=-
 \sum_{i_1,i_2} \sum_{l_1,l_2} \psi_{i_1}\psi_{i_2} det^{col}((M)_{(i_1i_2)(l_1l_2)}det^{col}((\BB)_{(l_1l_2)(j_2j_1)}
\eea
making the transformations in the reverse order we come to:
\bea
= - \psi_{j_2}^{M\BB} \psi_{j_1}^{M\BB}.
\eea
So $\psi_{j_2}^{M\BB}$ anticommute, so by Manin's property (proposition \pref{Coact-pr})
$M\BB$ is a \MM. \BX

Let us reformulate condition (\ref{grCSS1zer}) in several ways.

It is easy to see that  (\ref{grCSS1zer}) is equivalent to:
\bea
\sum_{l=1,...,m} M_{al} [ \BB_{lj}, M_{bp}] - M_{bl} [ \BB_{lj}, M_{ap}]  = 0.
\eea

{\Lem Assume that  we can find such  matrices $A^{pj}_{**}$, that:
\bea \label{grCSS1zerPsiM}
\forall l,p,j:~~~~[ \BB_{lj}, \psi^M_p] = \sum_{v} \psi^{M}_v A^{pj}_{vl},
\eea
then it is straight-forward to see that condition  (\ref{grCSS1zer}) is equivalent to:
\bea
\forall p,j:~~~~ A^{pj}_{vl}=  A^{pj}_{lv},
\eea
i.e.  $\forall p,j$ matrix $A^{pj}_{**}$ is symmetric.
}
This is quite  transparent condition, however some non-explicitness is hidden in the matrices $A^{pj}$.

The simplest is the following:
{\Lem
Assume that  we can find elements $f_{lpj}$, such that
\bea \label{grCSS1zerPsiPsi}
\forall l,p,j:~~~~[ \BB_{lj}, \psi^M_p] = \psi^M_l f_{lpj} ,
\eea
then condition (\ref{grCSS1zer}) is obviously satisfied.
}

This case corresponds to the previous with $A^{pj}$ being some diagonal matrices.

{\Ex ~~}
Consider $\CC[x_{ij}]$, the matrix $M$: $M_{ij}=x_{ij}$, and the
the operators $R_{lp}= \sum_{k} x_{kl} \partial_{kp}$.
One can easily see that $[R_{lj}, \psi^M_p]= \delta_{pj} \psi_l^M$.
Consider the matrices $M$, $\BB$:
\bea
M_{ij}=x_{ij}, \BB_{lj} =\sum_{p} f_{ljp}(x_{ij})  R_{lp},
\eea
then one can see that they satisfy the condition \ref{grCSS1zerPsiPsi} above.
So $det(M\BB)= det(M) det(\BB)$ and more generally Cauchy-Binet formulae holds true.


\subsection{Capelli-Caracciolo-Sportiello-Sokal  case  }
Let first recall the \href{http://gdz.sub.uni-goettingen.de/index.php?id=11&no_cache=1&IDDOC=26783&IDDOC=26783&branch=&L=1}
{Capelli}
 identity \cite{Ca87}, then
its remarkable generalization  \href{http://arxiv.org/abs/0809.3516}{(S. Caracciolo, A. Sportiello, A. Sokal)}
 \cite{CSS08} and explain how it can be naturally derived within our formalism.

Consider the polynomial algebra $\CC[x_{ij}]$, and the matrices:
\bea
M= \left( \begin{array}{ccc}
x_{11} & ... & x_{1m} \\
... & ... & ... \\
x_{n1} & ... & x_{nm} \\
\end{array} \right), ~~~
\BB = \left( \begin{array}{ccc}
\frac{\partial}{\partial x_{11}}  & ... & \frac{\partial}{\partial x_{s1}} \\
... & ... & ... \\
\frac{\partial}{\partial x_{1m}}  & ... & \frac{\partial}{\partial x_{sm}} \\
\end{array} \right). ~~~
\eea

{\Th (\href{http://gdz.sub.uni-goettingen.de/index.php?id=11&no_cache=1&IDDOC=26783&IDDOC=26783&branch=&L=1}
{A. Capelli}
  \cite{Ca87})
If $n=m=s$, then:
\bea
det^{col} ( M\BB + diag(n-1,n-2,...,1,0))= det^{col} ( M) det^{col} (\BB).
\eea
And more generally ($n,m,s$ are arbitrary)  the following Cauchy-Binet formulae holds true.
Let $I=(i_{1}<i_2<...<i_r)$,  $J=(j_{1},...,j_r)$,
be two multi-indexes $i_a\le n, j_a\le s$, $r\le n,s$
then:
\bea
det^{col} ( (M\BB)_{IJ} + diag(r-1,r-2,...,1,0))= \sum_{L=(l_1<l_2<...<l_r)}
det^{col} ( M_{IL}) det^{col} (\BB_{LJ}).
\eea
}

Recently the following unexpected and general  result which includes the
Capelli identity as a particular case has been obtained.
It concerns  matrices satisfying certain commutation condition:

{\Def \label{DefCSScond} Let us say that two matrices $M,\BB$  of sizes $n\times m$ and $m\times s$ respectively satisfy the
Caracciolo-Sportiello-Sokal condition (CSS-condition for brevity),
if the following is true:
\bea \label{CSSc1}
[M_{ij}, \BB_{kl} ] = - \delta_{jk} Q_{il},
\eea
for some elements $Q_{il}$.
}

Or "in words": elements  in $j$-th column of $M$ commute with elements in $k$-th row of $\BB$
unless $j=k$, and in this case commutator of the elements $M_{ik}$ and $\BB_{kl}$ depends only on $i,l$,
but does not depend on $k$.
(See \cite{CSS08} formula (1.14) page 4, our $\BB$ is {\em transpose} to their $B$).

{\Th (\cite{CSS08} proposition 1.2' page 4.)
Assume $M$ is $n\times m$ \MM, $\BB$ is $m\times s$ matrix  (not necessarily Manin),
and matrices $M$ and $\BB$ satisfy CSS-condition, then.

If $n=m=s$:
\bea
det^{col} (M\BB+Q~diag(n-1,n-2,...,1,0))= det^{col} (M) det^{col} (\BB ),
\eea
where matrix $Q$ matrix with elements $Q_{ij}$.

More generally for arbitrary $n,m,s$ the following Cauchy-Binet formulae holds true.
Let $I=(i_{1}<i_2<...<i_r)$,  $J=(j_{1},...,j_r)$,
be two multi-indexes $i_a\le n, j_a\le s$, $r\le n,s$
then:
\bea
det^{col} ( (M\BB)_{IJ} + Q_{IJ}diag(r-1,r-2,...,1,0))= \sum_{L=(l_1<l_2<...<l_r)}
det^{col} ( M_{IL}) det^{col} (\BB_{LJ}).
\eea
}

\PRF
Let us show that CSS-theorem naturally arises from our theorem \pref{ThGrasCSS}.
To do this we need to check conditions 1 and 2 page \pageref{grCSS1}.
Let us recall the notations: $\psi_i$ are Grassman variables (i.e. $\psi_i^2=0,
\psi_i \psi_j  = - \psi_j \psi_i$); $\psi_i$ commute with $M_{ij}$ and $\BB_{kl}$.
By $\psi_i^M$ we denote $\sum_{k} \psi_{k} M_{ki}$ (see notation \pref{NotPsiA}).

Condition 1 page \pageref{grCSS1} reads:
\bea \label{grCSS1prCSS}
\forall p,j~\exists \psi^Q_j \in \RRR \otimes \Lambda[\psi_1,...,\psi_n]:~~~~
\sum_{l=1,...,m} \psi^M_l [ \BB_{lj}, \psi^M_p] = \psi^M_p \psi^Q_j.
\eea
What is the simplest and most natural way to obtain that $\sum_{l} \psi^M_l A_l   $ be proportional
 to $\psi^M_p$ ? The answer is clear: $A_l=\delta_{lp} B_l$, so let us
require, that there exists $\psi^Q_j$ such that:
\bea \label{GrCSScondReform} \label{CSSc3} \label{CSSc2}
[ \BB_{lj}, \psi^M_p] = \delta_{lp} \psi^Q_{j}.
\eea
{\Lem \label{reformulCSSlem1} It is easy to see that condition (\ref{GrCSScondReform}) is exactly equivalent
to CSS-condition (\ref{CSSc1}).}

So CSS-condition implies that condition 1 page \pageref{grCSS1} is satisfied
for matrices $M, \BB$.

Now, an unexpected fact holds true:
{\Lem \label{reformulCSSlem2} The CSS-condition automatically implies that condition 2 page \pageref{grCSS2} is also satisfied.
(i.e. $\psi_i^M$ and $\psi_j^Q$ anticommute).
}

{\bf Proof of the lemma.} If $n=1$, then  $0= \psi_i^M \psi_j^Q= -  \psi_j^Q  \psi_i^M$, so anticommutativity
holds by trivial reasons. Assume $n>1$,
take: $ l\ne i$ (it is possibly since $n>1$),  use (\ref{CSSc3}): $ \psi_i^Q= [\BB_{lj}, \psi^M_l ]$:
\bea
\psi^M_i \psi_j^h+ \psi_j^h \psi^M_i=
\psi^M_i [\BB_{lj}, \psi^M_l ] +  [\BB_{lj}, \psi^M_l ] \psi^M_i = \\ =
\psi^M_i \BB_{lj} \psi^M_l - \psi^M_i  \psi^M_l \BB_{lj}   +  \BB_{lj} \psi^M_l  \psi^M_i -
 \psi^M_l \BB_{lj} \psi^M_i=  ~~~~\mbox{ use (\ref{CSSc2}): }  [\psi^M_i, \BB_{lj}]=0  \\ =
\BB_{lj} \psi^M_i  \psi^M_l - \psi^M_i  \psi^M_l \BB_{lj}   +  \BB_{lj} \psi^M_l  \psi^M_i -
 \psi^M_l \BB_{lj}  \psi^M_i  =
 \eea
by proposition \pref{Coact-pr} $\psi^M_i$ anticommute, so cancel
$\BB_{lj} \psi^M_i  \psi^M_l+ \BB_{lj} \psi^M_l  \psi^M_i$:
\bea
  =
 - \psi^M_i  \psi^M_l \BB_{lj}    -
 \psi^M_l \BB_{lj}  \psi^M_i
=
\eea
by Manin's property (proposition \pref{Coact-pr}) $\psi^M_i$ anticommute, so:
\bea
= \psi^M_l \psi^M_i   \BB_{lj}    -
 \psi^M_l \BB_{lj}  \psi^M_i =
\psi^M_l [\psi^M_i,    \BB_{lj}] \stackrel{by~ (\ref{CSSc2})} {=} 0.
\eea
Lemma is proved. \BX

So CSS-condition implies both conditions in our theorem
\pref{ThGrasCSS}. So the CSS theorem follows from theorem \ref{ThGrasCSS} above. \BX

{\Ex ~} It is easy to see that for arbitrary functions $f_{ij}(x_{11},...,x_{nm})$
 the matrices below satisfy the CSS-condition:
\bea
M= \left( \begin{array}{ccc}
x_{11} & ... & x_{1m} \\
... & ... & ... \\
x_{n1} & ... & x_{nm} \\
\end{array} \right), ~~~
\BB = \left( \begin{array}{ccc}
\frac{\partial}{\partial x_{11}} + f_{11}(x_{ij}) & ... & \frac{\partial}{\partial x_{n1}} + f_{n1}(x_{ij}) \\
... & ... & ... \\
\frac{\partial}{\partial x_{1m}} + f_{1m}(x_{ij})  & ... & \frac{\partial}{\partial x_{nm}}+ f_{nm}(x_{ij}) \\
\end{array} \right),
\eea
and so we can apply the theorem for them.

{\Rem ~} It is obvious that the Capelli identity is a particular case of the CSS-theorem.
The first one has been widely studied and generalized
(see
\href{http://gdz.sub.uni-goettingen.de/no_cache/dms/load/img/?IDDOC=163108}{(R. Howe and T. Umeda)}
\cite{HU91} for classical reference on the subject,
and  \cite{CSS08} for quite a complete
list of references), however,  all the generalizations have been related with the Lie
algebras, super algebras, quantum groups.
It seems it was widely believed that such an identity is intimately related to
these in a sense exceptional structures. CSS-theorem shows that it is not true,
it is  actually
a particular case of the more general statement about noncommutative matrices,
which actually has nothing to do with Lie algebras or whatever.
Let us also remark that \cite{CF07} section 4.3.1 page 19 contains a very simple proof of the Capelli identity
and actually of its generalization
\href{http://arxiv.org/abs/math.QA/0610799}{(E. Mukhin, V. Tarasov, A. Varchenko)}
\cite{MTV-Cap} based on the Schur complement theorem  for \MMs.

{\Rem ~}
The CSS conditions can also be reformulated
\footnote{
$[M\otimes 1 , 1 \otimes \BB]= \sum_{ijkl} [M_{ij}, \BB_{kl}] E_{ij} \otimes E_{kl}
= - \sum_{ijkl} Q_{il} \delta_{jk} E_{ij} \otimes E_{kl}=
- \sum_{ijl} Q_{il}  E_{ij} \otimes E_{jl}$
and
$-(Q\otimes 1) P = - (\sum_{il} Q_{il} E_{il} \otimes 1)(\sum_{jk} E_{kj}\otimes E_{jk})
= -(\sum_{ij} Q_{il}  \sum_{jk}  \delta_{lk} E_{ij}\otimes E_{jk})
= -(\sum_{ijl} Q_{il}    E_{ij}\otimes E_{jl})$.
}
 with the help of matrix
notations:
\bea
[M\otimes 1 , 1 \otimes \BB] = -(Q\otimes 1) P = -P (1\otimes Q),
\eea
where $P$ is a permutation operator: $P(a\otimes b)=b\otimes a$.
However, we will not use it. Matrix (or Leningrad) notations are discussed
in section \pref{MatrSect}.

\subsection{Turnbull-Caracciolo-Sportiello-Sokal  case}
In 1948 Turnbull \cite{Tu48} proved a Capelli-type identity for symmetric matrices,
\href{http://arxiv.org/abs/math/9309212}{D. Foata and D. Zeilberger} \cite{FZ93} page 6
 gave a combinatorial proof of this identity.
\href{http://arxiv.org/abs/0809.3516}{S. Caracciolo, A. Sportiello, A. Sokal}
\cite{CSS08} proposition 1.4 page 7 proposed a generalization of this result as well.
Here we will deduce it from our theorem \pref{ThGrasCSS}.

Consider the polynomial algebra $\CC[x_{ij}]$, and the {\em symmetric} matrices:
\bea
M= \left( \begin{array}{cccccc}
x_{11} & x_{12} & x_{13} & ... & x_{1n} \\
x_{12} & x_{22} & x_{23} & ... & x_{2n} \\
x_{13} & x_{23} & x_{33} & ... & x_{3n} \\
... & ... & ... & ... & ... \\
x_{1n} & x_{2n} & x_{3n} & ... & x_{nn} \\
\end{array} \right), ~~~
\BB = \left( \begin{array}{cccccc}
\frac{\partial}{\partial x_{11}}  & \frac{\partial}{\partial x_{12}}  & \frac{\partial}{\partial x_{13}}  & ... & \frac{\partial}{\partial x_{1n}} \\
\frac{\partial}{\partial x_{12}}  & \frac{\partial}{\partial x_{22}}  & \frac{\partial}{\partial x_{23}}  & ... & \frac{\partial}{\partial x_{2n}} \\
\frac{\partial}{\partial x_{13}}  & \frac{\partial}{\partial x_{23}}  & \frac{\partial}{\partial x_{33}}  & ... & \frac{\partial}{\partial x_{3n}} \\
... & ... & ... & ... &  ... \\
\frac{\partial}{\partial x_{1n}}  & \frac{\partial}{\partial x_{2n}}  & \frac{\partial}{\partial x_{3n}}  & ... & \frac{\partial}{\partial x_{nn}} \\
\end{array} \right).
\eea

{\Th (H. W. Turnbull \cite{Tu48})
\bea
det^{col} ( M\BB + diag(n-1,n-2,...,1,0))= det^{col} ( M) det^{col} (\BB).
\eea
}

{\Def \label{TCSSdef} Let us say that matrices $M$ and $\BB$
  satisfy the
Turnbull-Caracciolo-Sportiello-Sokal condition (TCSS-condition for brevity),
if the following is true:
\bea \label{TCSSc1}
[M_{ij}, \BB_{kl} ] = - h ( \delta_{jk} \delta_{il} + \delta_{ik} \delta_{jl}),
\eea
for some element $h$.
}

Or "in words": element $M_{ij}$ commute with all elements of $\BB$
 except elements $\BB_{ik}$ and $\BB_{ki}$,
and this non-zero commutator is equal to $-h$, for some element $h$.
(See \cite{CSS08} formula (1.25) page 7, our $M$ is their $A$ and our $\BB$ is their $B$).

{\Lem \label{LemHcommute} Assume $n>1$. Assume $M$ is a symmetric $n\times n$ matrix with commuting entries (i.e. $\forall i,j,k,l:~~[M_{ij},M_{kl}]=0$),
 $M$ and $\BB$ satisfy TCSS-condition above, then:
\bea
\forall i,j:~~~~[M_{ij}, h]=0.
\eea
}
Indeed, consider $M_{ij}$, take $M_{ab}$ such that $(ij)\ne (ab)$ and $(ij)\ne (ba)$ (this is possible
since $n>1$, then $h=[\BB_{ab}, M_{ab}]$, by TCSS-condition $[M_{ij},\BB_{ab}]=0$ and by assumption
$[M_{ij}, M_{ab}]=0$, so $[M_{ij},h]= [M_{ij},[\BB_{ab}, M_{ab}]]=0$.

{\Th (\cite{CSS08} proposition 1.4 page 7.)
Assume $M$ is  $n\times n$ {\em symmetric} matrix with commuting entries (i.e. $\forall i,j,k,l:~~[M_{ij},M_{kl}]=0$)
\footnote{we can requite $M$ to be a \MM, but for the field of characteristic not equal to $2$,
symmetric \MM~ is matrix with commuting entries },
$\BB$ is $n\times o$ matrix  (not necessarily Manin),
and matrices $M$ and $\BB$ satisfy TCSS-condition, then.

If $n=o$:
\bea
det^{col} (M\BB+h~diag(n-1,n-2,...,1,0))= det^{col} (M) det^{col} (\BB ).
\eea

For arbitrary $n,o$, more generally  the following Cauchy-Binet formulae holds true.
Let $I=(i_{1}<i_2<...<i_r)$,  $J=(j_{1},...,j_r)$,
be two multi-indexes $i_a\le n, j_a\le n$, $r\le n,o$
then:
\bea
det^{col} ( (M\BB)_{IJ} + h~ diag(r-1,r-2,...,1,0))= \sum_{L=(l_1<l_2<...<l_r)}
det^{col} ( M_{IL}) det^{col} (\BB_{LJ}).
\eea
}
\PRF
We will deduce the theorem  above from our theorem \pref{ThGrasCSS}.
To do this we need to check conditions 1 and 2 page \pageref{grCSS1}.
Let us recall the notations: $\psi_i$ are Grassman variables (i.e. $\psi_i^2=0,
\psi_i \psi_j  = - \psi_j \psi_i$) and $\psi_i$ commute with $M_{ij}$ and $\BB_{kl}$.
By $\psi_i^M$ we denote $\sum_{k} \psi_{k} M_{ki}$ (see notation \pref{NotPsiA}).

From the TCSS-condition above we see that:
\bea
[\BB_{lj}, \psi^M_p]= \delta_{lp} h  \psi_j +   \delta_{jp} h \psi_l.
\eea

Let us look at the our condition 1 page \pageref{grCSS1}:
\bea
\label{LocTCSSfml}
\sum_{l=1,...,n} \psi^M_l [ \BB_{lj}, \psi^M_p] =
\sum_{l=1,...,n} \psi^M_l (\delta_{lp} h \psi_j +   \delta_{jp} h \psi_l)
\eea

{\Lem \label{LemSymMatr} Obviously for any symmetric matrix $M$ (not necessarily \MM):
\bea
\sum_{l=1,...,n} \psi_l \psi^M_l =0.
\eea
}
Indeed,  $ \sum_{l=1,...,n}\psi_l \psi^M_l =  \sum_{l=1,...,n} \sum_{j=1,...,n}   \psi_l \psi_j M_{jl}
=
\sum_{l<j}  \psi_j \psi_l  (M_{jl}-M_{lj}) =0$.

In our case $\psi_i$ and $\psi^M_j$ anticommute so the sum
$\sum_{l=1,...,n} \psi^M_l \delta_{jp} h \psi_l$ is zero.
And so we continue
(\ref{LocTCSSfml}):
\bea
\label{LocTCSSfml2}
=\sum_{l=1,...,n} \psi^M_l \delta_{lp} h \psi_j = \psi^M_p  h \psi_j.
\eea
Hence condition 1 page \pageref{grCSS1} is satisfied for $\psi^Q_j= h \psi_j$.

Condition 2 page \pageref{grCSS2} requires that $ \psi^M_p$,   $ h \psi_j$.
Assume that $n>1$, then  by lemma \ref{LemHcommute} above we know that $h$, commute with
$M_{ij}$, by definition $ \psi^M_p = \sum_k \psi_k M_{kp} $, where $\psi_i$ commute with $M_{kl}$.
These implies the condition 2, for $n>1$.

 Hence our conditions 1,2 are satisfied so applying  theorem
\pref{ThGrasCSS} we obtain the theorem above for $n>1$.
For $n=1$ the theorem above is absolutely tautological identity: $M_{11} \BB_{11}= M_{11} \BB_{11}$.
So theorem is proved. \BX

\subsection{Generalization to permanents}
We have seen that the identity above has a natural formulation and proof in terms of
the Grassman algebra, so we can also look for the similar result for the  algebra of polynomials,
since both algebras play an equal role in the definition of \MMs.
Here we will briefly discuss
 analogs for permanents of the theorems above.
Since proofs are absolutely similar we will give only the formulations
and some key comments.
One can actually consider the case of super-\MMs~ and then
both cases of determinants and permanents are the particular cases of it,
but we do not want to overload the text going into the super-theory.

\subsubsection{Preliminaries}

{\Def Let us recall that the {\em column permanent} of a square $n\times n$ matrix $A$ is:
\bea
per^{col} A = \sum_{\sigma \in S_n} A_{\sigma(1)1} A_{\sigma(2)2}... A_{\sigma(n)n},
\eea
i.e. in the product $A_{\sigma(1)1} A_{\sigma(2)2}... A_{\sigma(n)n}$ we first
take elements from the first column, than second ...
Here $S_n$ is the permutation group.
}

The definition is absolutely similar to the determinant, but
sign of the permutation does not contribute.

Let $A$ be an $n\times m$ matrix. Consider multi-indexes $I=(i_1,...,i_{r_1})$,
$J=(j_1,...,j_{r_2})$.
Let us recall (notation \pref{NotAIJ}) that we denote by $A_{IJ}$ the following $r\times r$ matrix:
\bea
(A_{IJ})_{ab}= A_{i_aj_b}.
\eea

In formulas below it will be necessary to use the following normalized version of permanents.
{\Def Let $A$ be an $n \times m$ matrix over some (not necessarily commutative ring),
let $I=(i_1\le ... \le i_r)$,
$J=(j_1,...,j_r)$, $\forall a:~j_a\le m$ (we do not require ordering nor $j_a\ne j_b$).
Let us call by the normalized column permanent $perm^{col}_{norm}(A_{IJ})$ the following:
\bea
perm^{col}_{norm}(A_{IJ})= \frac{1}{(2!)^{n_2}(3!)^{n_3}...} perm^{col}(A_{IJ}),
\eea
where $n_p$ is defined as follows: $n_p=v$, means that some numbers $a_1, ... a_v$ enter
the sequence $I$ with multiplicity exactly $p$.
}

For example for $I=(i_1,i_1,i_2,i_2,i_3,i_3,i_3)$
the factor will be $(2!)^2(3!)$. Note that  multiplicity in $J$ does not contribute.

{\Ex ~}
\bea
perm^{col}_{norm} (A_{(11)(11)})=\frac{1}{2!}
perm^{col}\left(\begin{array}{cc}
A_{11} & A_{11}\\
A_{11} & A_{11} \\
\end{array}\right)= \frac{1}{2}(A_{11}A_{11}+A_{11}A_{11})=(A_{11})^2,\\
\mbox{more generally:~~~~} perm^{col}_{norm} (A_{(aa..a)(bb...b)})=(A_{ab})^r.
\eea

{\Not \label{NotX} ~}
Consider some elements $x_1,...,x_n$ and an $n \times o$ matrix $A$,
for brevity we  denote by $x^A_i$ the element $\sum_{k=1,...,n} x_k A_{ki}$
i.e. just the application of the matrix $A$ to row-vector $(x_1,...,x_n)$:
\bea
(x_1^A,...,x_m^A)=  (x_1,...,x_n) A.
\eea

The lemma below is quite obvious. It is an analogue of lemma \pref{DetGrasmLem}
 for Grassman variables and determinants.

{\Lem Consider commuting variables $x_i$, $i=1,...,n$  and  an $n\times m$ matrix $A$;
multi-index $J=(j_1,j_2,...,j_r)$ is arbitrary (i.e. it is {\em not} assumed that
$j_a \ne j_b$, nor $j_a <j_b$).

Assume that $x_i$ commute with $A_{kl}$, then:
\bea \label{xVsPermFml1}
x^A_{j_1} x^A_{j_2} ... x^A_{j_r}=\sum_{L=(l_1\le l_2\le l_3\le ...\le l_r\le n)} x_{l_1} x_{l_2}...
x_{l_r} perm_{norm}^{col}(A_{LJ}).
\eea

Without assumption $[x_i, A_{kl}]=0$ we can write the same in the following way:
\bea \label{xVsPermFml2}
\sum_{I=(i_1,..,i_r:~ 1\le i_a \le m)} x_{i_1} x_{i_2} ... x_{i_r} A_{i_1j_1} A_{i_2j_2}... A_{i_rj_r}
=\sum_{L=(l_1\le l_2\le l_3\le ...\le l_r\le n)} x_{l_1} x_{l_2}...
x_{l_r} perm_{norm}^{col}(A_{LJ}).
\eea
It implies  the formula for the
  expansion of the column permanent with respect to the first column:
\bea \label{xVsPermFml3}
\sum_{l_1=1,...,n} x_{l_1} \sum_{L^-=(l_2\le l_3\le ...\le l_r\le n)}  x_{l_2}...
x_{l_r} A_{l_1j_1} perm_{norm}^{col}(A_{L^-J^-})
=
\sum_{L=(l_1\le l_2\le l_3\le ...\le l_r\le n)} x_{l_1} x_{l_2}...
x_{l_r} perm_{norm}^{col}(A_{LJ}),
\eea
where $J^{-}=(j_2,j_3,...,j_r)$.
}

\subsubsection{Cauchy-Binet type formulas for permanents}

{\Prop (Easy Cauchy-Binet formula for permanents). Consider an $m\times n$ \MM~$M$
and an arbitrary $m\times o$ matrix $\BB$, such that $[M_{ij}, \BB_{kl}]=0$,
consider an arbitrary multi-index $J=(j_1,...,j_r)$ and ordered multi-index $I=(i_1\le l_2 \le ... \le i_r)$, $i_a\le n$,
then:
\bea
perm^{col}_{norm}((M^t \BB)_{IJ})= \sum_{L=(l_1\le l_2 \le ... \le l_r), l_a\le m}
 perm^{col}_{norm}((M^t)_{IL}) perm^{col}_{norm}((\BB)_{LJ}).
\eea
}
\PRF It is quite obvious let us nevertheless write it up.
Consider $\CC[x_1,...,x_n]$, such that $[x_i, M_{kl}]=0$ and $[x_i, \BB_{kl}]=0$.
By  (\ref{xVsPermFml2}) :
\bea \label{PermLocFm111}
x_{j_1}^{M^t\BB} x_{j_2}^{M^t\BB}...x_{j_r}^{M^t\BB}=
\sum_{I=(i_1\le i_2\le i_3\le ...\le i_r\le n)} x_{i_1} x_{i_2}...
x_{i_r} perm_{norm}^{col}((M^t\BB)_{IJ}),
\eea
on the other hand, $x_{j}^{M^t\BB}= \sum_{k} x_{k}^{M^t} \BB_{kj}$,
and $x_{j}^{M^t}$ commute by Manin's property (proposition \pref{Coact-pr}), so  again by (\ref{xVsPermFml2}):
\bea
x_{j_1}^{M^t\BB} x_{j_2}^{M^t\BB}...x_{j_r}^{M^t\BB}=
\sum_{L=(l_1\le l_2\le l_3\le ...\le l_r), l_r\le m} x_{l_1}^{M^t} x_{l_2}^{M^t}...
x_{l_r}^{M^t} perm_{norm}^{col}((\BB)_{LJ})=
\eea
and again by (\ref{xVsPermFml2}):
\bea \label{PermLocFm222}
=
\sum_{L=(l_1\le l_2\le l_3\le ...\le l_r), l_r\le m~~}
\sum_{I=(i_1\le i_2\le i_3\le ...\le i_r\le n)} x_{i_1} x_{i_2}...
x_{i_r} perm_{norm}^{col}((M^t)_{IL})
perm_{norm}^{col}((\BB)_{LJ}).
\eea
Comparing (\ref{PermLocFm111})  and (\ref{PermLocFm222}) we come to the desired proposition. \BX

{\bf Condition 1.}
\bea \label{xCSS1}
\forall p,j~\exists x^Q_j \in \RRR \otimes \CC[x_1,...,x_n]:~~~~\sum_{l=1,...,m} x^M_l [ \BB_{lj}, x^M_p] = x^M_p x^Q_j,
\eea
note that $x^Q_j$ does not depend on $p$.

As usually  we denote by $Q$ an $n \times s$ matrix  corresponding to elements $x_i^Q$:
\bea \label{xQfmlLOC}
(x_1^Q,...,x_s^Q)=(x_1,...,x_n) Q
\eea

{\bf Condition 2.}
$x^M_i$ and  $x^Q_j$ commute:
\bea \label{xCSS2}
\forall i,j~:~~~~x^M_i x^Q_j = x^Q_j x^M_i .
\eea

{\Th \label{ThxCSS} Assume $M$ is $m\times n$ \MM, $\BB$ is $m\times s$ matrix  (not necessarily Manin),
and matrices $M^t$ and $\BB$ satisfy conditions (\ref{xCSS1}, \ref{xCSS2}) above.
Then
the following Cauchy-Binet formulae holds.
Let $I=(i_{1}\le i_2\le ...\le i_r)$,  $J=(j_{1},...,j_r)$,
\footnote{conditions  $j_{a}<j_b$, $j_a\ne j_b$ are {\em not}  required}
be two multi-indexes $i_a\le n, j_a\le s$, $r\le n,s$
then:
\bea \label{xCBfmlTh2}
perm_{norm}^{col}((M^t\BB)_{IJ}-Q_{IJ}~diag(r-1,r-2,...,1,0) )
= \sum_{L=(l_1\le l_2\le ...\le l_r), l_r\le m}  perm_{norm}^{col}((M^t)_{IL}) perm_{norm}^{col}(\BB_{LJ}),
\eea
here $A_{IJ}$ is  a matrix such that $(A_{IJ})_{ab}=A_{i_aj_b}$,
(see notation \pref{NotAIJ});
where $Q$ is a matrix corresponding to elements $x_i^Q$ by formula (\ref{xQfmlLOC})
and elements $x_i^Q$ arise in the condition 1 above;
by $diag(a_1,a_2,...)$ we denote diagonal matrix with $a_i$ on the diagonal.
}

Let us recall definition \pref{DefCSScond} that two matrices $M,\BB$  of sizes $n\times m$ and $m\times o$ respectively
  satisfy the
Caracciolo-Sportiello-Sokal condition (CSS-condition for brevity),
if the following is true:
\bea \label{copyCSSc1}
[M_{ij}, \BB_{kl} ] = - \delta_{jk} Q_{il},
\eea
for some elements $Q_{il}$.

Or "in words": elements  in $j$-th column of $M$ commute with elements in $k$-th row of $\BB$
unless $j=k$, and in this case commutator of the elements $M_{ik}$ and $\BB_{kl}$ depends only on $i,l$,
but does not depend on $k$.
(See \cite{CSS08} formula (1.14) page 4, our $\BB$ is {\em transpose} to their $B$).

{\Th \label{CBpermTheorCSS} Assume $n>1$.
Assume $M$ is $m\times n$ \MM, $\BB$ is $m\times s$ matrix  (not necessarily Manin),
and matrices $M^t$ and $\BB$ satisfy CSS-condition above, then
the following Cauchy-Binet formulae holds true.
Let $I=(i_{1}\le i_2\le ...\le i_r)$,  $J=(j_{1},...,j_r)$,
be two multi-indexes $i_a\le n, j_a\le s$, $r\le n,s$
then:
\bea
perm_{norm}^{col} ( (M^t\BB)_{IJ} - Q_{IJ}diag(r-1,r-2,...,1,0))= \sum_{L=(l_1\le l_2\le ...\le l_r), l_r\le m}
perm_{norm}^{col} ( (M^t)_{IL}) perm_{norm}^{col} (\BB_{LJ}),
\eea
where matrix $Q$ matrix with elements $Q_{ij}$.
}

Let us also give an analogue of lemma \pref{XreformulCSSlem2}:
{\Lem \label{XreformulCSSlem2} For $n>1$ CSS-condition automatically implies
$x_i^M$ and $x_j^Q$ commute (and hence that condition 2 page \pageref{xCSS2} is satisfied).
}

{\Rem ~} Theorem above holds true for $n=1$ under additional requirement $[M_{11}, [M_{11}, \BB_{1i}]]=0$.

Let us modify  TCSS-condition (\ref{TCSSc1}) for antisymmetric matrix $M$.
We impose the following condition:
\bea \label{copyTCSSc1}
[M_{ij}, \BB_{kl} ] = - h ( \delta_{jk} \delta_{il} - \delta_{ik} \delta_{jl}),
\eea
for some element $h$.

Or "in words": element $M_{ij}$ commute all elements of $\BB$
 except elements $\BB_{ik}$ and $\BB_{ki}$,
and these non-zero commutators are equal to $\pm h$ respectively, for some element $h$.

{\Th Assume $n>2$.
Assume $M$ is  $n\times n$ {\em antisymmetric} matrix with commuting entries (i.e. $\forall i,j,k,l:~~[M_{ij},M_{kl}]=0$)
\footnote{we can require $M$ to be a \MM, but for the field of characteristic not equal to $2$,
antisymmetric \MM~ is matrix with commuting entries },
$\BB$ is $n\times s$ matrix  (not necessarily Manin),
and matrices $M$ and $\BB$ satisfy TCSS-condition for antisymmetric matrices, then
the following Cauchy-Binet formulae holds true.
Let $I=(i_{1}\le i_2\le ...\le i_r)$,  $J=(j_{1},...,j_r)$,
be two multi-indexes $i_a\le n, j_a\le n$, $r\le n,s$
then:
\bea
perm_{norm}^{col} ( (M\BB)_{IJ} - h~ diag(r-1,r-2,...,1,0))= \sum_{L=(l_1\le l_2\le ...\le l_r)}
perm_{norm}^{col} ( M_{IL}) perm_{norm}^{col} (\BB_{LJ}).
\eea
}

An analogue of lemma \pref{LemHcommute} is the following:
{\Lem \label{xLemHcommute} Assume $n>2$. Assume $M$ is a antisymmetric $n\times n$ matrix with commuting entries (i.e. $\forall i,j,k,l:~~[M_{ij},M_{kl}]=0$),
 $M$ and $\BB$ satisfy TCSS-condition above, then:
\bea
\forall i,j:~~~~[M_{ij}, h]=0.
\eea
}

{\Rem ~} Theorem above holds true for $n=2$ also under additional requirement $[M_{12},h]=0$.

An analogue of lemma \pref{LemSymMatr} is the following:
{\Lem \label{xLemSymMatr} Obviously for any antisymmetric matrix $M$ (not necessarily \MM):
\bea
\sum_{l=1,...,n} x_l x^M_l =0.
\eea
}

\subsubsection{Toy model}

Let us provide some toy model for the identities above, which is actually a particular
case of the theorem \pref{CBpermTheorCSS} (the case $n=1$).
The proof is extremely simple and it is actually a good illustration
of the proof of the main theorem \pref{ThGrasCSS}.

{\Prop \label{ToyId} Assume elements $M$ and $\BB$ satisfy the following condition: $[M,[M,\BB]]=0$,
denote by $Q=[M,\BB]$,
then for any $r$:
\bea
(M\BB - (r-1)Q) (M\BB - (r-2)Q)... (M\BB - Q) (M\BB )=(M)^{r} (\BB)^{r}.
\eea
}

{\Ex ~}
\bea
(\partial_z z - (r-1)) (\partial_z z - (r-2))... (\partial_z z - 1) (\partial_z z )=\partial_z^{r}  z^{r}.
\eea

\PRF
For $r=1$ it is a tautology.
Consider the left hand side of the equality and by induction:
\bea
(M\BB - (r-1)Q) (M\BB - (r-2)Q)... (M\BB - Q) (M\BB )=(M\BB - (r-1)Q) (M^{r-1} \BB^{r-1})=
\eea
(formula above is analogous to (\ref{IndFml000}) in the proof of the main theorem)

\bea \label{fmlLLLLLLLLLLLL}
=
M M^{r-1} \BB \BB^{r-1}+ M [\BB, M^{r-1} ]  \BB^{r-1}  -            (r-1)Q (M^{r-1} \BB^{r-1})=
\eea
(formula above is analogous to (\ref{calcGrCss}) in the proof of the main theorem)

Since $[M,Q]=0$:
\bea
  [\BB, M^{r-1} ]  =            (r-1)Q M^{r-1} ,
\eea
(formula above is analogous to lemma \pref{AntiCommLem}  in the proof of the main theorem)

So we continue  (\ref{fmlLLLLLLLLLLLL}):
\bea
=
M M^{r-1} \BB \BB^{r-1}+ M (r-1)Q M^{r-1}  \BB^{r-1}  -            (r-1)Q (M^{r-1} \BB^{r-1})=
M M^{r-1} \BB \BB^{r-1}.
\eea
It is right hand side of the desired equality. Proposition is proved. \BX

{\Rem ~} Capelli identities related to permanents can be also  found
in
\href{http://arxiv.org/abs/0809.3516}{(S. Caracciolo, A. Sportiello, A. Sokal)}
\cite{CSS08} proposition 1.5  page 9 (due to  Turnbull),
\href{http://www.springerlink.com/content/rt5r3313732p48j7}{(M. Nazarov)}  \cite{Na91},
and
\href{http://lanl.arxiv.org/abs/q-alg/9602027}{(A.~Okounkov)}
\cite{OkB96}.
\footnote{A. Okounkov considers immanants i.e.
$\sum_{\sigma \in S_n} \chi(\sigma) A_{\sigma(1)1} A_{\sigma(2)2}... A_{\sigma(n)n}$,
where $\chi$ is some character of the symmetric group. Permanent and determinant
are particular cases of the immanant for $\chi=1$ and $\chi=(-1)^{sgn(\sigma)}$ respectively.
}
Our result is clearly different from the first mentioned result
and relations with the others are not clear.

\section{Further properties \label{PropSect}}
In this section we discuss other properties of commutative matrices
which can be extended to the case of \MMs. Some of them are new like
the multiplicativity property of the determinant, the relation of
the determinant and the Gauss decomposition, conjugation to the
second normal
(also called Frobenius e.g. \href{http://www.dleex.com/read/?7443}{(Wilkinson)} \cite{WBook} page 15)
form. Others can be already
found, somewhat scattered, in the literature. We include them in
order to provide a complete list of properties established at the
moment so far \MMs~ and to add some details, comments or different
proofs of these results.

\subsection{Cayley-Hamilton theorem and the second normal (Frobenius) form \label{CH-ss} }
The Cayley-Hamilton theorem
can be considered one of the basic results in linear algebra. It was
generalized in \cite{CF07} (theorem 3) to the case of \MMs. Let us
recall it and present some corollary about conjugation to the second
normal (Frobenius) form. Some bibliographic notes are
in section \pref{BibNotSec}.

{\Th \label{CH-th} Let $M$ be  $n\times n$ \MM. Denote by
$\sigma_i,$ $i=0,...,n$ the coefficients of  its characteristic
polynomial: $det^{column}(t-M)=\sum_{i=0,...,n} (-1)^i \sigma_i
t^{n-i}$, then: \bea
\sum_{i=0,...,n} (-1)^i \sigma_i M^{n-i}=0, \mbox{ ~i.e. ~ } det^{column}(t-M)|_{t=M}^{right~substitute}=0.
\eea
If $M^t$ is a \MM, then one can obtain a similar result, using {\em left} substitution and the
row determinant:\\
$det^{row}(t-M)|_{t=M}^{left~substitute}=0$.
}

{\Rem ~} In the commutative case $\sigma_i$ is the $i$-th elementary
symmetric function of the eigenvalues ($\sigma_i=\sum_{1\le
j_1<...<j_i\le n} \lambda_{j_1}\lambda_{j_2}...\lambda_{j_i}$). In
general: $\sigma_1=Tr(M), \sigma_n=det(M)$, $\sigma_k=Tr \Lambda^k
M$. \bea
det^{column}(t-M)&=&\sum_{i=0,...,n} (-1)^i \sigma_i t^{n-i}=  \\
&=&
t^n+t^{n-1}(-1) \sigma_1+t^{n-2}(+1) \sigma_2+...
+t(-1)^{n-1}\sigma_{n-1} +(-1)^{n}\sigma_{n}.
\eea

\PRF  Proposition \ref{ad-inv} shows that $t-M$ admits a classical
(left) adjoint matrix $(M-t~1_{n\times n})^{adj}$, such that
 \bea (M-t~1_{n\times n})^{adj}
(M-t~1_{n\times n})=det(M-t~1_{n\times n}) 1_{n\times n}, \eea
where, as usual, we denote by  $1_{n\times n}$ the identity matrix
of size $n$. The standard idea of proof is very simple: we want to
substitute $M$ where $t$ stands; the LHS of this equality vanishes
manifestly, hence we obtain the desired equality $det(M-t~1_{n\times
n})|_{t=M}=0$. The only issue we need to clarify is {\em how} to
substitute $M$ into the equation and why the substitution preserves
the equality.

Let us denote by $Adj_k(M)$ the matrices defined by:
$\sum_{k=0,...,n-1} Adj_k(M) t^k= (M-t~1_{n\times n})^{adj}$. The equality
above is an equality of polynomials in the variable $t$:
\bea
(\sum_k Adj_k(M) t^k) (M-t ~ 1_{n\times n})= \sum_k Adj_k(M) M  t^k - \sum_k Adj_k(M)   t^{k+1} =
\\ = det (M-t~1_{n\times n})= (-1)^n \sum_{i=0,...,n} (-1)^i \sigma_i t^{n-i}.
\eea
This means that the coefficients
of $t^i$ of both sides of the relation coincide. Hence we can
substitute $t=M$ in the equality, substituting "from the right":
\bea
\sum_k Adj_k M  M^k - \sum_k Adj_k M^{k+1} =
 (-1)^n \sum_{i=0,...,n} (-1)^i \sigma_i t^{n-i} \big\vert_{t=M}.
\eea The left hand side is manifestly zero, so we obtain the desired
equality: $\sum_{i=0,...,n} (-1)^i \sigma_i M^{n-i}=0.$ \hfill\BX



Let us present a corollary on conjugation of a \MM to the so-called
Frobenius normal form (e.g. \href{http://www.dleex.com/read/?7443}{(Wilkinson)} \cite{WBook} page 15), also called the second normal form (e.g. \cite{Kir}).

{\Cor \label{FrobCor}
 Let $M$ be a  $n\times n$ \MM~ with elements in an associative
ring \kol, and let $\sigma_i,$ $i=0,...,n$ be the coefficients of
its characteristic polynomial, that is,
$det^{column}(t-M)=\sum_{i=0,...,n} (-1)^i \sigma_i t^{n-i}$.

Let $v=(v_1,...,v_n)$ be a vector with elements in \kol~ such that
$\forall ~k,l$, $v_k$ commutes with  $\sigma_l$, (this happens, for
instance, if $v_k\in \CC$). Let, finally,  \bea \label{C}  D= \left(
\begin{array}{c}
v \\
v M \\
v M^{2}\\
\ldots\\
v M^{n-1}
\end{array}\right)
=
\left(
\begin{array}{cccc}
v_1&  ... & v_n \\
\sum_i v_i M_{i1}& ... & \sum_i v_i M_{in}   \\
\sum_{i,j} v_i M_{ij}M_{j1}& ... & \sum_{i,j} v_i M_{ij}M_{jn}   \\
\ldots & ... & ...\\
\ldots & ... & ... 
\end{array}\right).
\eea Then it holds \bea \label{hyp} D ~  M  = M_{Frob} ~ D,
\mbox{~~where~~}
M_{Frob}=
\left(
\begin{array}{cccccc}
0 & 1 & 0 & 0 &... &0 \\
0 & 0 & 1 & 0 &... &0 \\
... & ... & ... & ...&... &... \\
0 & 0 & 0 & ... & 0 &1 \\
(-1)^{n+1}\sigma_n & (-1)^{n}\sigma_{n-1} & (-1)^{n-1}\sigma_{n-2} & ... & -\sigma_2 & \sigma_1
\end{array}\right) .
\eea }
{\Rem ~} In words, this Corollary says that, under the
commutativity conditions $[v_k, \sigma_j]=0$, a Manin matrix $M$ can
be conjugated to its normal Frobenius form. Before giving the proof
of the corollary let us remark the following. It is easy to see that
an arbitrary matrix over a noncommutative ring \kol~ which is
embeddable in a noncommutative field can be conjugated in the form
above. (Indeed, this is equivalent to the fact that $n+1$ vectors in
a $n$ dimensional vector space are linearly dependent over an
arbitrary field (no need of commutativity) and applying this fact to
the vectors $v, Mv, M^2,...,M^{n}v$ one gets the claim. The
coefficients of linear dependence precisely appear in the last row
of the matrix $M_{Frob}$.)
However, in general the coefficients of the linear dependence will
depend on the vector $v$ and they are rational (not polynomial)
functions of the matrix entries. (I.e. they belong  to the  field of
fractions of \kol, but not to the original ring.)
We notice that for \MMs, the theorem holds in the same form as in
the case of ordinary matrices.


\PRF By definition, the $i$-th ($i=1,...,n$) row  $DM$ equals
$vM^i$. The same is true for the first $n-1$ rows of$M_{Frob} D$
equals to $vM^i$. One only needs to check the equality of the $n$-th
row. It equals to $vM^n$ for $DM$ and the Cayley-Hamilton theorem,
together with the  condition $[v_k,\sigma_l]=0$ precisely provide
that the same expression appears  in $M_{Frob} D$. Indeed the $n$-th
row of $M_{Frob} D$ equals $\sum_{l=0,...,n-1}
(-1)^{n-l+1}\sigma_{n-l} vM^{l}$; thanks to the commutativity
condition $[v_k,\sigma_l]=0$ we can rewrite it as $v
(\sum_{l=0,...,n-1} (-1)^{n-l+1}\sigma_{n-l} M^{l})$. By the
Cayley-Hamilton theorem this equals to $(v M^n)$. The corollary is
proved. \BX


{\Ex ~ \label{ExFrobForm2} }

Consider the $2\times 2$ case, $v=(0,1)$,
denote $M= \left(
\begin{array}{cc}
a & b \\
c & d
\end{array}\right)$.

\bea
D=\left(
\begin{array}{cccc}
0 & 1 \\
c & d
\end{array}\right), ~~~~~
M_{Frob}=\left(
\begin{array}{cccc}
0& 1 \\
-(ad-cb) & a+d
\end{array}\right).
\eea

\bea
D~ M=
\left(
\begin{array}{cccc}
0& 1 \\
c & d
\end{array}\right)
~
\left(
\begin{array}{cccc}
a & b \\
c & d
\end{array}\right)
=
\left(
\begin{array}{cccc}
c & d \\
ca+ dc  & cb+ d^2
\end{array}\right).
\eea

\bea
M_{Frob} D=\left(
\begin{array}{cccc}
0& 1 \\
-(ad-cb) & a+d
\end{array}\right)
\left(
\begin{array}{cccc}
0 & 1 \\
c & d
\end{array}\right)
=
\left(
\begin{array}{cccc}
c & d \\
ac+dc & -(ad-cb)+a d +d^2
\end{array}\right)
= \\
=
\left(
\begin{array}{cccc}
c & d \\
ac+dc & cb+ d^2
\end{array}\right)
= \mbox{ ~~ (use: ac=ca) ~~ } =
\left(
\begin{array}{cccc}
c & d \\
ca+dc & cb+ d^2
\end{array}\right).
\eea

\subsection{Newton and MacMahon-Wronski identities
\label{NMWss} } The aim of this section is to discuss a
generalization of the Newton and MacMahon-Wronski identities to the
case of \MMs. As we shall see, they hold true exactly in the same
form as in the commutative case, which we herewith recall.

There are three basic families of symmetric functions in $n$
variables: \begin{enumerate}
\item
$\esf_{k}=\sum_{1\le i_1<i_2<...<i_k\le n } \prod_{p=1,...,n}
\lambda_{i_p} , i=1,...,n$; they are called the {\em elementary}
symmetric functions. \item $\csf_k=\sum_{0\le i_1,...,i_n :
i_1+...+i_n=k} \prod_{p=1,...,n} \lambda_{p}^{i_p}$, $k>0$ - the
so-called {\em complete}.
\item
$\psf_k=\sum_{p=1,...,n} \lambda_{p}^{k}$, $k>0$ - the so-called
power sums.
\end{enumerate}
They can be rewritten in the matrix language as follows: \bea
&&  \esf_{k}=Tr \Lambda^k M, ~~~~\left( \sum_{k=0,...,n} (-t)^k \esf_{k}=det(1-t M)\right),\\
&& \csf_{k}= Tr S^k M, \\ && \psf_{k}= Tr (M^k). \eea Here $M$ is a
matrix with entries in $\CC$ with eigenvalues
$\lambda_1,...,\lambda_n$. $S^k M$ is the symmetric tensor power of
$M$ while  $\Lambda^k M$ is the antisymmetric power.

{\bf Theorem.} {\em  $\esf_k,\csf_k,\psf_k$ are related by the
following set of identities for all $k>0$: \bea \label{MNWfml1}
\mbox{MacMahon-Wronski:} &&
0=\sum_{l=0,...,k} (-1)^{l} \csf_l  \esf_{k-l}= \sum_{l=0,...,k} (-1)^{l}   \esf_{k-l} \csf_l,\\
\mbox{Newton:}&& -(-1)^k k    \esf_k =\sum_{i= 0,...,k-1}  (-1)^{i} \esf_{i} \psf_{k-i} ,\\
\mbox{Second Newton:}&&  k    \csf_k =\sum_{i= 0,...,k-1} \psf_{k-i}
\csf_{i}. \eea } Our main goal is to explain the following:
\\ {\bf
Claim.} {\em The formulas above hold true when $M$ is a \MM.}

{\Rem ~} In the case of \MMs~ the order in products in the formulas
above is important. The MacMahon-Wronski identity has been first
obtained in \href{http://arxiv.org/abs/math.QA/0303319}{(S.
Garoufalidis, T. Le, D. Zeilberger)} \cite{GLZ}, the Newton one in
\cite{CF07}. Here we will collect these results.

The identities above can be easily reformulated in terms of generating functions:

{\bf Corollary.} {\em Let $M$ be a \MM.
Denote by $\esft,\csft , \psft$ the following generating functions:
\bea
\esft & = & \sum_{k=0,...,n} (-t)^k \esf_{k}=det(1-t M),\\
\csft  & = & \sum_{k=0,...,\infty } t^k \csf_{k}= \sum_{k=0,...,\infty } t^k Tr S^{k}M,\\
\psft  & = & \sum_{k=0,...,\infty } t^k \psf_{k+1}= \sum_{k=0,...,\infty } t^k Tr (M^{k+1})  =   Tr \frac{M}{1-tM}.
\eea
Then the relations between $\esf_l, \csf_l, \psf_l$ can be written as follows:
\bea
\mbox{MacMahon-Wronski:} &&
1= \esft \csft= \csft \esft ,\\
\mbox{Newton:}&& -\partial_t \esft = \esft \psft,\\
\mbox{Second Newton:}&&   \partial_t \csft = \psft \csft .
\eea
}

{\Rem~} In the commutative case the Newton identities can be
reformulated in the other forms. This might not be the case for
\MMs, since, as we shall see in the sequel, for generic \MMs\  we
have: \[det ( e^M )\ne  e ^{Tr(M)},\text{  and  }  (det (1+M)) \ne
e^{Tr (ln (1+M))}.\]

\subsubsection{Newton identities}
%

Let us recall the result of \cite{CF07} giving a detailed proof. Some bibliographic notes are
in section \pref{BibNotSec}.

{\Th \label{NewtonTh} Let $M$ be  $n\times n$ \MM.
Denote by $\esf_i,$ $i=0,...,n$ coefficients of  its characteristic polynomial:
$det^{column}(1-tM)=\sum_{i=0,...,n} (-t)^i \esf_i $, conventionally, let
$\esf_k =0$, for $k>n$; {\bf then:}
\bea \label{Newt-fml1}
 - \partial_t det^{column}(1-tM) = (det^{column}(1-tM))
 \sum_{k=0,...,\infty} t^k Tr (M^{k+1}),\\
\label{Newt-fml2}
\Leftrightarrow \forall k \ge 0:
-(-1)^k k    \esf_k =\sum_{i= 0,...,k-1}  (-1)^{i} \esf_{i} Tr (M^{k-i}).
\eea
If $M^t$
is a \MM, then:
$\partial_t det^{row}(1-tM) =
\left( 1/t \sum_{k=0,...,\infty} t^k Tr (M^{k+1}) \right) (det^{row} (1-tM))$.
}

{\Rem ~} Using the generating functions $\esft =det(1-t M), \psft =
Tr \frac{M}{1-tM}$, one rewrites: \bea
 -\partial_t \esft = \esft \psft.
\eea
 One can also rewrite as: $ \partial_t det^{column}(t-M) =\frac {1}{t} (det^{column}(t-M))
 \sum_{k=0,...,\infty} Tr (M/t)^k$.

{\Rem ~} In the commutative case $\esf_i$ is $i$-th elementary
symmetric function of the eigenvalues ($\sum_{j_1<...<j_i} \prod
\lambda_{j_k}$). In general: $\esf_1=Tr(M), \esf_n=det(M)$,
$\esf_k=Tr \Lambda^k M$. \bea
det^{column}(1-tM)&=&\sum_{i=0,...,n} (-t)^i \esf_i =  \\
&=&
1+(-1) t \esf_1+(+1) t^{2} \esf_2+...
+(-1)^{n-1}t^{n-1}\esf_{n-1} +(-1)^{n} t^{n}\esf_{n}.
\eea


Formulas are identical in the Manin case, provided one pays
attention to the order of terms: $\esf_i Tr M^p$ if $M$ is a \MM\
($Tr M^p \esf_i $ if $M^t$ is a \MM).

\PRF The proof of the Theorem is somewhat a standard one. First we
need the following simple lemma:
{\Lem \label{TrAdjLem} Consider an arbitrary matrix $M$ (i.e., not
necessarily a \MM), and define its adjoint matrix $M^{adj}$  in the
standard way as follows:
 $M^{adj}_{kl}=(-1)^{k+l}det^{column}(\widehat{M}_{l k})$ where
$\widehat{M}_{l k}$ is the $(n-1)\times (n-1)$ submatrix of $M$ obtained
removing the l-th row and the k-th column.
{\bf Then:}
\[
Tr (t+M)^{adj}=\partial_t det^{column}(t+M).
\]
}

{\bf Proof of the Lemma.}  The proof is trivial, but let us nevertheless write it up. \\
\bea
\partial_t det^{column}(t+M) =
 \sum_{i=0,...,n-1} (n-i) t^{n-1-i} \sum_{principal\ i\times i\ submatrices\ K\subset M } det^{column}( K),
\eea where a principal submatrix is the one formed by the elements
obtained on the intersection of the rows and columns labeled by the
same set of indices $j_1,...,j_i$.
\bea
 Tr (t+M)^{adj}=\sum_{l=1,...,n} det^{column} \widehat{(t+M)}_{ll} = \\
=\sum_{l=1,...,n}  \sum_{i=0,...,n-1}  t^{n-1-i}  \sum_{principal~i\times i ~submatrices~K~of~\widehat{M}_{ll} }
 det^{column}(K).
\eea

Clearly, any principal submatrix of $\widehat{M}_{ll}$ is a
principal submatrix of $M$, so  one has the same terms in both sums.
Moreover submatrices of size $i$ appear as coefficients of
$t^{n-1-i}$ in both sums.

So we only need to observe that $det^{column}(K)$ for an $i\times i$
submatrix $K$ enters with the same multiplicity in both sums. In the
first sum the multiplicity is manifestly $n-i$. Let us look at the
second sum. The principal submatrix of size $i$ clearly is a
submatrix of $(n-i)$ principal submatrices $ \widehat{(t+M)}_{ll} $,
for example  submatrix of size 1, say $M_{11}$, is a submatrix of $
\widehat{(t+M)}_{22}, \widehat{(t+M)}_{33},
...,\widehat{(t+M)}_{nn}$. So we get that the desired coefficient is
$n-i$ and the lemma is proved. \BX

{\Rem ~} The same arguments can be applied when dealing with the
row-determinant, as well, symmetrized determinant, and so on and so
forth.

Let us finalize the proof of the theorem. We have the following
chain of relations: \bea \label{Nfml1} 1/t \sum_{k=0,...,\infty}
Tr\big((-M/t)^k\big) = Tr\frac{1}{t+M} \stackrel{Cramer}{=}Tr
\Bigl((det^{column}(t+M))^{-1}(t+M)^{adj} \Bigr) =
\nn\\
=
(det^{column}(t+M))^{-1} Tr (t+M)^{adj}\stackrel{~~Lemma~ \ref{TrAdjLem}~~}{=}
(det^{column}(t+M))^{-1}
\partial_t det^{column}(t+M).
\eea
This identity gives the identities (\ref{Newt-fml2}), i.e., \bea
 \forall k \ge 0:
-(-1)^k k    \esf_k =\sum_{i= 0,...,k-1}  (-1)^{i} \esf_{i} Tr
(M^{k-i}), \eea which in turn is equivalent \footnote{In the
commutative case one can pass from (\ref{Nfml1}) to (\ref{Nfml2}) as
follows.
Substitute $M=N^{-1}$ in (\ref{Nfml1}) and use that $det(t+N^{-1}) =
det(Nt+1)det(N^{-1})$, one gets:
$Tr\frac{N}{Nt+1}=(det^{column}(Nt+1))^{-1} det(N) det(N^{-1})
\partial_t det^{column}(Nt+1)= (det^{column}(Nt+1))^{-1} \partial_t det^{column}(Nt+1) $
changing $N$ to $-N$ one gets (\ref{Nfml2}).  The open question is
whether $det(t+N^{-1}) = det(Nt+1)det(N^{-1})$ for \MMs\ as well. }
to formula (\ref{Newt-fml1}): \bea \label{Nfml2}
 - \partial_t det^{column}(1-tM) = (det^{column}(1-tM))
 \sum_{k=0,...,\infty} t^k Tr (M^{k+1}).
\eea Theorem \ref{NewtonTh} on the Newton identities is thus proved.
\hfill \BX.

{\Rem ~} The case  in which $M^t$ is a \MM~ can be treated in a
similar way.
 {\Ex
~\label{Newt22} } It is instructive to explicitly perform the
computation in the  $2\times 2$ case. Let
\[
M=\left(\begin{array}{cccc}
a & b \\
c & d \\
\end{array}\right),\] where $M$ is a generic matrix (i.e., not necessarily a \MM).
We have: \bea
Tr(M^2)-\esf_1 Tr(M)+2\esf_2= \\
Tr
\left(\begin{array}{cccc}
a^2+bc & ab+bd \\
ca+dc & cb+d^2 \\
\end{array}\right)
-(a+d)^2+2(ad-cb)
= \\ =
(a^2+bc+cb+d^2)-(a^2+ad+da+d^2) + 2(ad-cb)
= \\ =
(bc)-(da)+ (ad-cb)
= [a,d]+[b,c], \eea and, similarly, \[
 Tr(M^3)-\esf_1
Tr(M^2)+\esf_2 Tr(M) = ([a,d] + [b,c ])a+[c,a]b+[b,d]c. \]
Notice
that $-3\esf_3$ does not appear in the last formula since
$\esf_k=0,$ for $k>2$ for $2\times 2$ matrices. We see that Manin's
relations imply that the expressions above are zeroes.

\paragraph{A No-go fact \label{nogoNewtsect}}

In the commutative case there is the well-known identity $det (e^N)
= e ^{TrN}$, which can be readily seen by diagonalization of $N$.
Substituting $N=log(1-tM)$ one obtains: $ (det (1-tM)) =  e^{Tr (ln
(1-tM))}$. Here we show that these identities do not hold in the
case of \MMs; actually, they do not hold even in more restrictive
case of Cartier-Foata matrices.

Let us remark that in the commutative case the Newton identity
\[-\partial_t det^{column}(1-tM) =(det^{column}(1-tM))
 \sum_{k=0,...,\infty} Tr t(Mt)^k\]
 easily follows from the identities above.
Indeed, deriving the identity $ (det (1-tM)) =  e^{Tr (ln (1-tM))}$,
with respect to $t$, one obtains:
\begin{equation}\label{ngfact}
\partial_t (det (1-tM)) = Tr
(ln'(1-tM)) e^{Tr (ln (1-tM))}=Tr (ln'(1-tM)) (det (1-tM)).
\end{equation}
Using $Tr(ln'(1-tM))= - Tr (M(1-tM)^{-1})$, one arrives to
$\partial_t (det (1-tM)) = - Tr (M (t-M)^{-1} ) (det (1-tM))$, which
is the Newton identity from theorem \ref{NewtonTh} above.

One may try to argue in the opposite direction, but the crucial
point is the commutativity which is absent for \MMs, and was used in
the first equality of (\ref{ngfact}). So it is not guaranteed that
the exponential relations hold true in the general case of a \MM, as
indeed it is proven by the next two counterexamples.

{\CEx \label{CEx1} ~} Consider a $2\times 2$ matrix $M$, \bea
M=\left(\begin{array}{cccc}
a & b \\
c & d \\
\end{array}\right),
\eea and assume it is a Cartier-Foata matrix, that is, elements from
different rows commute. Introduce a formal scalar variable $\epsilon$; 
clearly enough, $1_{2\times 2}+ \epsilon M$ is again a Cartier-Foata
matrix.

{\bf Fact}.
 \bea det(exp(1+\epsilon M)) \ne
exp(Tr(ln(1+\epsilon M))).
\eea Actually, the equality holds up to order 2, but not at order
3.\\ \PRF The coefficient of $\epsilon^3$ in the left hand side is
equal to: \bea
\frac{1}{6}(a^3+ bca+abc+bdc+cab+dcb+cbd+d^3+ \\
+ 3acb+3ad^2+3 a^2d+3bcd-3cab-3cbd-3cab-3dcb), \eea  while that on
the right hand side is $\frac{1}{6}(a+d)^3$. In the Cartier-Foata
case $[a,d]=0$, so the difference is given by: \bea
\frac{1}{6}( bca+abc+bdc+cab+dcb+cbd+ \\
+ 3acb+3bcd-3cab-3cbd-3cab-3dcb) = \\
=\frac{1}{6}( bac+bcd- abc-bdc) =\frac{1}{6}( [b,a] c+b[c,d]). \eea
Since pairs of elements $a,b$ and $c,d$ may be taken to generate
free associative algebras and arranged to give a Cartier-Foata
matrix, we see that this does not vanish. \BX

{\CEx ~} Consider $M$, $\epsilon$ as above; then \bea det(1+\epsilon
M) \ne e^{Tr(ln(1+\epsilon M))}. \eea The equality holds up to order
2.

The proof of this fact goes exactly as that of Counterexample 1.




\subsubsection{MacMahon-Wronski relations \label{MacMahSect}  }
For the sake of completeness  let us briefly discuss the so-called
MacMahon-Wronski formula for \MMs. It was first obtained in
\cite{GLZ}. In the
language of symmetric functions this identity relates the elementary
and the complete symmetric functions.
We will provide some more detailed bibliographic notes
in section \pref{BibNotSec}.

For an $n\times n$-matrix $M$  over a commutative ring  the MacMahon
- Wronski
identity reads \bea \label{MacMfml} 1/det(1-M)=\sum_{k=0,...,\infty}
Tr {S^k M}, \eea where $S^k M$ is $k$-th symmetric power of $M$. It
can be easily verified by diagonalizing the matrix $M$. As it was
mentioned before, $Tr {S^k M}$ is the {\em complete} symmetric
function of the eigenvalues $\lambda_i$ of $M$, i.e.
$\csf_k=\sum_{0\le i_1,...,i_n : i_1+...+i_n=k} \prod_{p=1,...,n}
\lambda_{p}^{i_p}$.\\ {\em {\bf Theorem}  \cite{GLZ}
The MacMahon-Wronski identity  holds true for \MMs~
\footnote{ one can find in \cite{GLZ} more general case of q-\MMs}
provided one defines} 
\bea \label{Skfml} Tr S^k M =
1/k!\sum_{l_1,...l_k: 1\le l_i \le n} perm^{row}
\left(\begin{array}{cccc}
M_{l_1\,l_1} & M_{l_1\,l
_2}& ... & M_{l_1\,l_k} \\
M_{l_2\,l_1} & M_{l_2\,l_2}& ... & M_{l_2\,l_k} \\
... & ... & ... & ... \\
M_{l_k\,l_1} & M_{l_k\,l_2}& ... & M_{l_k\,l_k}
\end{array}\right).
\eea
\noindent We remark that the range of summation in this formula
allows repeated indexes: $l_{i_1} = l_{i_3} =...$. The permanent of
a \MM~ was defined by  formula \pref{M-perm-f} as follows: \bea perm
M=perm^{row} M=\sum_{\sigma\in S_n} \prod_{i=1,...,n}
M_{i\,\sigma(i)}. \eea Due to the property that the permanent of a
\MM~ does not change under any permutation of columns, one can
rewrite the formula above with the summation without repeated
indexes: \bea Tr S^k M = 1/k!\sum_{1\le l_1\le...\le l_k \le n} n_1!
n_2! ...n_l! perm^{row} \left(\begin{array}{cccc}
M_{l_1\,l_1} & M_{l_1\,l_2}& ... & M_{l_1\,l_k} \\
M_{l_2\,l_1} & M_{l_2\,l_2}& ... & M_{l_2\,l_k} \\
... & ... & ... & ... \\
M_{l_k\,l_1} & M_{l_k\,l_2}& ... & M_{l_k\,l_k}
\end{array}\right),
\eea where $n_i$ is the set of multiplicities of the set
$(l_1,...,l_k)$, i.e. $n_i$ is multiplicity of the number $i$ in the
set  $(l_1,...,l_k)$.

One sees that the definition of traces of symmetric powers
is the same as in the commutative case, with the proviso in mind to
use row permanents.

As an immediate consequence of formula \ref{MacMfml} we get
{\Cor
 For all $p>0$ it holds
\bea 0=\sum_{l=0,...,p} (-1)^{l} \csf_l  \esf_{p-l}=
\sum_{l=0,...,p} (-1)^{l}   \esf_{p-l} \csf_l, \eea where $\csf_k=Tr
S^k M$ and $\esf_k$ are the coefficients of the
 characteristic polynomial $det(1-tM)=\sum_{l=0,...,n} (-1)^l t^l
 \esf_l$.}

{\Ex ~} Let us explicitly write the relation:
$\csf_2-\csf_1\esf_1 +\esf_2 =0, $
on a  $2\times 2$ matrix $ M=\left(\begin{array}{cccc}
a & b \\
c & d \\
\end{array}\right).$

We have: \bea \csf_1= a+d, ~~~~~ \esf_1=a+d, ~~~~~ \esf_2=ad-cb,
\eea while
\bea \csf_2= \frac{1}{2}( perm \left(\begin{array}{cccc}
a & a \\
a & a \\
\end{array}\right)
+
perm
\left(\begin{array}{cccc}
a & b \\
c & d \\
\end{array}\right)
+
perm
\left(\begin{array}{cccc}
d & c \\
b & a \\
\end{array}\right)
+
perm
\left(\begin{array}{cccc}
d & d \\
d & d \\
\end{array}\right))=\\
= \frac{1}{2}( 2a^2+ ad+bc+da+cb+2d^2)=
\eea
using Manin's
relation: $[a,d]=[c,b]$:
\bea
 = ( a^2+ ad+bc +d^2).
\eea
 Thus: \bea
\csf_2-\csf_1\esf_1 +\esf_2  = ( a^2+ ad+bc +d^2) -(a+d)^2+ad-cb =
\\ = ( a^2+ ad+bc +d^2) -(a^2+ad+da+d^2)+ad-cb =  ( bc )
-(da)+ad-cb =   [a,d] + [b,c ]= 0. \eea


\subsubsection{Second Newton identities \label{SymNewtSect}   }

{\Th Let $M$ be a \MM. Let  $\csf_k=Tr S^k M$ be defined by \pref{Skfml},
 $\csft=\sum_{k=0,...,\infty } t^k \csf_{k}$, and
$\psft= \sum_{k=0,...,\infty } t^k Tr (M^{k+1})  =   Tr \frac{M}{1-tM}$.
Then the following identities hold:
\[
\partial_t \csft = \psft \csft,
\Leftrightarrow \forall k>0:\>  k    \csf_k =\sum_{i= 0,...,k-1}    Tr (M^{k-i}) \csf_{i}.
\]
If $M^t$
is a \MM, similar formulas hold with the reverse
$\partial_t S(t) = S(t) T(t)$, and the use of column-permanents in formula \pref{Skfml}.
}

\PRF The proof follows from the Newton and MacMahon-Wronski identities above. Indeed,
with the definition $\esft=det(1-tM)$, the MacMahon-Wronski identities read
\[
\esft=\csft^{-1},\]
while the Newton identities read $-\partial_t \esft=\esft \psft$. Substituting $ -\partial_t \esft=
-\partial_t \csft^{-1} = + \csft^{-1} (\partial_t \csft) \csft^{-1}$, we get
\[ \csft^{-1} (\partial_t \csft) \csft^{-1} = \csft^{-1} T(t),\] and hence
$ \partial_t \csft = T(t) \csft $.\BX

\subsection{Pl\"{u}cker relations \label{Pluck-ss} }
Here we recall the
simplest version of the Pl\"{u}cker identities for \MMs~ (see
\cite{Manin}). They actually follow immediately from the coaction
characterization of \MMs~ (proposition \ref{Coact-pr}, page
\pageref{Coact-pr}).
We will provide some bibliographic notes on various noncommutative Pl\"{u}cker coordinates
in section \pref{BibNotSec}.

{\Prop \label{Plucker} Consider a $2\times 4$  matrix $A$, assume
$A^t$ is a \MM,
 let $\pi_{ij}$ be minors made from $i$-th and $j$-th columns (minors are calculated as row-determinants, since
$A^t$ is a \MM, not $A$):
\bea
A=\left(\begin{array}{cccc}
 a_1 & a_2& a_3 & a_4 \\
 a_1' & a_2' & a_3' & a_4'
\end{array}\right) \\
\mbox{Then: }
(\pi_{12} \pi_{34} +  \pi_{34} \pi_{12} )- (\pi_{13} \pi_{24} +  \pi_{24} \pi_{13} )
+(\pi_{14} \pi_{23} +  \pi_{23} \pi_{14} )=0
\eea}

\PRF The proof is the same as in the commutative case. Consider the
Grassman algebra $\CC[\psi_1,....,\psi_4]$, and the  variables
$\tpsi_1,\tpsi_2$ defined as: \bea \left(\begin{array}{c}
\tpsi_1 \\
\tpsi_2
\end{array}\right) = A
\left(\begin{array}{c}
\psi_1 \\
\psi_2 \\
\psi_3 \\
\psi_4 \\
\end{array}\right).
\eea
It is clear that $ \tpsi_1 \wedge \tpsi_2= \sum_{i<j} \pi_{ij} \psi_i \wedge \psi_j$
By proposition \ref{Coact-pr} page \pageref{Coact-pr} ~~
$\tpsi_1,\tpsi_2$ are again Grassman variables,
 so $ (\tpsi_1\wedge \tpsi_2)^2=0$.
Writing this equation explicitly one arrives to the Pl\"{u}cker relations.
\BX

\subsection{Gauss decomposition and the determinant \label{Gaus-ss}}
Here we show that the determinant of a \MM~ can be expressed via the
diagonal part of the Gauss decomposition exactly in the same way as
in the commutative case.

{\Prop \label{Gauss-pr} Let $M$ be a \MM, assume it can be
factorized into Gauss form: \bea M=\left(\begin{array}{ccc}
 1 &{}&x_{\alpha \beta}\\
{}&\ddots&{}\\
0&{}&1
\end{array}\right)
\left(\begin{array}{ccc}
y_1&{}&0\\
{}&\ddots&{}\\
0&{}&y_n
\end{array}\right)
\left(\begin{array}{ccc}
1&{}&0\\
{}&\ddots&{}\\
z_{\beta \alpha}&{}&1
\end{array}\right),
\eea
then
\bea
det(M)= y_n ... y_1.
\eea
Analogously for
\bea
M=\left(\begin{array}{ccc}
 1 &{}&0\\
{}&\ddots&{}\\
x_{\alpha \beta}' &{}&1
\end{array}\right)
\left(\begin{array}{ccc}
y_1'&{}&0\\
{}&\ddots&{}\\
0&{}&y_n'
\end{array}\right)
\left(\begin{array}{ccc}
1&{}&z_{\beta \alpha}'\\
{}&\ddots&{}\\
0&{}&1
\end{array}\right),
\eea
it is true:
\bea
det(M)= y_1' ... y_n'.
\eea

}

{\Ex ~}
\bea
&& M= \left(\begin{array}{ccc}
a  & b \\
c  & d \\
\end{array}\right)
=
\left(\begin{array}{ccc}
1  & bd^{-1} \\
0  & 1 \\
\end{array}\right)
\left(\begin{array}{ccc}
a-bd^{-1}c  & 0 \\
0  & d \\
\end{array}\right)
\left(\begin{array}{ccc}
1  & 0 \\
 d^{-1} c   & 1 \\
\end{array}\right),  \\
&& det(M)=ad-cb=d(a-bd^{-1} c)
,\\
&& M= \left(\begin{array}{ccc}
a  & b \\
c  & d \\
\end{array}\right)
=
\left(\begin{array}{ccc}
1  & 0 \\
ca^{-1}  & 1 \\
\end{array}\right)
\left(\begin{array}{ccc}
a  & 0 \\
0  & d- c a^{-1} b \\
\end{array}\right)
\left(\begin{array}{ccc}
1  & a^{-1} b \\
 0   & 1 \\
\end{array}\right),  \\
&& det(M)=ad-cb=a(d-ca^{-1} b).
\eea

\PRF   According to \href{http://arxiv.org/abs/q-alg/9705026}
{I. Gelfand, V. Retakh} \cite{GR97}
page 14, theorem 2.2.5 
for any noncommutative $M$
it is true that $y_k=|M_{(k)}|_{kk}$, where
$M_{(k)}$ is submatrix of $M$ with $k,k+1,k+2,...,n$ rows and columns,
and $|N|_{kk}$ is  quasideterminant, which is by definition  inverse to a corresponding element
of $N^{-1}$. By Cramer's formula (see proposition  \ref{ad-inv}, page \pageref{ad-inv})
$|M_{(k)}|_{kk}=(det^{col}(M_{(k)})^{-1} det^{col}(M_{(k+1)}))^{-1} $,
hence straightforward multiplication and cancellation gives the desired result.

The same arguments for the second equality (see
\href{http://arxiv.org/abs/math.QA/0208146} {I. Gelfand, S. Gelfand,
V. Retakh, R. Wilson} \cite{GGRW02} page 41 theorem 4.9.7). \BX
{\Rem ~} Though we do not discuss applications to integrability in
this paper, but let us mention that the result applied to the
Yangian generating matrix $e^{-\p}T(z)$ (which is a \MM) gives the
useful fact relating the $qdet(T(z))$ and the diagonal of the Gauss
decomposition.

\subsection{Bibliographic notes \label{BibNotSec} }

{\bf The Cayley-Hamilton theorem.}
The Cayley-Hamilton-like theorems for noncommutative matrices
are an object of numerous papers. Enough to say that "Cayley-Hamilton" occurs 179 times during the MathSciNet search.
There are also some similar theorems in non-matrix settings.
\href{http://link.aip.org/link/?JMAPAQ/12/2099/1}{A. J. Bracken and H. S. Green}
\cite{BrackenGreen71,Green71}
found the first examples of related identities for classical semisimple Lie algebras.
The subject was further developed in subsequent papers by Australian group
(see e.g. \cite{BrienCantCarey77,Gould85}). This type of identities plays an important
role in applications (e.g. \cite{JM95},\cite{MacfarlanePfeiffer99}).
\href{http://dx.doi.org/10.1063/1.529152}{M. Gould, R. Zhang, A. Bracken}
\cite{GouldZhangBracken911}
(see formula 29 page 2300) extended results to $U_q(g)$
for semisimple $g$, see also \cite{GouldZhangBracken912}.
{\href{http://www.springerlink.com/content/h11t67772521q615/
}{H. Ewen, O. Ogievetsky,  and J. Wess
}}
\cite{EwenOgievetskyWess91} section 4 lemma 4.1 contains CH identity for $Fun(GL_{p,q}(2))$.
\href{http://dx.doi.org/10.2977/prims/1195165907}{M. Nazarov and V. Tarasov}
 \cite{NT94} (see also \cite{Molev02} section 4.3 page 37, \cite{M07}) found
 new approach to Bracken-Green type identities via the Yangian
and relation with the Capelli determinant was understood.
\href{http://www.ams.org/distribution/mmj/vol1-1-2001/kirillov.pdf}{A. Kirillov}
\cite{Kir} (see also \cite{R05})
generalized  the CH identities related to $U(gl(n))$.
\href{http://www.ams.org/proc/1998-126-11/S0002-9939-98-04557-2/home.html}
{T. Umeda}
 \cite{Umeda} section 3, page 3174  and  \href{http://dx.doi.org/10.1006/jabr.2001.8824}{M. Itoh}
\cite{Itoh2}
gave  another more direct
approach to CH theorem for semisimple Lie algebras and generalizations.
\href{http://dx.doi.org/10.1080/00927879908826430}{ I. Kantor, I.
Trishin} \cite{KT99} (see also \cite{UM94}) considered the case of
super-matrices. A comprehensive study of Cayley-Hamilton and
related identities was attempted in \GIOPS.
A non-trivial character of these identities (in general) is that
instead of the matrix power $M^k$ one considers
the so-called quantum matrix powers introduced first by
\href{http://dx.doi.org/10.1016/0370-2693(90)90677-X}{J.-M. Maillet}  \cite{Maillet90}.
This quantum powers are important from the point of view of quantum integrable
systems since their traces provide the commuting elements (integrals of motion),
while the traces of usual powers do not commute in general.
In all these papers the CH theorem states the linear dependence
of power (or quantum powers) of matrices where coefficients of linear
dependence are elements from the basic ring. Another version of CH theorem
was proposed by
\href{http://dx.doi.org/10.1016/S0022-4049(97)00039-X}{J. J. Zhang}  \cite{Zh93}
 for the case of quantum group $Fun_q(GL_n)$.
In this paper the coefficients of linear dependence are diagonal matrices,
with different (in general) elements on the diagonal, but they are not elements of the ground ring like in theorems above.

A quasideterminant version of the CH theorem and the Capelli formula has appeared in \cite{GR92} (for matrices with coefficients in an arbitrary ring). $gl_n -$ case was detailed in \cite{GKLLRT94} section 8.6 page 96.
The paper by
\href{http://dx.doi.org/10.1023/A:1027316209695}{O. Ogievetsky, A. Vahlas}
 \cite{OgievetskyVahlas03}
 compares the two formulations of the CH theorem (for quantum and for usual powers).
The first formulation is more useful from the point of view of integrable systems.

There are also works of more ring theoretic spirit.
\href{http://dx.doi.org/10.1016/0021-8693(87)90073-1}{C. Procesi} \cite{P87}
proves that an appropriate version of CH identity in ring \kol~ with a trace is necessary and sufficient condition
for existence of  embedding of \kol~ into matrices $Mat(C)$, where $C$ is commutative ring.
The paper
\href{http://dx.doi.org/10.1090/S0002-9939-97-03868-9}{J. Szigeti} \cite{Sz97}
 discusses another generalization of CH theorem for
rings \kol~ such that $\exists n: \forall x_1,...,x_n: [x_1[x_2...[x_{n-1}, x_n]...]]=0$, his result
states that there exists a polynomials $\chi_n$ such that $\chi_n(A)=0$ for any matrix $A$ over \kol.
(See also \href{http://dx.doi.org/10.1016/S0022-4049(97)00184-9}
{M. Domokos} \cite{Do96} for developments and concrete examples).
In \cite{Sz06} the generalization to the case of matrices with values in $[R,R]$ is discussed.

{\bf The Newton identities.}
The identities for \MMs~ seems to appear first in \cite{CF07}.
 For some other classes of matrices with noncommutative entries they
can be found in
\href{http://arxiv.org/abs/hep-ph/9407124}{(I. M. Gelfand, D. Krob, A. Lascoux,
B. Leclerc, V. S. Retakh and J.-Y. Thibon)}
\cite{GKLLRT94},
\GIOPS,
\href{http://www.ams.org/proc/1998-126-11/S0002-9939-98-04557-2/home.html}{(T. Umeda)}
 \cite{Umeda},
\href{http://www.ingentaconnect.com/content/els/00218693/1998/00000208/00000002}
{(M. Ito)} \cite{Ito},
see also  \href{http://arxiv.org/abs/math.QA/0211288}{(A.  Molev)}
 \cite{Molev02} section 4.1 page 36 and bibliographic notes on  page 43, \cite{M07} section 7.1,
for quantum matrices see {\href{http://arxiv.org/abs/math/0302063}{(M. Domokos, T. H. Lenagan)}}
\cite{DL03}.
The Perelomov-Popov formulas \cite{PP66} can also seen as a form of the Newton
identities.
Applying the Newton identities above to the
example defined by formula \ref{Lax-gl}a page \pageref{Lax-gl}
one can probably derive these results.

{\bf The MacMahon-Wronski identity.}
As we have mention the identity for q-\MMs~has been discovered in
\href{http://arxiv.org/abs/math.QA/0303319}{(S. Garoufalidis, T. Le, D. Zeilberger)} \cite{GLZ},
(see also \cite{EtingofPak06,KonvalinkaPak06,FoataHan06}).
The most natural and simple proof based on Koszul duality has been obtained
in
\href{http://arxiv.org/abs/math/0603169}{(Phung Hai, M. Lorenz)}
\cite{HaiLorenz06}.
For \cite{FRT} quantum group matrices $RTT=TTR$ and some more general algebras
the identity has been already known in some form
\href{http://lanl.arxiv.org/abs/math/9809170}{(A. Isaev, O.
Ogievetsky, P. Pyatov)}
 \cite{GIOPSW} ~page 9, (see also \cite{GIOPS95}-\cite{GIOPS05}).
For some other noncommutative matrices the identity is discussed in
\href{http://www.springerlink.com/content/2l45hauppwnjc6be/}{(T.
Umeda)} \cite{Umeda03}, where similar ideas on the Koszul duality are used.


{\bf The Pl\"{u}cker coordinates.}
 For generic noncommutative matrices quasi-Pl\"{u}cker coordinates were studied by
\href{http://arxiv.org/abs/q-alg/9705026}{I. Gelfand, V. Retakh}
\cite{GR97} section 2.1 page 9,
(\href{http://arxiv.org/abs/math.QA/0208146}{I. Gelfand, S. Gelfand, V. Retakh, R. Wilson}
 \cite{GGRW02} section 4 page 33).
Note that in the commutative case they are not the standard Pl\"{u}cker coordinates,
but their ratios.  Quasi-Pl\"{u}cker coordinates were
further studied in .
\href{http://arxiv.org/abs/math.QA/0602448}{(A. Lauve)}
\cite{Lauve06}.
It would be interesting to clarify relations between Pl\"{u}cker coordinates above and quasi-Pl\"{u}cker ones.

Due to relations \href{http://math.uoregon.edu/~arkadiy/andrei.ps}{(A. Berenstein, A. Zelevinsky)}
\cite{BerensteinZelevinsky93}
for the dual canonical (crystal) basis,
the quantum minors and  relations between them for  \cite{FRT} quantum group matrices $RTT=TTR$
 are widely studied:
\cite{LeclercZelevinsky98,Skoda05,Goodearl05,Rita03,CalderoMarsh05,Scott00,KellyLenaganRigal02}.
But these important results are specific to $RTT=TTR$ matrices, one cannot
expect that they hold true for more general class of \GT~ matrices.
It might be natural to ask whether any of the properties of q-Pl\"{u}cker coordinates survives in
a "twicely" wider class of (q)-\MMs ? (Surely not all of them do).



\section{Matrix (Leningrad) form of the defining relations for \MMs
\label{MatrSect}}

In this section we present the defining relations for \MMs~ in the
matrix (Leningrad) notations as well as  some applications. Such
notations are an almost universally used tool in the theory of
quantum groups and of quantum integrable system. We shall herewith
frame the definition and the main properties of \MMs~ within this
formalism. The main benefit is that some formulas (e.g., the
commutation relations between the generators) will be most compactly
written. Also, we will show how some of our statements can be
conveniently translated and used in this formalism.

At first we shall collect some notions coming from the "Leningrad school"'s approach to these issues.
Then we shall consider the case of \MMs~ and finally give a few
applications.

\subsection{A brief account of matrix (Leningrad) notations}
The notations are briefly discussed in various texts,
let us only mention (L. Faddeev, L. Takhtajan) \cite{FT87},
\href{http://books.google.com/books?id=YQSKPnFzDOEC&printsec=frontcover&dq=A+guide+to+quantum+groups&sig=0wLtF19lwtBDVu0lHFwd7kKXspk}
{(V. Chari, A. Pressley)} \cite{CP94} section 7.1C page 222.
Here we provide definitions and examples.


Let \kol~ be an associative algebra over $\CC$. An $n\times n$ matrix
$A\in Mat_n[\mathcal{K}]$ with entries $A_{ij}$ in a non-commutative algebra
\kol~ can be considered as an element: \bea A=\sum_{ij} A_{ij} \otimes
E_{ij} \in \mathcal{K}\otimes Mat_n, \eea where $E_{ij}$ are  standard "matrix
units", i.e. those matrices whose $(i,j)$-th element is $1$, and all the
others are zero. $Mat_n$ is an associative algebra of $n\times n$
matrices over $\CC$.

One considers the tensor product\footnote{All tensor products are
taken over $\CC$.} $\mathcal{K}\otimes Mat_n \otimes Mat_n$ and introduces the
following notations: \bea
\one A= A \otimes 1 =\sum_{ij} A_{ij}\otimes E_{ij} \otimes 1_{n\times n} \in \mathcal{K} \otimes Mat_n \otimes Mat_n,\\
\two B= 1 \otimes B =\sum_{ij} B_{ij}\otimes 1_{n\times n} \otimes
E_{ij}  \in \mathcal{K} \otimes Mat_n \otimes Mat_n, \eea where $1_{n\times
n}$ is identity matrix of size $n\times n$.
For further use, we notice that  the permutation matrix $P\in
Mat_n \otimes Mat_n$, defined by \bea
 P (e_i\otimes e_j)= e_j\otimes e_i.
\eea can be written as $\displaystyle{P = \sum_{ij} E_{ij}\otimes
E_{ji}},$ and satisfies the properties: \bea
&&P^2=1, \\
&&(A \otimes 1) P =  P (1 \otimes A), ~~~ (1 \otimes A) P =  P (A \otimes 1).
\eea
A crucial observation is that, in general $(A \otimes 1) ( 1
\otimes B)\neq(1 \otimes B) ( A \otimes 1)$,
and all $n^4$ commutators between the elements $A_{ij}$
and $B_{kl}$ are encoded in the expression $[A \otimes 1 , 1 \otimes
B]$.
Indeed,
 \bea
 ~ [A \otimes 1 , 1 \otimes B]= \sum_{ijkl} [ A_{ij}, B_{kl}] \otimes E_{ij}\otimes E_{kl}\in \mathcal{K} \otimes Mat_n \otimes
 Mat_n,\\
(A \otimes 1) ( 1 \otimes B) = \sum_{ijkl} A_{ij}B_{kl} \otimes E_{ij}\otimes E_{kl} \in \mathcal{K} \otimes Mat_n \otimes Mat_n,\\
(1 \otimes B) ( A \otimes 1)   = \sum_{ijkl} B_{kl} A_{ij} \otimes E_{ij}\otimes E_{kl} \in \mathcal{K} \otimes Mat_n \otimes Mat_n.
\eea
\subsubsection{Matrix  (Leningrad) notations in  $2\times 2 $ case }
It is useful  to exemplify  the matrix (Leningrad) notations by $2
\times 2$ examples. Although these notions are  standard, it might be convenient
for the reader to reproduce them here.

$Mat_n\otimes Mat_n$ can
be identified with $Mat_{n^2}$, provided an order of basis elements
$e_i\otimes e_j\in \CC^n\otimes \CC^n$ is chosen.
As it is customary, we order the basis in the tensor product  $\CC^2
\otimes \CC^2$ as $e_1 \otimes e_1, e_1 \otimes e_2, e_2
\otimes e_1, e_2 \otimes e_2$. Thus \bea A \otimes 1=
\left(\begin{array}{cccc}
A_{11}  & 0  & A_{12}  & 0 \\
0 & A_{11}   & 0  & A_{12} \\
A_{21} & 0   & A_{22}  & 0  \\
0 & A_{21}   & 0  & A_{22}
\end{array}\right), ~~~~~
1 \otimes B=
\left(\begin{array}{cccc}
B_{11}  & B_{12}  & 0   & 0 \\
B_{21} & B_{22}   & 0  & 0 \\
0 & 0  & B_{11}  & B_{12}  \\
0 & 0  & B_{21}  & B_{22}
\end{array}\right),
\eea

\bea
(A\otimes 1) (1 \otimes B) = \left(\begin{array}{cccc}
A _{11} B  & A_{12} B \\
A _{21} B  & A_{22} B
\end{array}\right)=
\left(\begin{array}{cccc}
A _{11} B_{11} & A_{11} B_{12}  & A_{12} B_{11} & A_{12} B_{12} \\
A _{11} B_{21} & A_{11} B_{22}  & A_{12} B_{21} & A_{12} B_{22} \\
A _{21} B_{11} & A_{21} B_{12}  & A_{22} B_{11} & A_{22} B_{12} \\
A _{21} B_{21} & A_{21} B_{22}  & A_{22} B_{21} & A_{22} B_{22}
\end{array}\right)
,\\
(1 \otimes B) (A\otimes 1)  = \left(\begin{array}{cccc}
B A _{11}   & B A_{12}  \\
B A _{21}   & B A_{22}
\end{array}\right)=
\left(\begin{array}{cccc}
B_{11} A _{11}  & B_{12} A_{11}   & B_{11} A_{12}  & B_{12} A_{12}  \\
B_{21} A _{11}  & B_{22} A_{11}   & B_{21} A_{12}  & B_{22} A_{12}  \\
B_{11} A _{21}  & B_{12} A_{21}   & B_{11} A_{22}  & B_{12} A_{22}  \\
B_{21} A _{21}  & B_{22} A_{21}   & B_{21} A_{22}  & B_{22} A_{22}
\end{array}\right)
.
\eea

The permutation matrix $\PP=E_{11}\otimes E_{11} +E_{12}\otimes
E_{21} +E_{21}\otimes E_{12} +E_{22}\otimes E_{22} $ reads \bea \PP=
\left(\begin{array}{cccc}
1 & 0  & 0 & 0 \\
0 & 0  & 1 & 0 \\
0 & 1  & 0 & 0 \\
0 & 0  & 0 & 1
\end{array}\right).
\eea

\subsection{Manin's relations in the matrix form \label{PLss}}
Here we present a basic lemma which encodes the definition of \MMs~ in matrix
notations. All the commutation relations between $M_{ij}$ are
encoded in one equation, whose  form  does not depend on the size $n$ of matrices.
The lemma  has been suggested to us by P. Pyatov.

Let M be a $n\times n$ matrix with the elements in an associative algebra \kol over $\CC$.
{\Lem \label{P-lem}
A matrix $M$ is a \MM~ {\bf iff} any of the following equivalent formulas hold:
\bea
[ M\otimes 1, 1 \otimes M ] = \PP [ M\otimes 1, 1 \otimes M ],\\
\frac{(1-\PP)}{2} (M\otimes 1)~ (1 \otimes M) \frac{(1-\PP)}{2}= \frac{(1-\PP)}{2} (M\otimes 1)~ (1 \otimes M), \\
\frac{(1+\PP)}{2} (1 \otimes M) ~ (M\otimes 1) \frac{(1+\PP)}{2}=  (M\otimes 1)~ (1 \otimes M) \frac{(1+\PP)}{2},\\
(1-\PP) (1 \otimes M) ~ (M\otimes 1) (1+\PP)= 0.
\eea
}

It should perhaps be noticed that $(1-\PP)/2$ and $(1+\PP)/2$ are two orthogonal idempotents, namely the antisymmetrizer and the symmetrizer.

Here we allowed ourselves some abuse of notations, denoting by $1$ the identity matrix
in $Mat_n$ (e.g.,  $1\otimes M$) as well as the identity matrix in $Mat_n\otimes Mat_n$
(e.g., $1\pm \PP$).

Letting $\PP\in Mat_n\otimes Mat_n$ be the  permutation matrix $\PP (a\otimes b) =  b \otimes a$,
the formulas above are equalities in the associative algebra $ \mathcal{K} \otimes Mat_n \otimes Mat_n$.
Namely $M\otimes 1$ is a shorthand notation for $\sum_{i,j} M_{ij} \otimes E_{ij}\otimes 1  $ and
$1\otimes M$ is $\sum_{i,j}  M_{ij} \otimes 1 \otimes E_{ij}  $ where $E_{ij}$ are standard "matrix units".

{\bf Proof.}
Let us prove $(1) \Leftrightarrow $ ($M$ is a \MM).
\bea
\label{LHS1111} [ M\otimes 1, 1 \otimes M ] = \sum_{i,j,k,l} [M_{ij}, M_{kl}] \otimes E_{ij} \otimes E_{kl},\\
\PP [ M\otimes 1, 1 \otimes M ] = (\sum_{a,b} E_{ab} \otimes E_{ba})
(\sum_{i,j,k,l} [M_{ij}, M_{kl}] \otimes E_{ij} \otimes E_{kl})= \\ =
(\sum_{a,b} \sum_{i,j,k,l} [M_{ij}, M_{kl}]\otimes  E_{ab} E_{ij} \otimes E_{ba} E_{kl})
=
(\sum_{a,b} \sum_{i,j,k,l} [M_{ij}, M_{kl}] \otimes  E_{aj} \delta_{bi} \otimes E_{bl} \delta_{ak})
= \\ =
( \sum_{i,j,k,l} [M_{ij}, M_{kl}] \otimes  E_{kj} \otimes E_{il} )
= \label{LHS11112}
( \sum_{i,j,k,l} [M_{kj}, M_{il}] \otimes  E_{i j} \otimes E_{kl} )
.
\eea
Thus $[ M\otimes 1, 1 \otimes M ] = \PP  [ M\otimes 1, 1 \otimes M ]$
is equivalent to $[M_{ij}, M_{kl}]= [M_{kj}, M_{il}]$, which is the definition of \MMs, and
$(1) \Leftrightarrow $ ($M$ is a \MM) is proved.

To derive that the properties 1,2,3,4 are equivalent to each other
is trivial use of the properties: $(M\otimes 1)P=P (1\otimes M)$, $(1\otimes M)P=P (M\otimes 1)$.
%
%
\BX


{\Cor  Let $M$ or $M^t$ be a \MM, then:}
\bea \label{MM-sq-zero-fml}
[ M\otimes 1, 1 \otimes M ]^2=0.
\eea
\PRF Let us consider $M$ is a \MM, the other case being similar.
\bea
[ M\otimes 1, 1 \otimes M ] = \PP [ M\otimes 1, 1 \otimes M ], \mbox{ let us  square this equality: } \\
~~  [ M\otimes 1, 1 \otimes M ]^2 = \PP [ M\otimes 1, 1 \otimes M ] \PP [ M\otimes 1, 1 \otimes M ],\\
\PP [ M\otimes 1, 1 \otimes M ] \PP=   [\PP M\otimes 1 \PP , \PP 1 \otimes M  \PP]=
 [1\otimes M, M\otimes 1 ]=- [ M\otimes 1, 1 \otimes M ],\\
\mbox{ So:  } [ M\otimes 1, 1 \otimes M ]^2=-[ M\otimes 1, 1 \otimes M ]^2.
\eea
\BX

We will show below see corollary \ref{cur-cor}, page \pageref{cur-cor}:
{\Prop  Let $M$ be a two-sided invertible \MM, then
\bea
0=\PP  [\one M^{-1} , \two M^{-1}]   [ \one M,  \two M ]
= \PP  [\one M , \two M]  \one M^{-1}  \two M^{-1} [ \one M, \two M ].
\eea
}
It is actually equivalent to the theorem that the inverse to a two-sided invertible  \MM~ is again a \MM.


The following generalization of lemma \ref{P-lem} has been discovered by A.Silantiev:
{\Prop \label{prop_AMMM_AMMMA}
Let $M$ be an $n\times n$ \MM. Then $\forall m$:
\begin{align}
  A_m \one M \two M \cdots \oneM M =A_m \one M \two M \cdots \oneM M A_m,\\ \label{AMMM_AMMMA}
   \one M \two M \cdots \oneM M S_m =S_m \one M \two M \cdots \oneM M S_m,
\end{align}
where we consider the tensor product $\mathcal{K}\otimes Mat_n^{\otimes m}$,
$\oneK M$  is $ \sum_{ij} M_{ij} \otimes 1 \otimes ... \otimes E_{ij} \otimes 1 \otimes ... \otimes 1$,
where $E_{ij}$ stands on $k$-th position.
$ S_m, (A_m) $ is (anti)-symmetrizer in $(\CC^n)^{\otimes m}$. I.e. permutation group $S_m$
naturally acts on $(\CC^n)^{\otimes m}$ and $A_m$ is an image of $\sum_{\sigma\in S_m} (-1)^\sigma \sigma$
and $S_m$ is  an image of $\sum_{\sigma\in S_m} \sigma$ under this action.
}

The proof will be provided in \cite{CFRS}.

\subsection{Matrix
(Leningrad) notations in  Poisson case  \label{PoissonManinSec} }
Recall (see subsection \ref{Pois-mm} page \pageref{Pois-mm})
that by a Poisson-Manin we call a matrix with entries in some Poisson algebra,
such that $\{ M_{ij}, M_{kl} \} = \{ M_{kj}, M_{il} \} $.

Let us also give a Poisson version of lemma \pref{P-lem}:
{\Lem \label{PoisPyat} Matrix $M$ is a  Poisson-\MM~ iff:
\bea
\{ \one M \stackrel{\otimes}{,}  \two M \} = \PP  \{ \one M \stackrel{\otimes}{,}  \two M \}.
\eea
Matrix $M^t$ is a  Poisson-\MM~ iff:
\bea
\{ \one M \stackrel{\otimes}{,}  \two M \} =  \{ \one M \stackrel{\otimes}{,}  \two M \} \PP .
\eea
}

Where we have used matrix (Leningrad) notation $\{ \one A \stackrel{\otimes}{,}  \two B \}$
for the Poisson case, which is defined as follows:
\bea
\{ \one A \stackrel{\otimes}{,}  \two B \}
\stackrel{def}{=}
\{ A_{ij}, B_{kl} \} \otimes E_{ij}\otimes E_{kl} .
\eea


We have proved above (theorem \ref{Inverse-conj} page \pageref{Inverse-conj}) that
inverse to a \MM~ is again a \MM, under the condition that $M$ is two sided invertible.
Let us give a Poisson version of this theorem, its proof demonstrates efficiency of
matrix notations in calculations.

{
{\Th
   Assume that $M$  is invertible   Poisson-Manin matrix
then $M^{-1}$ is again Poisson-Manin matrix.
}



{\bf Proof.} 
Due to the lemma \ref{PoisPyat} page \pageref{PoisPyat} above it is enough to prove that:
\bea
\{ \one M^{-1} \stackrel{\otimes}{,}  \two M^{-1} \} = \PP  \{ \one M^{-1} \stackrel{\otimes}{,}  \two M^{-1} \}.
\eea
This can be achieved by the straightforward calculation:
\bea
\{ \one M^{-1} \stackrel{\otimes}{,}  \two M^{-1} \}
= -\two M^{-1}  \{ \one M^{-1} \stackrel{\otimes}{,} \two M^{-1} \} \two M^{-1}
= \two M^{-1}  \one M^{-1} \{ \one M\stackrel{\otimes}{,} \two M \} \one M^{-1} \two M^{-1}
=
\\ \mbox{ use "Pyatov's lemma" \ref{PoisPyat}  : } ~~ \{ \one M\stackrel{\otimes}{,} \two M \}= \PP \{ \one M\stackrel{\otimes}{,} \two M \}
~~~ \mbox{ so we get: }
\\ =
 \two M^{-1}  \one M^{-1} \PP \{ \one M\stackrel{\otimes}{,} \two M \} \one M^{-1} \two M^{-1}
=
\\ \mbox{ let us use that $\forall$ matrices $A$ : } ~~ \one A \PP = \PP \two A,  \two A \PP = \PP \one A
~~~ \mbox{ so we get: }
\\ \label{fml-ll}
=
\PP  \one M^{-1}  \two M^{-1} \{ \one M\stackrel{\otimes}{,} \two M \} \one M^{-1} \two M^{-1}
=
\\ \mbox{ let us use  } \one A \two B = \two B \one A \nn\\
  \mbox{(It is true for matrices with commutative elements).} \nn \\
  \mbox{(This is the only point in arguments which does not work in the noncommutative case).} \nn \\
 =
\PP  \two M^{-1} \one M^{-1}   \{ \one M\stackrel{\otimes}{,} \two M \} \one M^{-1} \two M^{-1}
 =
 \PP  \{ \one M^{-1} \stackrel{\otimes}{,}  \two M^{-1} \}.
\eea
\BX

Let us mention a curious corollary, which arises if one uses the same arguments
as above in the noncommutative case (not a Poisson one) and also
uses the established fact that $M^{-1}$ is again a \MM.

{\Lem\label{cur-cor} Let $M$ be two-sided invertible \MM, then:
\bea 0=\PP  [\one M^{-1} , \two M^{-1}]   [ \one M,  \two M ] = \PP
[\one M , \two M]  \one M^{-1}  \two M^{-1} [ \one M, \two M ]. \eea
}

\section{Conclusion and open questions
\label{OpenSect}
}

Let us make concluding remarks and mention open problems.

In the present paper we considered a class of matrices with noncommutative entries
and demonstrated that most of the  linear algebra properties can be transferred to this class.
We refer to  \cite{CF07,CFRS, CM},  \cite{CFRy}, \cite{RST, CSS08} for applications and
related issues.

The natural question whether one can extend such theory to the other classes
of noncommutative matrices. We already have positive results
on q-analogs \cite{CFRS} and more general matrices appearing in Manin's
framework, we hope to publish this in some future.
However what seems to be unclear - what can be applications
 of the "noncommutative endomorphisms" of general algebras,
what it can give for studies of Sklyanin algebras or Calabi-Yau algebras \cite{Gi06}?

Manin's framework is quite general, nevertheless it does not seem to cover
 all the examples of noncommutative matrices
which appear in applications and for which some sporadic interesting results
has been already obtained.

First class of examples comes from the theory of quantum integrability
- almost all integrable  systems have Lax matrices, and so their quantum versions
provides matrices with noncommutative entries. We expect that
linear algebra can be developed in all these cases and it will have
important applications \cite{CT06-1}, however it is not so clear how to
attack this problem.

Second: related questions appear in quantum group and Lie algebra theory:
$gl(n)$-related matrices in the non-vector representations \cite{Kir}, \cite{R05};
quantum groups for non-$gl(n)$-case;
twisted current algebras and twisted Yangians (e.g. \cite{MR07});
reflection equation and more general quadratic algebras
\href{http://dx.doi.org/10.1016/0370-2693(91)91566-E}{(L. Freidel and J. M. Maillet)} \cite{FM91}.
Let us mention quite an interesting matrix related to symmetric group $S_n$
appearing in
\href{http://dx.doi.org/10.1088/0305-4470/20/10/009}{(M. Gould)}
\cite{Go87} page 1 and \href{http://arxiv.org/abs/math/0006111}{(P. Biane)}
\cite{Bi00} page 3 from completely different
points of view, \cite{Go87} also contains generalizations to the more general finite groups.

The diversity of interesting noncommutative matrices suggest that it might
natural not to work case by case, but rather
ask a general question: given a matrix with noncommutative entries is it possible to understand whether its
 proper determinant  exists or not ? If yes, how to develop the linear algebra ?

So we see that there is a field for the further research which
concerns generalization of \MMs. Let us also mention some more concrete questions
related to
\MMs~ themselves.

\subsection{Tridiagonal matrices and duality in Toda system}
Let us recall a matrix identity for tridiagonal matrices and rise a question
on its extension to \MMs. The identity is very well-known in integrability theory,
it implies that  Toda classical integrable system has two different Lax representations
one by $n\times n$ matrices and another by $2\times 2$ matrices.
In the language of integrability theory our question is about quantization of this
identity.

Let us consider  commutative variables $x_i$ and $p_i$; consider the following matrices:
\bea
\Lcal_{n\times n} (v) =
\left(\begin{array}{ccccc}
  -p_1 & 1 & \ldots & 0 & v^{-1}e^{x_{1n}} \\
  e^{x_{21}} & -p_2 & \ldots & 0 & 0 \\
  \ldots & \ldots & \ldots & \ldots & \ldots \\
  0 & 0 & \ldots & -p_{n-1} & 1 \\
  v & 0 & \ldots & e^{x_{n,n-1}} & -p_n
\end{array}\right),
\qquad x_{jk}\equiv x_j-x_k.
\eea

\bea
 L_{2\times 2} (u)=
 \left(\begin{array}{cc}
     u+p_n & -e^{x_n} \\
     e^{-x_n} & 0 \end{array}\right)
 \left(\begin{array}{cc}
     u+p_{n-1} & -e^{x_{n-1}} \\
     e^{-x_{n-1}} & 0
   \end{array}\right)...
 \left(\begin{array}{cc}
     u+p_{1} & -e^{x_{1}} \\
     e^{-x_{1}} & 0
   \end{array}\right)
\eea

The following identity is well-known in integrability theory (e.g.
\href{http://arxiv.org/abs/nlin/0009009}{(E.K. Sklyanin)} \cite{Sk00} section 2.5 page 13 formula (46)).
It is also known in the other fields of mathematical physics
\href{http://arxiv.org/abs/0712.0681}{(L. Molinari)}
\cite{Mo07} theorem 1 page 5 (here its    generalization for block
tridiagonal matrices has been found).

{\Lem The following is true :
\bea
(-1)^{n-1}det\bigl(u-\Lcal_{n\times n}(v)\bigr)=v^2 ~det(1-v^{-1} L_{2\times 2}(u)).
\eea
}

Now consider the quantum case. I.e. consider noncommuting variables $\hat p_i$ and $\hat x_i$,
such that they satisfy the following commutation relations: $[\hat p_i, \hat x_j]= \delta_{ij}$,
$[\hat p_i, \hat p_j]= [\hat x_i, \hat x_j]=0$.
And consider matrices $\widehat{\Lcal}_{n\times n}(u)$, $\widehat{L}_{2\times 2}(u)$
defined by the same formulas as above substituting $\hat x_i$ for $x_i$ and respectively
 $\hat p_i$  for $p_i$.

 One can see directly or look at \cite{CF07} section 3.2 page 14, that:
\bea
 e^{-\partial_u} \widehat{ L}_{2\times 2} (u) ~~~~ \mbox{ - is a \MM.}
\eea

{\bf Question.} Is it possible to find an appropriate definition for
$"det\bigl(u-\widehat{\Lcal}_{n\times n}(e^{\partial_u})\bigr)"$ such that:
\bea
(-1)^{n-1}"det\bigl(u-\widehat{\Lcal}_{n\times n}(e^{\partial_u})\bigr)"{=}e^{2\partial_u}~ det^{col} (1-e^{-\partial_u} \widehat{L}_{2\times 2}(u)) ~~~{\bf ? }
\eea

The matrix $u-\widehat{\Lcal}_{n\times n}(e^{\partial_u})$ does not seem to be related to \MMs~
at least in any simple way. So solution to this problem may lead to the  development of linear algebra
for noncommutative matrices beyond Manin's case, which is highly desired.

{\Rem ~} For  readers who are not familiar with Lax matrices
let us make  the following remark (see also \cite{CF07} section 3 page 9).
Lax matrices $L(z)$ should satisfy several properties, one of them is:
the characteristic polynomial $det(\lambda -L(z))= \sum_{ij} H_{ij} z^i \lambda ^j$
should produce all Liouville integrals of motion,
i.e. $\{ H_{ij}, H_{kl} \}=0$ and any integral of motion is a function of $H_{ij}$.
 Lax matrices
exist for majority of integrable systems.
One integrable system may have several Lax matrices.
Matrices  $L_{2\times 2} (u)$,  $\Lcal_{n\times n} (u)$ are such two examples of Lax matrices for
the one system -  Toda classical integrable system.
Moreover in quantum case it is quite natural to look
 for a kind of determinantal formula: $"det(\hat \lambda - L(\hat z))"$
to produce all  quantum integrals of motion: $[\hat H_{ij} , \hat H_{kl}] = 0$, and, possibly, to
satisfy other important properties (see \cite{CT06-1}).
Formula
$det^{col} (1-e^{-\partial_u} \widehat{L}_{2\times 2}(u))$
is such a formula for quantum Toda system.

\subsection{Fredholm type formulas}
In the commutative case the following formulas can be found in \cite{C97} (page 1).
Fredholm's formulas for the solution of an integral equation is a particular case of them
(\cite{C97} section 4, page 10.)

Let $\CC [x_1, x_2,...,x_n] $
 be the algebra of polynomials. Let  ${\bf C} (x_1, x_2,...) $ be
the   operator of multiplication on a polynomial $ C (x_1, x_2,...) $;~~ let
${\bf A}(\partial_{x_1},\partial_{x_2},...) $ be a  polynomial of operators
$\partial_{x_1}=\frac{\partial}{\partial{x_1}},$
$\partial_{x_2}=\frac{\partial}{\partial{x_2}},... $.

{\Conj Consider  $n \times n$ \MM~$M$, assume that $[M_{ij},x_k]=0$ and
 define its action on $\CC [x_1, x_2,...,x_n] $
in a standard way: $x_i \to \sum_j M_{ij} x_j$ and $M(x_{i_1}...x_{i_k})=M(x_{i_1}) ... M(x_{i_k})$, then:
\bea Tr_{\CC [x_1, x_2,..., x_n]} \left ( M   {\bf A} (\partial_{x_1},\partial_{x_2},....) {
\bf C} (x_1, x_2,...)\right)=
\left (Tr_{\CC [x_1, x_2,...,x_n]} M \right)\left <  { A} (\partial_{x_1},\partial_{x_2},....)|
\frac{1}{1 -M}C (x_1, x_2,...)\right >\label{tr_1},
\eea}
Here pairing $ <...|... > $ is defined by
the formula: $ < (\partial_{x_1}) ^{i_1}... (\partial_{x_n})
^{i_n}|(x_1) ^{j_1}... (x_n) ^{j_n} >=\delta_{i_1}^{j_1}...
\delta_{i_n}^{j_n} i_1!... i_n! $.
The trace of the operator $M$: $\CC [x_1, x_2,...] \to$
$\CC [x_1, x_2,...]  $ is  the sum of diagonal elements
in the natural basis $ (x_1) ^{i_1} (x_2) ^{i_2}....(x_k) ^{i_k}$.

Analogously for the anticommuting variables
$\xi_i\xi_j= -\xi_j\xi_i$:
{\Conj ~}
\bea Tr_{\L [\xi_1,\xi_2,...]}\left (M {\bf A} (\partial_{\xi_1},
\partial_{\xi_2},....) {\bf C } (\xi_1,\xi_2,...)\right)=
\left (Tr_{\L [\xi_1,\xi_2,...]}M \right)\left < A (\partial_{
\xi_1},\partial_{\xi_2},....)|\frac{1}{1+M}C (\xi_1,\xi_2,...)
\right>\label{tr_2}\label{first_fml}\eea

{\Cor
Let ${\bf v } = \sum_k v_k x_k , {\bf w} = \sum_k w_k \partial_{x_k} $,
$ C (e^{\bf v}): \CC [x_1, x_2,..., x_n]\to \CC [x_1, x_2,..., x_n] $
operator of multiplication on $(e^{\bf v})$
and $A (e^{\bf w}):\CC [x_1, x_2,..., x_n]\to \CC [x_1, x_2,..., x_n]$
exponential of differentiation operator, then

\bea\frac{Tr_{\CC [x_1, x_2,..., x_n]}\left (M A (e^{\bf w}) C (e^{\bf v})\right)
}{Tr_{\CC [x_1, x_2,..., x_n]}\left (M\right)}
=e^{\left < {\bf w}|\frac{1}{1 -M} {\bf v}\right >}
\label{main_sl_exp}\eea
}
(See \cite{C97} 3.4 page 9).

\subsection{Tensor operations, immanants, Schur functions ...  }
In the commutative case one can consider tensor powers $V^{\mu}$ of $V=\CC^n$
and corresponding tensor powers $M^{\mu}$ of any matrix $M$, indexed
by a Young diagram $\mu$.
Can one extend this to the case of \MMs ?

The obvious problem is the following: symmetric powers of a \MM~$M$ can be defined
by the right action on $x_i$, while antisymmetric powers  by the left action on $\psi_i$.
So it is not clear whether the natural way to mix left and right actions exists or not.

Schur functions $S(\l_1,...,\l_n)$ can be considered as traces $Tr(M^{\mu})$, where
$\l_i$ are eigenvalues of $M$. They satisfy plenty of relations.
Can one generalize them to the case of \MMs ?
$Tr(M^{\mu})$ are written explicitly in terms of sum of principal $\mu$-immanants of $M$
($\sum_{i_1,i_2,...} \sum_{\sigma\in S_n} \mu(\sigma)\prod M_{i_p,\sigma(i_p)}$,
where $\mu(\sigma)$ is a character of irrep of $S_N$ corresponding to Young diagram $\mu$).
The determinant and permanent are particular cases of immanants.
The problem above show up itself again:  the permanent was defined by the row-expansion,
while the determinant via column expansion. One can try  to consider
symmetrized immanants. (For \MMs~ symmetrized-determinant equals to column-determinant
and symmetrized-permanent equals to row-permanent, since determinant (permanent) behaves well
under permutations of columns (rows)).

Progress in this question may be applied to quantum immanants theory
 by A. Okounkov, G. Olshansky \cite{Ok96-1,OkB96,OkOlsh}, (see also \cite{M07})
and to the related  so-called "fusion" procedure in quantum integrable systems theory
(e.g.
\href{http://arxiv.org/abs/hep-th/0703147}{(V. Kazakov, A. Sorin, A. Zabrodin)}
 \cite{KSZ07},\\
\href{http://arxiv.org/abs/math.QA/0508506}{(D. Gurevich, P. Pyatov, P. Saponov)}
 \cite{GIOPS05}) via applications to examples of \MMs~ considered in \cite{CF07}.

Noncommutative Schur functions were
proposed in\\ \href{http://arxiv.org/abs/hep-ph/9407124}{(I. M. Gelfand, D. Krob, A. Lascoux,
B. Leclerc, V. S. Retakh and J.-Y. Thibon)}
\cite{GKLLRT94}. Is this recipe related to symmetrized (or whatever) immanants of \MMs?
Can one relate these Schur functions to quantum immanants from
\href{http://arxiv.org/abs/q-alg/9602028}{(A. Okounkov)}
\cite{Ok96-1}
via specialization to a \MM~ defined by formula   \ref{Lax-gl}a page \pageref{Lax-gl} ?

These questions seems to us quite important, but
we have not analyzed them yet.
One of the reasons is lack of time, experience, etc... while another is that we have not seen
any suggestion for the optimism, yet.

{\Quest ~} Is there something interesting about the co-product of $\esf_k, \csf_k, Tr(M^k)$ ?
(Hopf algebra structures are useful in the symmetric function theory (\cite{GKLLRT94} and references therein),
but co-product natural for \MMs~ is different from  used therein.

\subsection{Other questions}


Let us consider algebra generated by Manin's relations $[M_{ij}, M_{kl}]=[M_{kj}, M_{il}]$,
consider the moduli space of all its $k$-dimensional representations
(i.e. just the set of $k\times k$ matrices $A_{ij}$ satisfying the relations above, modula
conjugation).
It is some manifold (orbifold).  What can be said about it? 
That question might be of some interest since in a particular case it includes
the "commuting variety": $A,B: [A,B]=0$, which is a subject of intensive research.

Concerning the algebra  generated by $M_{ij}$ with the only
relations  $[M_{ij}, M_{kl}]=[M_{kj}, M_{il}]$,
it seems also little to be known.
Are their left or right zero divisors? If no - can one embed it into some field of fractions?
What is its Poincar\`e series with respect to the natural grading ?

Recently some non-linear algebra of multi-index multi-matrices began
to emerge \cite{ADMS}, where appropriate resultants play  role
of the determinant, it might be interesting to obtain some
noncommutative generalizations of their results in the spirit
of the present paper.

One may also try to study some properties of random \MMs.


\end{document}